\documentclass[a4paper,11pt]{article}
\usepackage{amsmath, amssymb, amsfonts, amstext, amsthm, textcomp, enumerate}
\usepackage[mathscr]{euscript}
\usepackage{mathtools}
\usepackage{setspace, geometry}
\newtheorem{thm}{Theorem}[section]
\newtheorem{lem}[thm]{Lemma}
\newtheorem{cor}[thm]{Corollary}
\newtheorem{pro}[thm]{Proposition}
\newtheorem{defn}[thm]{Definition}
\newtheorem{rem}[thm]{Remark}
\newtheorem{ex}[thm]{Example}

\usepackage{color}

\usepackage{amsmath}

\begin{document}
	\begin{center}
		{\bf\large\large A Unified Approach to Stein's Method for Stable Distributions}
	\end{center}
	\begin{center}
	{Neelesh S Upadhye* and Kalyan Barman**}
\end{center}
\begin{center}
	{ Department of Mathematics\\
		Indian Institute of Technology Madras\\
		Chennai-600036, India\\
		*Email: neelesh@iitm.ac.in\\
		**Email: kalyanbarman6991@gmail.com}
\end{center}
\begin{abstract}
\noindent
In this article, we first review the connection between L\'evy processes and infinitely divisible random variables, and the classification of infinitely divisible distributions. Using this connection and the L\'evy-Khinchine representation of the characteristic function, we establish a Stein identity for an infinitely divisible random variable. The classification and slight modification in approach give us a Stein identity for an $\alpha$-stable random variable with $\alpha\in (0,2).$ Using fine regularity estimates for the solution to Stein equation, we derive error bounds for $\alpha$-stable approximations. We then apply these results to obtain rates of convergence. Finally, we compare these rates with the results available in the literature. 
\end{abstract}

\vspace{2mm}


\section{Introduction}
\label{Intro}

\noindent
Over the last five decades, Stein's method has been an important tool for studying approximation problems. This method was first introduced for normal approximation by Charles Stein \cite{k2} in 1972. Thereafter, Chen \cite{k3} applied this method for Poisson approximation. Several extensions of this method for various well-known probability distributions are studied by several authors. For instance, Stein's method for the classical distributions, such as the gamma distribution \cite{k5}, Laplace distribution \cite{k6} and exponential distribution \cite{k4} is well-known. This technique is not limited to the classical distributions but also extends general families of distributions, such as the Pearson family \cite{k9}, variance-gamma \cite{k10} and discrete Gibbs measure family \cite{k22,k21}. For an overview of Stein's method, we refer the reader to \cite{k25}, and for distribution specific development of Stein's method, we refer to the web-page maintained by Yvik Swan: https://sites.google.com/site/steinsmethod/home and references therein.

\noindent
 In general, Stein's method can be divided into three parts. In the first step, one finds an appropriate Stein operator ${\cal A}$ for a given target random variable $X\sim F_{X}$, a given probability distribution such that $\mathbb{E}(\mathcal{A}f(X))=0$, for $f\in\mathcal{F}$ (an appropriate class of functions). In the second step, one solves a Stein equation

\begin{equation}\label{PP1:e1}
\mathcal{A}(f)(x)=h(x)-\mathbb{E}h(X)
\end{equation}

\noindent
with $h\in\mathcal{H}_{X}$ (a class of smooth functions). In the last step, one derives regularity estimates for solution to the Stein equation \eqref{PP1:e1} with respect to $h$, which are used to derive ``Stein factors'', by bounding the quantity $\left|\mathbb{E}(h(Y)-h(X))\right|$, where $Y\sim G_{Y}$ is a random variable of interest. The problem of getting an upper bound of $|\mathbb{E}(h(Y)-h(X))|$ is equivalent to bounding $|\mathbb{E}\mathcal{A}(f)(Y)|$ which gives an error bound in approximation of the distribution of a random variable $Y$ to the distribution of the target random variable $X$.

\noindent
Application of Stein's method for a random variable $X\sim F_{X}$ relies heavily on identifying a suitable operator $\mathcal{A}$ such that $\mathbb{E}(\mathcal{A}f(X))=0$, as $\mathcal{A}$ characterizes the behavior of the distribution $F_{X}$. In recent years, several techniques have been developed for finding Stein operator such as the density approach \cite{k11}, the generator approach \cite{k8}, the probability generating function approach \cite{k7} and so on.\vspace{2mm}

\noindent
Recently, Arras and Houdr\'e developed Stein's method for infinitely divisible distributions (IDD) with finite first moment (see \cite{k1}), and for multivariate self-decomposable distributions (with finite first moment (see \cite{PP1:arras0}), and without finite first moment (see \cite{k15})). An important subclass of these IDD and self-decomposable distributions is the $\alpha$-stable distributions. In this direction, Xu \cite{k14} developed Stein's method for symmetric $\alpha$-stable distributions with $\alpha\in (1,2)$. Jin et. al. \cite{jin} extended Xu's idea \cite{k14} and developed Stein's method for asymmetric $\alpha$-stable distributions with $\alpha\in (1,2)$. Chen et. al. \cite{k20} also developed Stein's method for asymmetric $\alpha$-stable distributions with $\alpha\in (1,2)$. Later, Chen et. al. \cite{multivariate1} extended it for multivariate case and developed Stein's method for multivariate $\alpha$-stable distributions with $\alpha\in(1,2)$. More recently, Chen et. al. \cite{k19} developed Stein's method for $\alpha$-stable distributions with $\alpha\in(0,1].$ An overview of these articles is discussed in more detail in Section \ref{PP1:preli}. \vspace{2mm}

\noindent
It is clear from the existing literature that the techniques for developing Stein's method for $\alpha$-stable distributions depend on ranges of $\alpha$ ($\alpha\in(0,1]$ and $\alpha\in (1,2)$) and are different. It raises an interesting question:

\noindent
\textit{\begin{enumerate}
		\item[] (Question) For $\alpha\in (0,2)$, can one unify the Stein's method for $\alpha$-stable distributions?
\end{enumerate}} 
\noindent
In this article, we obtain a Stein identity for an infinitely divisible random variable using the L\'evy-Khinchine representation of the characteristic function. With slight modification in approach, we also obtain a Stein identity for an $\alpha$-stable random variable with $\alpha\in (0,2).$ We solve our Stein equation in a unified way via the semigroup approach. Using the fine regularity of the solution to our Stein equation, we demonstrate error bounds for $\alpha$-stable approximations. Finally, we apply these results to obtain the convergence rates for $\alpha$-stable approximations using two known examples, and we also compare our rates with the existing literature.

\noindent
The organization of the article is as follows. In Section \ref{PP1:preli}, we discuss some preliminaries and known results. In Section \ref{PP1:mr}, we state our result on Stein identity for an infinitely divisible random variable, and in particular for an $\alpha$-stable random variable. Using regularity estimates of solution to our Stein equation, we compute bounds in appropriate probability metrics for $\alpha\in(0,1]$ and $\alpha\in(1,2)$ respectively. In Section \ref{PP1:application}, we discuss two applications of our results for $\alpha$-stable approximations and obtain the convergence rates. In Section \ref{PP1:proof}, we provide the proofs of results presented in Section \ref{PP1:mr}.

\section{Preliminaries and known results}\label{PP1:preli}
\noindent
In this section, we review the relationship between IDD and L\'{e}vy processes. We also establish the classification of $\alpha$-stable distributions based on L\'evy processes. Further, we discuss the results on convergence rates for $\alpha$-stable approximations.

\subsection{Infinite divisibility and L\'evy process}
\noindent
 Let us first define the concept of infinite divisibility.

\begin{defn} \cite[p.3]{and}\label{PP1:defidd}
	The distribution of a random variable $X$ is said to be infinitely divisible, if, for every $n \in {\mathbb N}$, 
	$$X \stackrel{d}{=} X_{n,1}+ \ldots+X_{n,n},$$
	where $X_{n,1},\ldots,X_{n,n}$ are independent and identically distributed (i.e., $X_{n,j} = X_n$, $j = 1,2,\ldots n$) and $X_n$ is called $n$-th factor of $X$.   
\end{defn}
\noindent
In other words, a distribution function $F_X$ is infinitely divisible if, for each $n \in {\mathbb N}$, $F_X$ is $n$-fold convolution of $F_{n,X_n}$ with itself (i.e., $F_X = F_{n,X_n}^{*n}$), where $F_{n,X_n}$ is $n$-th factor of $F_X$. 
This can also be summarized using a notion of characteristic exponent as follows:
Define $\eta(z) :=\log \phi_X (z)= \log {\mathbb E}(e^{izX})$, $z\in {\mathbb R}$ to be the characteristic exponent of a random variable $X$. Then, the distribution of random variable $X$ ($F_X$) is infinitely divisible, if, for each $n \in {\mathbb N}$, there exist a characteristic exponent $\eta_n(\cdot)$, such that $\eta(z) = n \eta_n(z)$, $z \in {\mathbb R}$.

\noindent
Next, we use this property of characteristic exponents for some familiar distributions and show that these are in fact infinitely divisible.

\noindent
\begin{ex}[Normal distribution \cite{walsh}] Let $X \sim {\cal N}(\beta, \sigma^2)$, where $\beta \in {\mathbb R}$, and $\sigma > 0$. Then it is well-known that the characteristic exponent of $X$ is given by
	\begin{eqnarray*}
		\eta(z) =  i \beta z - \frac{\sigma^2 z^2}{2}
		= n \left(i (\beta/n) z - \frac{(\sigma^2/n) z^2}{2}\right)= n \eta_n(z).
	\end{eqnarray*}
	Observe now that, for every $n \in {\mathbb N}$, $\eta_n(z)$ is the characteristic exponent of the random variable $X_n \sim {\cal N}(\beta/n, \sigma^2/n)$. Hence the distribution of $X$ is infinite divisible.
\end{ex}

\begin{ex}[Poisson Distribution \cite{cont}] Let $N \sim Poisson(\lambda)$, where $\lambda > 0$. Then the characteristic exponent of $N$ is given by
	$$\eta(z) =   \lambda(e^{i z} -1) = n((\lambda/n) (e^{i z} -1)= n \eta_n(z).$$ 
	Note that, for each $n \in {\mathbb N}$, $N_n \sim Poisson(\lambda/n)$. Hence the distribution of $N$ is infinitely divisible.
\end{ex}
\noindent
Next, we recall the stochastic processes associated with these examples, namely, Brownian motion and Poisson process, and explore their connection with infinite divisibility.
\begin{defn}[Brownian Motion \cite{walsh}] \label{defbm} A real-valued stochastic process $\{X_t\}_{t \ge 0}$ on a probability space $(\Omega, {\cal F}, {\mathbb P})$ is said to be a Brownian motion, if
	
	\begin{enumerate}
		\item[(i)] $X_0 = 0$ a.s.
		\item[(ii)] For any fixed $\omega \in \Omega$, $t \mapsto X_t$ is continuous a.s.
		\item[(iii)] For $0 \le s\le t$, $X_t - X_s \stackrel{d}{=} X_{t-s}$.
		\item[(iv)] For any partition of interval $[0,t]$, $0=t_0 < t_1< \cdots< t_n = t$, the increments $X_{t_1} - X_0, X_{t_2} - X_{t_1}, \ldots, X_{t_n} - X_{t_{n-1}}$ are independent.
		\item[(v)] For $t > 0$, $X_t \sim {\cal N}(0, t)$. 
	\end{enumerate}
\end{defn}

\noindent
From (v), it is clear that the process is generated from a standard  normal random variable $X \sim {\cal N}(0,1)$. Also, from (ii), we see that sample paths are continuous and not monotone.  

\begin{defn}[Poisson Process \cite{cont}] \label{defPP}  A non-negative integer-valued stochastic process $\{N_t\}_{t \ge 0}$ on a probability space $(\Omega, {\cal F}, {\mathbb P})$ is said to be a Poisson process, if
	
	\begin{enumerate}
		\item[(i)] $N_0 = 0$ a.s.
		\item[(ii)] For any fixed $\omega \in \Omega$, $t \mapsto N_t$ is right continuous with left limits.
		\item[(iii)] For $0 \le s\le t$, $N_t - N_s \stackrel{d}{=} N_{t-s}$.
		\item[(iv)] For any partition of interval $[0,t]$, $0=t_0 < t_1< \cdots< t_n = t$, the increments $N_{t_1} - N_0, N_{t_2} - N_{t_1}, \ldots, N_{t_n} - N_{t_{n-1}}$ are independent.
		\item[(v)] For $t > 0$, $N_t \sim Poisson(\lambda t)$. 
	\end{enumerate}
\end{defn}

\noindent
From (v), it is clear that the process is generated from $N \sim Poisson (\lambda)$. Also, from (ii), we can see that the sample paths are right continuous and non-decreasing.

\noindent
Observe now that these two processes may appear to be quite different from each other, but the distributions that generate these processes are infinitely divisible. Let us look at the processes closely. We can see that these two processes have common properties, such as right-continuous sample paths, stationary and independent increments (from (iii) and (iv)), and generated from infinitely divisible random variables (from (v)). Let us use these common properties and define a new class of processes known as L\'{e}vy processes.

\begin{defn}[L\'{e}vy Process \cite{and}]  A real-valued stochastic process $\{X_t\}_{t \ge 0}$ on a probability space $(\Omega, {\cal F}, {\mathbb P})$ is said to be a L\'{e}vy process, if
	
	\begin{enumerate}
		\item[(i)] $X_0 = 0$ a.s.
		\item[(ii)] For any fixed $\omega \in \Omega$, $t \mapsto X_t$ is right continuous with left limits.
		\item[(iii)] For $0 \le s\le t$, $X_t - X_s \stackrel{d}{=} X_{t-s}$.
		\item[(iv)] For any partition of interval $[0,t]$, $0=t_0 < t_1< \cdots< t_n = t$, the increments $X_{t_1} - X_0, X_{t_2} - X_{t_1}, \ldots, X_{t_n} - X_{t_{n-1}}$ are independent.
		\item[(v)] For $\varepsilon > 0$, $\lim_{h \to 0}{\mathbb P}(|X_{t+h} - X_t| \ge \varepsilon) = 0$.
	\end{enumerate} 
\end{defn}

\noindent
In short, a stochastic process can be characterized as L\'{e}vy process if its sample paths satisfy (ii) and have stationary and independent increments (from (iii) and (iv), respectively).

\noindent
Next, we focus on the relation between infinite divisibility and L\'{e}vy process. From the definition of L\'{e}vy process, it is clear that the distribution of $X_t$ is infinitely divisible. To see this, observe that 
\begin{equation}
X_t = (X_t - X_{(n-1)h})+ (X_{(n-1)h}- X_{(n-2)h})+ \cdots+ (X_h - X_0),\label{iid}
\end{equation}
where $h = t/n$ and $n \in {\mathbb N}$, and these increments are independent and identically distributed with $X_0 =0$. Hence, from Definition \ref{PP1:defidd}, the distribution of $X_t$ is infinitely divisible.  We can also use characteristic exponent to show that $X_t$ has IDD, for any $t >0$. To see this, let $\eta_t(z) = \log {\mathbb E}(e^{izX_t})$, $z \in {\mathbb R}$. Assume first that $t = m \in {\mathbb N}$ then, from \eqref{iid}, with $h = m/n$, $\eta_m(z) = n \eta_{m/n}(z)$. Similarly, for $t \in {\mathbb Q}_{+}$, set of positive rational numbers, say $t = m/n$, $\eta_t(z) = \eta_{m/n}(z) = (m/n) \eta_1(z)$ follows by choosing $h = n/m$ and \eqref{iid}. Now, for $t$ in positive irrationals, construct a decreasing sequence $\{t_n\}$  of positive rational numbers such that $t_n \rightarrow t$ as $n \rightarrow \infty$, then $\eta_t(z) = \lim_{n \rightarrow \infty} \eta_{t_n}(z) = \lim_{n \rightarrow \infty} t_n \eta_1(z) = t \eta_1(z)$. The last but one equality follows from continuity of sample paths (see, (ii) in the definition of L\'{e}vy process and dominated convergence theorem). Hence, using characteristic exponent, we have proved that $\eta_t (z) = t \eta_1(z)$, $z \in {\mathbb R}$ and $t > 0$. This shows that, for any $t > 0$,  $X_t$ has IDD with characteristic exponent $\eta_t(\cdot)$ and can be generated using the distribution of $X_1$ with characteristic exponent $\eta_1(\cdot)$. The above discussion can now be summarized in the following theorem.
\begin{thm} \cite[p.81]{cont}
	Let $\{X_t\}_{t \ge 0}$ be a real-valued L\'evy process. Then there exists an IDD $F$ such that $X_1 \sim F$.
\end{thm}  

\noindent
In particular, the representation of characteristic exponent of $F$ is given in following formulation. 

\begin{thm}\cite[p.5]{and}\label{t1}
	Let $\{X_t\}_{t\geq 0}$ be a real-valued L\'evy process. Then there exists a triplet $(\beta,\sigma^{2},\nu),$ where $\beta \in \mathbb{R}$, $\sigma \geq 0$ and $\nu$ is a measure concentrated on $\mathbb{R}\setminus \{0\}$ satisfying $\int_{\mathbb{R}}(1 \wedge z^{2})\nu(dz)<\infty,$ such that 
	
	$$\mathbb{E}(e^{izX_t})=e^{t \eta (z)},~~\text{for}~~z \in\mathbb{R},$$	$\text{with}~~\eta(z)=i\beta z-\frac{\sigma^{2}z^{2}}{2}+\displaystyle\int_{\mathbb{R}}\left(e^{iuz}-1-iuz\mathbf{1}_{(|u|\leq 1)}\right)\nu(du).$
	
\end{thm}

\noindent
Note that $\eta(\cdot)$ is the characteristic exponent of $F$ and the measure $\nu(\cdot)$ is known as L\'{e}vy measure (need not be a probability measure).

\noindent
This brings us to the important question. Given an IDD $F$, can we construct a L\'evy process $\{X_t\}_{t \ge 0}$ such that $X_1 \sim F$? To answer this question, in view of Theorem \ref{t1}, we need assure the existence of the triplet $(\beta, \sigma^2, \nu)$ associated with $F$. This can be seen in following theorem.

\begin{thm} \cite[p.3]{and}
	A distribution $F$ with characteristic exponent $\eta (\cdot)$ is infinitely divisible
	if and only if there exists a triplet $(\beta, \sigma^{2},\nu),$ where $\beta \in \mathbb{R}$, $\sigma\geq 0$ and $\nu$, the L\'evy measure on $\mathbb{R}\setminus \{0\}$ satisfying $\int_{\mathbb{R}}(1 \wedge z^{2})\nu(dz)<\infty,$ with
	\begin{equation}
	\eta(z)=i\beta z-\frac{\sigma^{2}z^{2}}{2}+\int_{\mathbb{R}}(e^{iuz}-1-iuz\mathbf{1}_{(|u|\leq 1)})\nu(du). \label{chexp}
	\end{equation}
\end{thm}
\noindent
The proofs of these theorems are quite lengthy and involved, we refer the interested readers to Sato \cite[p.41]{k16} for more detailed discussion. 

\noindent
We have now established the fact that, for any IDD $F$ with triplet $(\beta, \sigma^2, \nu)$, there exist a unique L\'{e}vy process $\{X_t\}_{t \ge 0}$. Let us understand the concept through the following examples.
\begin{enumerate}
	\item [Ex.1.] Let $X \sim {\cal N}(0,1)$. Then $\eta(z) = -z^2/2$. On comparison with (\ref{chexp}), we get the triplet $(\beta, \sigma^2, \nu) = (0,1,0)$ and associated L\'{e}vy process is Brownian motion as defined in Definition \ref{defbm}.
	
	\item [Ex.2.] Let $N \sim Poisson(\lambda)$, $\lambda >0$. Then $\eta(z) = \lambda(e^{iz} - 1)$. On comparison with (\ref{chexp}), we get the triplet $(\beta, \sigma^2, \nu) = (\lambda,0,\lambda \delta_1)$, where $\delta_1$ is the Dirac measure at $\{1\}$, and associated L\'{e}vy process is Poisson process as defined in Definition \ref{defPP}.

	\item [Ex.3.] Let $X \sim Gamma (\lambda, \gamma)$, where $\lambda>0, \gamma > 0$. Then $\eta(z) = - \gamma \log (1 -iz/\lambda)$. On comparison  with (\ref{chexp}), we get the triplet $(\beta, \sigma^2, \nu) = (\gamma(1 - e^{\lambda}),0, \nu_1)$, where $\nu_1(du) = (\gamma e^{-\lambda u}/u) du$. For more details on computation of the triplet, we refer the readers to \cite{and}. The associated L\'{e}vy process is known as gamma process. 
	
	\item [Ex.4.] Let $X\sim Cauchy(x_0, c)$, $x_0 \in {\mathbb R}, c > 0$. Then $\eta(z) = ix_0 z - c |z|$. On comparison with (\ref{chexp}), we get the triplet $(\beta, \sigma^2, \nu) =(x_0 - 2c\Gamma/\pi, 0, \nu_2)$ where $\Gamma = \displaystyle{\int_0^\infty \left(\frac{\sin u}{u^2} - \frac{\mathbf{1}_{\{u:|u|\leq 1\}}(u)}{u}\right) du}$, and $\nu_2(du) = (c/(\pi u^2)) du$. The associated L\'{e}vy process is known as 1-stable process. 
\end{enumerate}

\noindent
In the examples discussed above, we see that the IDD and associated L\'{e}vy process are uniquely characterized by the triplet $(\beta, \sigma^2, \nu)$. Also, the behavior of $\nu$ is different in each of the examples. For Poisson distribution, $\nu({\mathbb R}) = \int_{{\mathbb R}} \nu(du) = \lambda <\infty$, for normal and gamma distribution, $\nu({\mathbb R}) = 0$ and $\nu({\mathbb R}) = \infty$, respectively, and for Cauchy distribution, $\int_{\{u:|u| \le 1\}} u \nu_2(du) = \infty$. Also, observe that  the behavior $\sigma$ is also important for normal distribution. These two components of the triplet need to be further classified. Sato \cite[Definition 11.9]{k16} has classified L\'evy process $\{X_t\}_{t \ge 0}$ (infinitely divisible distribution ($X_1$)) into three different classes based on the triplet $(\beta, \sigma^2, \nu)$, as follows.

\begin{description}
	\item[(Type A)] $\sigma=0$ and $\nu(\mathbb{R})<\infty$. (e.g. Poisson process).
	
	\item[(Type B)] $\sigma=0,\nu(\mathbb{R})=\infty$, $\int_{\{|u|\leq 1\}}u\nu(du)<\infty$. (e.g. gamma process) 
	
	\item[(Type C)] $\int_{\{|u|\leq 1\}}u\nu(du)=\infty$ or $\sigma >0$. (e.g. 1-stable process or Brownian motion) 
\end{description}

\noindent
The examples studied here are by no means exhaustive. The class of IDD is very rich and include various well-known distribution like Student's t-distribution, Pareto distribution, $F$-distribution among many others. 

\noindent
Next, we focus on an important subclass of IDD, namely, non-Gaussian stable distributions. This class is characterized by the triplet $(\beta, \sigma^2, \nu) = (\beta, 0, \nu_{\alpha})$, with $\beta \in {\mathbb R}$ and the L\'{e}vy measure $\nu_{\alpha}$ is given by
\begin{equation}\label{u1}
\nu_{\alpha}(du)=\left(m_1\frac{1}{u^{1+\alpha}}\mathbf{1}_{(0,\infty)}(u)+m_2\frac{1}{|u|^{1+\alpha}}\mathbf{1}_{(-\infty,0)} (u)  \right)du,
\end{equation}
where $\alpha \in (0,2)$, $m_1,m_2 \in [0,\infty]$ and $m_1+m_2 >0$. (see \cite[p.32]{k27}). Next, we give characteristic exponent representation for non-Gaussian stable distribution.

\begin{defn} \cite[p.168]{k17} A real-valued random variable $X$ is said to have non-Gaussian stable (also called $\alpha$-stable) distribution, if there exists a triplet $(\beta,0,\nu_{\alpha}),$ such that for all $z \in \mathbb{R},$ the characteristic exponent is given by
	
	\begin{equation}\label{u3}
	\eta_{\alpha}(z)=\log \phi_{\alpha}(z)=iz\beta+ \displaystyle\int_{\mathbb{R}}\left(e^{izu} -1-izu\mathbf{1}_{\{|u|\leq 1\}}(u)  \right)  \nu_{\alpha}(du),
	\end{equation}
	\noindent
	where $\beta\in \mathbb{R},$ $\alpha \in (0,2),$ and  $\nu_{\alpha}$ is the L\'evy measure defined in (\ref{u1}), and is denoted by $X \sim \mathcal{S}(\alpha, \beta,m_1,m_2).$
\end{defn}

\noindent
Note here that $\beta \in {\mathbb R}$ is the location parameter and $\alpha \in (0,2)$ is stability parameter, useful in determining the decay of the tail of distribution of $X$.

\noindent
Observe next that, based on the classification of IDD summarized earlier in this section, we can classify $\alpha$-stable distributions as follows:

\begin{description}
	\item[(Type B)] $\alpha \in (0,1)$ (as $\nu_\alpha({\mathbb R}) = \infty$, but $\int_{\{|u| \le 1\}} u \nu_{\alpha}(du) <\infty$).
	\item[(Type C)] $\alpha \in [1,2)$. (As, $\int_{\{|u| \le 1\}} u \nu_{\alpha}(du) =\infty)$ .
\end{description} 

\noindent
Observe now that, for stable distributions of Type B ($\alpha \in (0,1)$), as  $\int_{\{|u| \le 1\}} u \nu_{\alpha}(du) <\infty$. The characteristic exponent in (\ref{u3}) can be rewritten as
\begin{align}\label{PP1:u4}
\eta_{\alpha}(z)= iz\beta_{1}+ \displaystyle\int_{\mathbb{R}}\left(e^{izu} -1  \right)  \nu_{\alpha}(du),
\end{align}
\noindent
where $\beta_1=\beta - \int_{\{|u|\leq 1\}}u \nu_{\alpha}(du).$ 

\noindent
Also, for stable distributions of Type C ($\alpha \in (1,2), \alpha \neq 1$), as $\int_{\{|u|\leq 1\}}u\nu_{\alpha}(du)=\infty,$ but $\int_{\{|u|> 1\}}u\nu_{\alpha}(du)<\infty.$ The characteristic exponent (\ref{u3}) can be rewritten as
\begin{align}\label{PP1:u5}
\eta_{\alpha}(z)=iz\beta_2+ \displaystyle\int_{\mathbb{R}}\left(e^{izu} -1-izu  \right)  \nu_{\alpha}(du),
\end{align}

\noindent
where $\beta_2=\beta+\int_{\{|u|> 1\}}u\nu_{\alpha}(du).$

\noindent
Next, we show the connection between our representation and the various other representations of characteristic exponents available for $\alpha$-stable distributions. Based on the well-known four parameter representation of $\alpha$-stable distributions (see, \cite{k26}), where the parameters $\alpha\in (0,2)$, $\gamma_{\alpha}\in\mathbb{R}$, $d_{\alpha}\geq0$ and $\theta \in[-1,1]$ denote the stability, shift, scale and skewness parameters, respectively. Note here that, on careful adjustments of the integrals in \eqref{u3} with respect to $\nu_\alpha$, one can obtain a well-known form of characteristic exponent (characteristic function) of $\alpha$-stable distributions of both types from the L\'{e}vy-Khinchine representation (\ref{u3}) (see, \cite{k26}). The explicit forms are given below:

\begin{description}
	\item [(Type B)] Let $\theta =\frac{m_1-m_2}{m_1+m_2}$, $\gamma_{\alpha} = \beta - \frac{(m_1-m_2)}{(1-\alpha)}$, and $d_{\alpha } = (m_1+m_2)\cos \frac{\pi}{2}\alpha \displaystyle\int_{0}^{\infty}\left(1-e^{-u}\right)\frac{du}{u^{1+\alpha}}$. Then 
	$$\eta_\alpha(z) = iz \gamma_{\alpha}-d_{\alpha}|z|^{\alpha}\left(1-i\theta\frac{z}{|z|}\tan\frac{\pi}{2}\alpha\right).$$
	
	\item[(Type C)] Here we classify further into two cases $\alpha = 1$ and $\alpha \in (1,2)$
	\begin{itemize}
		\item[($\alpha =1$).] Let $\theta =\frac{m_1-m_2}{m_1+m_2}$, $\gamma_\alpha=\beta + (m_1+m_2)\displaystyle\int_{0}^{\infty}\left(\frac{\sin u}{u^{2}}  -\frac{\mathbf{1}_{\{|u|\leq1 \}}(u)}{u}   \right)du$, and $d_\alpha = (m_1+m_2)\frac{\pi}{2}$. Then
		$$\eta_\alpha(z) = iz\gamma_{\alpha}-d_{\alpha}|z|(1+i\theta\frac{z}{|z|}\frac{2}{\pi}\text{log}|z|).$$
		\item[($\alpha \in (1,2)$).] Let $\theta =\frac{m_1-m_2}{m_1+m_2}$, $\gamma_{\alpha} = \beta - \frac{(m_1-m_2)}{(1-\alpha)}$, and $d_{\alpha } = (m_1+m_2)\cos \frac{\pi}{2}\alpha \displaystyle\int_{0}^{\infty}\left(1-e^{-u}\right)\frac{du}{u^{1+\alpha}}$. Then 
		$$\eta_\alpha(z) = iz \gamma_{\alpha}-d_{\alpha}|z|^{\alpha}\left(1-i\theta\frac{z}{|z|}\tan\frac{\pi}{2}\alpha\right).$$
	\end{itemize}
\end{description}


\noindent
 The derivation of these forms of characteristic exponents is given in the Appendix \ref{appendix}. Observe also that, for $X\sim \mathcal{S}(\alpha,0,m,m)$, characteristic exponent of a symmetric $\alpha$-stable random variable is given by 
\begin{align*}
\eta_{\alpha}^{s}(z)=\displaystyle\int_{\mathbb{R}}(e^{izu}-1-izu\mathbf{1}_{\{|u|\leq1\}} (u))\nu_{\alpha}(du)= -d_\alpha|z|^{\alpha} ,~~z\in\mathbb{R},  
\end{align*}
\noindent
where $d_\alpha>0$ is the scale parameter given by

\begin{align*}
d_\alpha=
\begin{cases}
2m\cos(\frac{\pi}{2}\alpha)\displaystyle\int_{0}^{\infty}(1-e^{-y})\frac{dy}{y^{1+\alpha}},~~\alpha\in (0,1),\\
2m\cos(\frac{\pi}{2}\alpha)\displaystyle\int_{0}^{\infty}(1-y-e^{-y})\frac{dy}{y^{\alpha}},~~\alpha\in(1,2),\\
\pi m ,~~\alpha=1.
\end{cases}
\end{align*}

\noindent
For brevity, if we set the scale parameter $d_\alpha=1$. Then, the characteristic exponent $\eta^{s}_{\alpha}$ simplifies to
\begin{align*}
\eta^{s}_{\alpha}(z)= -|z|^{\alpha},~~z\in\mathbb{R}.
\end{align*}
\noindent
From the above discussed explicit form of $\eta_{\alpha}$, it is clear that characteristic function of an $\alpha$-stable random variable is differentiable for $\alpha\in(0,1)\cup (1,2),$ but fails at $\alpha=1$. To fix this problem, we consider a tempered $1$-stable random variable $Y_\gamma$ (see, \cite{tappe}) with characteristic exponent $\eta_{1,\gamma}$, given by
\begin{equation}\label{PP1:cetemst}
\eta_{1,\gamma}(z)= iz\beta+\int_{\mathbb R}(e^{izu}-1-izu\mathbf{1} _{\{|u|\leq1\}}(u))\nu_{1,\gamma}(du),~~z \in \mathbb{R},
\end{equation}
\noindent
where tempering parameter $\gamma\in(0,\infty),$ and $\nu_{1,\gamma}$ is the L\'evy measure defined as

$$	\nu_{1,\gamma}(du):=\left(m_{1}\frac{e^{-\gamma u}}{u^{2}}\mathbf{1}_{(0,\infty)}(u)+m_{2}\frac{e^{-\gamma |u|}}{|u|^{2}}\mathbf{1}_{(-\infty,0)}(u)\right) du.
$$

\noindent
Note that $Y_\gamma$ is infinitely divisible and its characteristic exponent $\eta_{1,\gamma}$ is differentiable on $\mathbb{R}$. Also, it can be easily shown that as $\gamma\to 0^{+}$, $\eta_{1,\gamma}\to \eta_{1}$, the characteristic exponent of $1$-stable random variable $X$. This fact is used later to derive a Stein identity for 1-stable random variable $X$.\vspace{2mm}

\subsection{Probability metrics}
\noindent
Here, we review two well-known probability metrics used in this article.
\begin{enumerate}
	\item {\bf The Wasserstein-$\delta$ distance (\cite{k19})}. Let $$\mathcal{H}_{\delta} = \left\{h:\mathbb{R} \to (\mathbb{R}, d_{\delta}) \bigg| |h^{(k)}(v)-h^{(k)}(x)|\leq d_{\delta}(v,x),k=0,1\right\},$$
	\noindent
	where $d_{\delta}(v,x):=|v-x|\wedge |v-x|^{\delta}$, $h^{(1)}$ is the first derivative of $h$, with $h^{(0)}=h$ and the range of $h^{(k)}$ is endowed with metric $d_\delta$. Then, for any two random variables $V$ and $X$ the distance is given by
	\begin{equation}
\nonumber	d_{W_{\delta}}(V,X):=\sup_{h \in \mathcal{H}_\delta}|\mathbb{E}[h(V)]-\mathbb{E}[h(X)]|, ~~\delta<\alpha\leq1,
	\end{equation}
	This distance is useful for studying stable approximations of Type B ($\alpha \in (0,1)$) and Type C, 
	
	Case 1 ($\alpha=1$). In \cite{k19}, authors use $d_{W_{\delta}}$ distance with $\delta\in (0,\alpha)$, for obtaining non-integrable stable approximations (Type B, and Type C, Case 1). We also see that Chen et. al. \cite[Subsection 1.2]{k19} only consider the test function $h$ to be bounded Lipschitz, and endowed with the metric $d_{\delta}$. They also have discussed relation with other metrics in the existing literature. This is discussed in more detail in Section \ref{PP1:mr} and Section \ref{PP1:application}, respectively.   
	
	\item {\bf The Wasserstein-type distance (\cite{k1})}. Let $$\mathcal{H}_r = \left\{h:\mathbb{R}\to \mathbb{R}\bigg|h \mbox{ is $r$ times continuously differentiable and},\|h^{(k)}\|\leq 1, k =0, 1,\ldots, r \right\},$$ where $h^{(k)}$, $k=1, \ldots,r$, is the $k$-th derivative of $h$, with $h^{(0)}=h$ and $\|f\|=\sup_{x\in\mathbb{R}}|f(x)|$. Then, for any two random variables $Y$ and $Z$ the distance is given by
	\begin{equation*}
	d_{W_r}(Y,Z):=\sup_{h \in \mathcal{H}_r}|\mathbb{E}[h(Y)]-\mathbb{E}[h(Z)]|.
	\end{equation*}
	\noindent
   This distance is useful for studying stable approximations of Type C, Case 2 ($\alpha\in(1,2)$), see \cite{k1}.
   \noindent
   	Johnson and Samworth \cite{john} use Mallows $r$-distance $d_r$, $r>0$ for $\alpha$-stable approximations. Note that $d_r$ is the classical Wasserstein $r$-distance for $r\geq 1$. Moreover, these distances have the following order relationship.
   	
   	$$d_{W_{r}}(Y,Z)\leq d_{W_{1}}(Y,Z)\leq d_{1}(Y,Z)\leq d_{p}(Y,Z),~~r,p\geq 1.$$ 
   	\noindent
   	In \cite[ Subsection 2.3]{k1}, authors use $d_{W_2}$ distance for obtaining approximations to infinitely divisible distributions with first finite moment and they also have discussed relation with other metrics in the existing literature. We also use these relationships and discuss the consequences in Section \ref{PP1:mr} and Section \ref{PP1:application}, respectively.
  \end{enumerate}

\subsection{Literature review}

\noindent
Here, we review the known results on convergence rates for $\alpha$-stable approximations motivated by generalized CLT. The generalized CLT states that the sum of i.i.d. random variables when scaled and centered appropriately, converges to an $\alpha$-stable distribution. Mathematically, assume $(Y_n)_{n\geq 1}$ is a sequence of i.i.d. random variables. For any $n\in\mathbb{N}$, define $S_n=a_{n}(Y_1+Y_2+\ldots+Y_n)-b_{n}$, where $a_{n}\in (0,\infty)$ and $b_{n}\in\mathbb{R}$ are constants. Then,  the sum $S_{n}$ converges weakly to an $\alpha$-stable distribution with stability index $\alpha\in (0,2],$ see \cite{k17}. For $\alpha=2$, if we assume $\mathbb{E}Y_1=0$, $\mathbb{E}Y_{1}^{2}=1$, $\mathbb{E}|Y_{1}|^{3}< \infty$ and set $a_{n}=\frac{1}{\sqrt{n}}$, $b_{n}=0$, then the sum $S_n$ converges weakly to the standard normal distribution $F_{N(0,1)}$, and by Berry-Esseen theorem, it is shown that
\begin{align*}
\sup_{x\in\mathbb{R}}\left|P(S_n \leq x)-F_{N(0,1)}(x)   \right|&\leq \frac{C}{\sqrt{n}}\mathbb{E}|Y_{1}|^{3},
\end{align*} 
\noindent
where $C>0$ is some constant, see \cite{nourdin}. \vspace{1mm}

\noindent
The problem of convergence rates for  $\alpha$-stable approximations in the generalized CLT is studied by many authors using various approaches, see \cite{boon,k18,k23,kuske} for more details. However, we consider an interesting article written by Johnson and Samworth \cite{john} for comparing the results in this article, as it closely matches with our framework. In \cite{john}, authors consider the generalized CLT, and develop the convergence rate for $\alpha$-stable approximation in the Mallows $r$-distance $d_{r},$ where $r>0,$ using the following framework. Let $(Y_n)_{n\geq1}$ be a sequence of i.i.d. random variables with distribution function $F$ such that $F(y)=\frac{c_1+b(y)}{|y|^{\alpha}}$ for $y<0$ and $1-F(y)=\frac{c_2+b(y)}{|y|^{\alpha}}$ for $y>0$, where $c_1,c_2>0$ and $b(y)=O(|y|^{-d}),$ $d>0$. In \cite[Theorem 1.2]{john}, it is shown that the partial sum $S_n=n^{-\frac{1}{\alpha}}(Y_1+Y_2+\ldots+Y_n)$ converges weakly to an $\alpha$-stable distribution with the rate $n^{\frac{1}{\beta}-\frac{1}{\alpha}}$ in the $d_\beta$ distance, where $\beta\in (\alpha,2]$. This result is proved using the coupling technique and Lindeberg method. More precisely, for $(Z_n)_{n\geq1}$ an i.i.d. sequence of $\alpha$-stable distributed random variables, $n^{-\frac{1}{\alpha}}(Z_1+Z_2+\ldots+Z_n)$ has an $\alpha$-stable distribution, and hence
 
\begin{align*}
d_{\beta}^{\beta}\left(S_n, n^{-\frac{1}{\alpha}}\sum_{i=1}^{n}Z_i\right)&=n^{-\frac{\beta}{\alpha}}d_{\beta}^{\beta}\left(\sum_{i=1}^{n}Y_{i},\sum_{i=1}^{n}Z_{i}  \right)=n^{-\frac{\beta}{\alpha}}d_{\beta}^{\beta}\left(\sum_{i=1}^{n}Y_{i}^{*},\sum_{i=1}^{n}Z_{i}^{*} \right)\\
&\leq n^{-\frac{\beta}{\alpha}} \mathbb{E} \left| \sum_{ i=1}^{n}\left(Y_{i}^{*}-Z_{i}^{*}  \right)  \right|^{\beta},
\end{align*}
\noindent
where $(Y_{i}^{*},Z_{i}^{*})$ is a coupling of $Y_i$ and $Z_i$ with the property $\mathbb{E}|Y_{i}^{*}-Z_{i}^{*}|=d_{\beta}^{\beta}(Y_{i},Z_{i})$, and $\{ (Y_{i}^{*},Z_{i}^{*})\}_{1\leq i \leq n}$ are independent, see \cite[Eq.(3)]{john}. Finally, using Esseen inequality \cite[Eq.(11) and Eq.(12)]{john} and by \cite[Lemma 5.1]{john}, authors show that $d_{\beta}\left(S_n,n^{-\frac{1}{\alpha}}\sum_{i=1}^{n}Z_i\right)=O(n^{\frac{1}{\beta}-\frac{1}{\alpha}})$.  \vspace{1mm}

\noindent
As mentioned in Section \ref{Intro}, the goal of this article is to unify Stein's method for $\alpha$-stable distributions. Let us now review the articles in the literature related to $\alpha$-stable distribution approximation via Stein's method. Recently, Arras and Houdr\'e \cite{k1} develop Stein's method for IDD with finite first moment. Let us briefly discuss their ideas in \cite{k1}. In \cite[Theorem 3.1]{k1}, authors obtain a Stein characterization for IDD with finite first moment using covariance representation of functions of infinitely divisible random variables. We see that Stein characterizations of several distributions are followed by \cite[Theorem 3.1]{k1}. Using the Fourier method, they provide general upper bounds for $d_{\text{Kol}}(X_n, X)$, where both $X_n$ and $X$ are infinitely divisible random variables. In \cite[Proposition 5.1]{k1}, they also obtain generators for the self-decomposable distributions and prove a bound for self-decomposable distribution approximation in the Wasserstein-type distance $d_{W_{2}}$, see \cite[Section 6]{k1}. Applying \cite[Theorem 6.1 and Theorem 6.2]{k1}, they demonstrate bounds for $\alpha$-stable approximations with $\alpha\in (1,2)$ in the $d_{W_2}$ distance. \vspace{1mm}

\noindent
In this direction, Xu \cite{k14} develops Stein's method for symmetric $\alpha$-stable distributions with $\alpha\in (1,2)$ and prove a rate of convergence $n^{- \frac{2-\alpha}{\alpha}}$ for $\alpha$-stable approximations in the Wasserstein distance $d_W$. Let us briefly discuss ideas in \cite{k14}. Let $S_{n}=Z_1+Z_2+\ldots+Z_n$ be a sum of $n$ centered i.i.d random variables. By \cite[Theorem 1.4]{k14}, a Stein operator for symmetric $\alpha$-stable random variable is given by $$\mathcal{A}f(x)=\Delta^{\frac{\alpha}{2}}f(x)-\frac{1}{\alpha}xf^{\prime}(x),$$ 
where $\Delta^{\frac{\alpha}{2}}$ is the fractional Laplacian and $f\in\mathcal{F}$ (a class of functions $f$ with first and second derivatives bounded by a constant depending on $\alpha$ and that $\Delta^{\frac{\alpha}{2}}f$ is $\gamma$-H$\ddot{\text{o}}$lder continuous for any $0<\gamma<1$).

\noindent
Using $K$-function approach \cite{PP1:chen0}, Xu shows that,
\begin{align}\label{PP1:in0}
\mathbb{E}\left(S_{n}f^{\prime}(S_{n})  \right)=\sum_{i=1}^{n}\displaystyle\int_{-N}^{N}\mathbb{E}\left( K_{i}(t,N)f^{\prime\prime}(S_{n}(i)+t) \right)dt+R,
\end{align}

\noindent
where $N>0$ is an arbitrary number, $S_{n}(i)=S_{n}-Z_{i}$, $K_{i}(t,N)=\mathbb{E}\left(Z_{i}\mathbf{1}_{\{0\leq t\leq Z_{i}\leq N\}} -Z_{i}\mathbf{1}_{\{-N\leq Z_{i}\leq t\leq0\}}  \right),$ and $R$ is a remainder.

\noindent
Due to the heavy tail property of $Z_i$, Xu  also shows that,

\begin{align}\label{PP1:in1}
\Delta^{\frac{\alpha}{2}}f(S_n)=\displaystyle\int_{-N}^{N}\mathcal{K}_{\alpha}(t,N)f^{\prime\prime}(S_n+t)dt+R^{\prime},
\end{align}

\noindent
where $\mathcal{K}_{\alpha}(t,N)=\frac{m_{\alpha}}{\alpha(\alpha-1)}(|t|^{1-\alpha}-N^{1-\alpha})$ with $m_\alpha=\left(\displaystyle\int_{{\mathbb R}}\frac{1-\cos y}{|y|^{1+\alpha}} dy \right)^{-1}.$ 

\noindent
Using \eqref{PP1:in0} and \eqref{PP1:in1}, it can be shown that
\begin{align}\label{PP1:disbd1}
\mathbb{E}\left(\mathcal{A}f(S_{n})\right)=\sum_{i=1}^{n}\displaystyle\int_{-N}^{N}\mathbb{E} \left(\frac{\mathcal{K}_{\alpha}(t,N)}{n}- \frac{K_{i}(t,N)}{\alpha} \right)f^{\prime\prime}(S_{n}(i)+t)dt+R^{\prime\prime},
\end{align}

\noindent
where $R^{\prime\prime}$ is an another remainder. Hence,  
$$\left|\mathbb{E}\mathcal{A}f(S_{n})\right|\leq\left( \sum_{i=1}^{n}\displaystyle\int_{-N}^{N}\mathbb{E}\left| \frac{\mathcal{K}_{\alpha}(t,N)}{n}- \frac{K_{i}(t,N)}{\alpha} \right|dt\right)\| f^{\prime\prime}\|+\|R^{\prime\prime}\|,$$
\noindent
where $\|f^{\prime\prime} \|=\sup_{x\in \mathbb{{R}}}|f^{\prime\prime}(x)|.$ Therefore, to obtain a rate of convergence, it is sufficient to bound $\|f^{\prime\prime}\|$ and the remainder $\|R^{\prime\prime}\|$.\vspace{1mm}

\noindent
Jin et. al. \cite{jin} extend Xu's idea \cite{k14} and develop Stein's method for asymmetric $\alpha$-stable distributions with $\alpha\in (1,2)$. They also obtain kernel discrepancy type bound as \eqref{PP1:disbd1} and demonstrate a rate of convergence for $\alpha$-stable approximation, see \cite[Theorem 1.1]{jin}. Chen et. al. \cite{k20} also develop Stein's method for asymmetric $\alpha$-stable distributions with $\alpha\in (1,2)$. Using leave-one-out approach developed by Stein \cite{k2}, they derive a rate of convergence for $\alpha$-stable approximations (see, \cite[Theorem 1.4]{k20}). Later, Chen et. al. \cite{multivariate1} extend it for multivariate case and develop Stein's method for multivariate $\alpha$-stable distributions with $\alpha\in(1,2)$. In this sequel, Arras and Houdr\'e develop Stein's method for multivariate self-decomposable distributions with finite first moment \cite{PP1:arras0} and without finite first moment \cite{k15}. \vspace{1mm}

\noindent
More recently, Chen et. al. \cite{k19} develop Stein's method for $\alpha$-stable distributions with $\alpha\in(0,1].$ In \cite[Theorem 4 and Section 4]{k19}, they compute the rate of convergence in the generalized CLT for the partial sum of i.i.d. random variables in the domain of normal attraction of an $\alpha$-stable distribution. Due to lack of finite first moment, the strategy in deriving rate of convergence for $\alpha$-stable approximation differs significantly from the case $\alpha\in (1,2),$ obtained in \cite{k1,jin,k14}. Let us briefly discuss their ideas in \cite{k19}. Recall that $d_{W_{\delta}}(V,X):=\sup_{h \in \mathcal{H}_\delta}\left|\mathbb{E}h(V)-\mathbb{E}h(X)  \right|,\delta\in (0,\alpha)$, where $\mathcal{H}_{\delta}$ is the class of Lipschitz functions $h:\mathbb{{R}}\to (\mathbb{{R}},d_{\delta})$ such that $|h(x)-h(y)|\leq d_{\delta}(x,y)$ and the range $\mathbb{{R}}$ of the function $h$ is endowed with the metric $d_{\delta}(x,y)=|x-y|\wedge |x-y|^{\delta}$. Such a probability metric is suitable for $\alpha$-stable approximation with $\alpha\in (0,1]$, see \cite[Subsection 1.2]{k19}. Let $S_{n}=\frac{n^{-1/\alpha}}{\sigma_{\alpha}}\left( Y_{1}+Y_{2}+\ldots+Y_{n}\right)$, $\sigma_{\alpha}>0$ be a partial sum of i.i.d. random variables in the domain of normal attraction of an $\alpha$-stable distribution. For any $\alpha$-stable random variable $X$, Chen et. al. \cite{k19} show that,

\begin{align}\label{PP1:in2}
d_{W_{\delta}}(S_{n},X) \leq \sup_{h \in \mathcal{H}_\delta} \left|\mathbb{E}\left(\mathcal{L}^{\alpha,\theta}f(S_{n}) \right) -\frac{1}{\alpha}\mathbb{E}\left(S_{n}f^{\prime}(S_{n}) \right)  \right|,
\end{align}
\noindent
where $\theta \in [-1,1]$ is a skewness parameter, $\mathcal{L}^{\alpha,\theta}$ is a generator of an $\alpha$-stable L\'evy process, and $\mathcal{F}_{\alpha,\theta}$ is a class of smooth functions. Due to the lack of first finite moment, one has to compute carefully contribution from the term $\mathbb{E}\mathcal{L}^{\alpha,\theta}f(S_{n})$ and that from the term $\mathbb{E}S_{n}f^{\prime}(S_{n})$ for obtaining an upper bound in $d_{W_\delta}$ distance. In \cite[Section 4]{k19}, we see that the Taylor-like expansion and tail properties of $Y_1$ are the keys for obtaining the upper bound in \eqref{PP1:in2}. Also, the regularity estimates of the solution to Stein equation help them to obtain the rates of convergence for $\alpha$-stable approximations. 

\section{Main results}\label{PP1:mr}
\noindent 
In this section, we discuss the three important components of Stein's method for IDD, as mentioned in the Introduction.  
First, we start with Stein identity for infinitely divisible random variables. Let $\mathcal{S}(\mathbb{R})$ be the Schwartz space defined by
$$\mathcal{S}(\mathbb{R}):=\left\{f\in C^\infty(\mathbb{R}): \lim_{|x|\rightarrow \infty} |x^m\frac{d^{n}}{dx^{n}}f(x)|=0, \text{ for all } m,n\in \mathbb{N}\right\},$$
\noindent
where $C^\infty(\mathbb{R})$ is the class of infinitely differentiable functions on $\mathbb{R}$. It is important to note that the Fourier transform on $\mathcal{S}(\mathbb{R})$ is automorphism onto itself. This enables us to identify the elements of dual space $\mathcal{S}^{*}(\mathbb{R})$ with $\mathcal{S}(\mathbb{R}).$ In particular, if $f \in \mathcal{S}(\mathbb{R})$, and $\widehat{f}(u)=\int_{\mathbb{R}}e^{-iux}f(x)dx,\text{ }u\in\mathbb{R},$ then $\widehat{f}(u)\in \mathcal{S}(\mathbb{R}).$ Similarly, if $\widehat{f}(u)\in \mathcal{S}(\mathbb{R})$, and $f(x)=\int_{\mathbb{R}}e^{iux}\widehat{f}(u)du,\text{ }x\in\mathbb{R},$
then $f(x)\in \mathcal{S}(\mathbb{R})$, see \cite{stein}.

\noindent
Now, we state our first result on Stein identity for infinitely divisible random variables. 
\begin{thm} \label{PP1:Th1}
	Let $X\sim IDD(\beta,\sigma^{2},\nu)$ with characteristic exponent \eqref{chexp} and the characteristic function $\phi_{X}$ be differentiable.

	Then,
	\begin{equation}\label{PP1:stidenidd}
	\mathbb{E}\left((X-\beta)g(X)-\sigma^{2}g^{\prime}(X)-\displaystyle\int_{\mathbb R}u(g(X+u)-g(X)\mathbf{1}_{\{|u|\leq1\}}(u))\nu(du)\right)=0,~~g\in\mathcal{S}(\mathbb{R}).
	\end{equation}
\end{thm}
\begin{rem}
	\begin{enumerate}
		\item [(i)] Note that the differentiability of the characteristic function $\phi_X$ plays a crucial role in deriving the Stein identity \eqref{PP1:stidenidd}. Indeed, our approach in deriving Stein identity for infinitely divisible random variables does not follow easily when $\phi_{X}$ is not differentiable. For example, the characteristic function of a Cauchy random variable is not differentiable (see, Section \ref{PP1:preli}, Ex.4.). We need to modify our approach to handle this problem (see Theorem \ref{PP1:Th2}). 
		\item [(ii)] Arras and Houdr\'e \cite[Theorem 3.1]{k1} provide a Stein identity of infinitely divisible random variables using covariance representation for functions of infinitely divisible random variables. However, they assume finite first moment and the function space as bounded Lipschitz. Our proof of Theorem \ref{PP1:Th1} in Section \ref{PP1:proof} is without this assumption, and we consider the Schwartz space $\mathcal{S}(\mathbb{{R}})$ as a suitable function space. It is important to note that the assumption of the first finite moment is an artefact of the technique of covariance representation. Also, our Stein identity \eqref{PP1:stidenidd} exactly matches with the Stein identity given in \cite[Theorem 3.1]{k1}. 
		
		\noindent
		Also, observe that several random variables such as Poisson, negative binomial, normal, Laplace, and gamma can be viewed as infinitely divisible by choosing appropriate triplet $(\beta,\sigma^{2},\nu)$. Stein identities for these random variables can be easily derived using Theorem \ref{PP1:Th1}. In particular, $\alpha$-stable random variables are also infinitely divisible, but the derivation of Stein identity is not straightforward (see, Chen et. al. \cite{k19,k20}).	
	\end{enumerate}

\end{rem}


\noindent
Next, we establish a Stein identity for $\alpha$-stable random variables. As noted in Section \ref{Intro}, the characteristic function of $\alpha$-stable variable is differentiable for $\alpha\in (0,1)\cup (1,2),$ but fails at $\alpha=1$. Hence, we need to modify our approach for $\alpha=1$, using tempered 1-stable random variable $Y_\gamma$ in deriving the Stein identity, as mentioned in Section \ref{PP1:preli}.


\begin{thm} \label{PP1:Th2}
	Let $X\sim\mathcal{S}(\alpha,\beta,m_1,m_2)$ with  characteristic exponent (2). Then,
	\begin{equation}\label{stidstable}
	\mathbb{E}\left((X-\beta)g(X)-\displaystyle\int_{\mathbb R}(g(X+u)-g(X)\mathbf{1}_{\{|u|\leq 1\}}(u))u\nu_{\alpha}(du)\right)=0,~~g \in \mathcal{S}(\mathbb{R}).
	\end{equation} 
\end{thm}

\noindent
Note here that $\sigma=0$, as $X$ has a non-Gaussian stable distribution. The following corollary immediately follows for symmetric $\alpha$-stable random variables by setting $m_1=m_2 = m$, $\beta =0$ and adjusting $\nu_\alpha$. 
\begin{cor}
	Let $X \sim \mathcal{S}(\alpha,0,m,m)$, for $\alpha\in(0,2)$. Then a Stein identity for $X$ is given by
	\begin{equation}\label{e13}
	\mathbb{E} \left(Xg(X)-m\int_{0}^{\infty}\frac{g(X+u)-g(X-u)}{u^{\alpha}}du  \right)=0,~~g \in \mathcal{S}(\mathbb{R}) .
	\end{equation}
\end{cor}

\noindent
In the following remark, we discuss and review the Stein identities available in the literature and in $(iv)$, we compare them with our Stein identities.

\begin{rem}
\begin{enumerate}
	\item[(i)] Chen et. al. \cite[Proposition 2.4]{k19} derive a Stein identity for an $\alpha$-stable random variable with $\alpha\in (0,1],$ using Barbour's generator approach \cite{k8}, where the scale and the location parameters are chosen to be  1 and 0 respectively, for $\alpha\in(0,1)$, and further, for $\alpha=1$, the skewness parameter is also set to zero. Using truncation technique, Arras and Houdr\'e \cite[Theorem 3.1 and Theorem 3.2]{k15} also derive Stein identities for $\alpha$-stable random variables with $\alpha\in (0,1)$ and $\alpha=1$, respectively. 
	
	\item [(ii)] Chen et. al. \cite[Theorem 1.2]{k20} derive a Stein identity for an $\alpha$-stable random variable with $\alpha\in(1,2)$, using Barbour's generator approach \cite{k8}, where the scale and the location parameters are chosen to be  1 and 0 respectively. Chen et. al. \cite{multivariate1} also extend their idea in deriving a Stein identity for multivariate $\alpha$-stable random vectors with $\alpha\in (1,2)$, using Barbour's generator approach \cite{k8}, where the location parameter is chosen to be 0.

	\item [(iii)]  Xu \cite[Theorem 4.1]{k14} derive a Stein identity for a symmetric $\alpha$-stable random variable with $\alpha\in(0,2)$, using invariant measure property of L\'evy-type operators \cite{invariant}. Arras and Houdr\'e \cite[Examples 3.3, (viii)]{k1} also derive a Stein identity for symmetric $\alpha$-stable random variables with $\alpha\in(1,2)$, using covariance representation of functions of infinitely divisible random variables \cite{covrep}.

	 \item [(iv)] Our Stein identity given in \eqref{stidstable} is derived using the L\'evy-Khinchine representation of the characteristic exponent given in \eqref{u3} without any assumption on the scale, location and skewness parameter. Under the assumptions of Chen et. al. \cite{k19,k20},  their Proposition 2.4 and Theorem 1.2 can be retrieved from Theorem \ref{PP1:Th2}. Using Proposition \ref{PP1:appendixPro4}, we see that Stein identities given in \cite[Theorem 3.1 and Theorem 3.2]{k15} exactly match with our Stein identities. For the symmetric case, our Stein identity given in (\ref{e13}) is comparable ($g$ replaced with $g^{\prime}$) to the Stein identities given in \cite[Example 3.3, (viii) and Remark 5.3, (iv)]{k1} and \cite[Theorem 4.1]{k14}.

\end{enumerate}	
\end{rem}

\noindent
As noted in Section \ref{Intro}, the linchpin of Stein's method is the Stein operator $\mathcal{A}$, and the properties of $\mathcal{A}$ play a crucial role in the success of this method. In this context, we adopt the following definition of Stein operator from \cite{gaunt}. 

\begin{defn}\cite[p.1]{gaunt}
	An operator $\mathcal{A}$ is said to be a Stein operator, if $\mathcal{A}$ acts on a class of test functions $\mathcal{G}$ such that $\mathbb{E}(\mathcal{A}g(X))=0$ for all $g\in\mathcal{G},$ where the random variable $X\sim F$.  	
\end{defn}

\begin{rem}
	It is now clear from Theorem \ref{PP1:Th1} that, for an infinitely divisible random variable $X$, $\mathcal{A}_{X}(g)(x)$ $:=(-x+\beta)g(x)+\sigma^{2}g^{\prime}(x)+\displaystyle\int_{\mathbb R}u(g(x+u)-g(x)\mathbf{1}_{\{|u|\leq1\}}(u))\nu(du)$ is an operator acting on $\mathcal{S}(\mathbb{R})$ such that $\mathbb{E}\left(\mathcal{A}_{X}(g)(X)\right)=0$ for all $g\in \mathcal{S}(\mathbb{R})$. Then, by the above definition, $\mathcal{A}_{X}$ is a Stein operator for an infinitely divisible random variable $X$. Also, for any $g\in\mathcal{S}(\mathbb{{R}})$, $\mathcal{A}^{\alpha}_{X}(g)(x):=(-x+\beta)g(x)+\displaystyle\int_{\mathbb R}u(g(x+u)-g(x)\mathbf{1}_{\{|u|\leq1\}}(u))\nu_{\alpha}(du)$ is a Stein operator for an $\alpha$-stable random variable. Observe also that $\mathcal{A}^{\alpha}_{X}$ is an integral operator, where domain of the operator is $\mathcal{F}=\overline{\mathcal{S}(\mathbb{R})}$, the closure of $\mathcal{S}(\mathbb{R})$ (see, \cite{k14} and references therein \cite{stein} for more details).
\end{rem}

\noindent
The existing literature on Stein's method for $\alpha$-stable distributions (see, \cite{k1,k19,k20,multivariate1,jin,k14}) suggests a variety of techniques for deriving a Stein operator depending on the stability parameter $\alpha\in$ $(0,1]$ or $(1,2)$. As mentioned in Section \ref{Intro}, the purpose of this article is to unify Stein's method for $\alpha$-stable distributions. To achieve this, let us use the Stein operator $\mathcal{A}^{\alpha}_{X}$ to set up Stein equation. For any $h\in\mathcal{H}_{X}$ (a class of smooth functions), Stein equation of an $\alpha$-stable random variable $X$ is 
\begin{equation}\label{PP1:a15}
\mathcal{A}^{\alpha}_{X}(g)(x)= h(x)-\mathbb{E}(h(X)). 
\end{equation}
\noindent
 To solve \eqref{PP1:a15}, we use well-known semigroup approach (see, \cite{k8}), and this can be motivated as follows. Recall first that, for $X\sim \mathcal{S}(\alpha,0,m,m)$ with $d_{\alpha}=1$, characteristic function simplifies to 
\begin{align*}
\phi^{s}_{\alpha}(z)=\exp \left( -|z|^{\alpha} \right),~~z\in\mathbb{R}.
\end{align*}
\noindent
Also, observe that, for any $z\in \mathbb{R},$ $\phi^{s}_{\alpha}(z)=\phi^{s}_{\alpha}(e^{-t}z)\phi^{s}_{\alpha}((1-e^{-t})z),~~t\geq0,$ where $\phi^{s}_{\alpha}(e^{-t}z)$ and $\phi^{s}_{\alpha}((1-e^{-t})z)$ denote the characteristic functions of $e^{-t}X$ and $(1-e^{-t})X$ respectively. Note that $e^{-t}X$ and $(1-e^{-t})X$ are symmetric $\alpha$-stable random variables. Let us now generalize this idea for non-symmetric case. One can define a characteristic function, for all $z\in \mathbb{R},$ and $t\geq 0,$ by
\begin{align}\label{PP1:a19}
\phi_t(z):=\frac{\phi_{\alpha}(z)}{\phi_{\alpha}(e^{-t}z)}=\displaystyle\int_{\mathbb{R}}e^{iz u}F_{X_{(t)}}(du),
\end{align}
\noindent
where $F_{X_{(t)}}$ is the distribution function of $X_{(t)}$ and $\phi_{\alpha}$ is the characteristic function of $\alpha$-stable random variable given in \eqref{u3}. The property given in \eqref{PP1:a19} is also known as self-decomposability (see, \cite{k16}).

\noindent
Following Barbour's approach \cite{k8} and using \eqref{PP1:a19}, we choose a family of operators $(P^{\alpha}_{t})_{t\geq 0},$ for all $x\in \mathbb{R},$ as
\begin{equation}\label{PP1:a20}
P^{\alpha}_t(g)(x):=\frac{1}{2\pi}\int_{\mathbb R}\widehat{g}(z)e^{iz xe^{-t}}\phi_{t}(z) dz,~~g\in \mathcal{F}.
\end{equation} 

\noindent
Using \eqref{PP1:a19}, we get
\begin{align}
\nonumber  P^{\alpha}_t(g)(x)&=\frac{1}{2\pi}\int_{\mathbb R}\int_{\mathbb{R}}\widehat{g}(z)e^{iz xe^{-t}}e^{iz u}F_{X_{(t)}}(du)dz\\
\nonumber &=\frac{1}{2\pi}\int_{\mathbb R}\int_{\mathbb{R}}\widehat{g}(z)e^{iz (u+xe^{-t})}F_{X_{(t)}}(du)dz\\
&=\displaystyle\int_{\mathbb{R}}g(u+xe^{-t})F_{X_{(t)}}(du),\label{PP1:a21}
\end{align}
\noindent
where the last step follows by applying inverse Fourier transform.


\begin{pro}\label{PP1:proSem}
	The family of operators $(P^{\alpha}_{t})_{t\geq 0}$ given in \eqref{PP1:a20} is a $\mathbb{C}_0$-semigroup on $\mathcal{F}$.
\end{pro}
\noindent
\textbf{Proof.} For each $g\in\mathcal{F}$, it is easy to show that $P^{\alpha}_{0}g(x)=g(x) ~~\text{and}~~ \lim_{t\to\infty}P^{\alpha}_{t}(g)(x)=\displaystyle\int_{\mathbb R}g(x)F^{\alpha}_{X}(dx)$. 
\noindent
Now, for any $s,t\geq0$, we have
\begin{align}\label{PP1:a22}
\phi_{t+s}(z)=\frac{\phi_{\alpha}(z)}{\phi_{\alpha}(e^{-(t+s)}z)}
=\frac{\phi_{\alpha}(z)}{\phi_{\alpha}(e^{-s}z)}\frac{\phi_{\alpha}(e^{-s}z)}{\phi_{\alpha}(e^{-(t+s)}z)}
=\phi_{s}(z)\phi_{t}(e^{-s}z)
\end{align}
\noindent
Using \eqref{PP1:a22}, we have
\begin{align}\label{PP1:a022}
LHS=P_{t+s}^{\alpha}(g)(x)=\frac{1}{2\pi}\int_{\mathbb R}\widehat{g}(z)e^{iz xe^{-(t+s)}}\phi_{t+s}(z)dz=\frac{1}{2\pi}\int_{\mathbb R}\widehat{g}(z)e^{iz xe^{-(t+s)}}\phi_{s}(z)\phi_{t}(e^{-s}z) dz.
\end{align}
\noindent
We need to show that $P_{t+s}^{\alpha}(g)(x)=P^{\alpha}_{t}(P^{\alpha}_{s}g)(x)$ for all $g\in \mathcal{F}$.

\begin{align*}
\text{RHS}&=P^{\alpha}_{t}(P^{\alpha}_{s}(g))(x)\\ &=\frac{1}{2\pi}\int_{\mathbb R}\widehat{P^{\alpha}_{s}(g)}(z)e^{iz xe^{-t}}\phi_{t}(z)dz\\
&=\frac{1}{2\pi}\int_{\mathbb R}\left(\int_{\mathbb{R}}e^{-ivz}P_{s}^{\alpha}(g)(v)dv  \right)e^{iz xe^{-t}}\phi_{t}(z)dz\\
&=\frac{1}{(2\pi)^{2}}\int_{\mathbb R}\left(\int_{\mathbb{R}}e^{-ivz}\left(\int_{\mathbb{R}}\widehat{g}(w)e^{iwe^{-s}v}\phi_{s}(w)dw  \right)dv  \right)e^{iz xe^{-t}}\phi_{t}(z)dz\\
&=\frac{1}{(2\pi)^{2}}\int_{\mathbb R}\widehat{g}(w)\phi_{s}(w)\int_{\mathbb R}e^{iz xe^{-t}}\phi_{t}(z) \left(\int_{\mathbb R}e^{iv(e^{-s}w-z)}dv  \right)dz dw\\
&=\frac{1}{(2\pi)^{2}}\int_{\mathbb R}\widehat{g}(w)\phi_{s}(w)\int_{\mathbb R}e^{iz xe^{-t}}\phi_{t}(z) 2\pi \delta (e^{-s}w -z) dz dw~~\text{(where $\delta$ is the Dirac-$\delta$ measure )}\\
&=\frac{1}{2\pi}\int_{\mathbb R}\widehat{g}(w)\phi_{s}(w)e^{i e^{-s}wxe^{-t}}\phi_{t}(e^{-s}w)dw\\
&=\frac{1}{2\pi}\int_{\mathbb R}\widehat{g}(z)e^{iz xe^{-(t+s)}}\phi_{s}(z)\phi_{t}(e^{-s}z) dz\\
&=P_{t+s}^{\alpha}(g)(x)=\text{LHS}~~(\text{from }\eqref{PP1:a022}),
\end{align*}   
\noindent
and the desired conclusion follows.\vspace{2mm}

\noindent
Next, we find the generator of the semigroup defined in \eqref{PP1:a20}.
\begin{lem}
	Let $(P^{\alpha}_{t})_{t\geq 0}$ be a $\mathbb{C}_0$-semigroup as defined in \eqref{PP1:a20}. Then, its generator $\mathcal{T}_{\alpha}$ is given by
	
	\begin{equation*}
	\mathcal{T}_{\alpha}g(x)=(-x+\beta)g^{\prime}(x)+\int_{\mathbb{R}}\left(g^{\prime}(x+u)-g^{\prime}(x)\mathbf{1}_{\{|u|\leq1 \}}(u)\right)u\nu_{\alpha}(du),~~g\in\mathcal{S}(\mathbb{R}),
	\end{equation*}
	where $\alpha\in (0,1)\cup (1,2).$	
\end{lem}
\noindent
\textbf{Proof.} The proof of this lemma is split into two parts.

\noindent
$\mathbf{(i)~~\alpha\in(0,1)}:$ For all $g\in \mathcal{S}(\mathbb{R}),$
\begin{align}
\nonumber \mathcal{T}_{\alpha}(g)(x)&=\lim_{t\to0^{+}}\frac{1}{t}\left(P^{\alpha}_{t}(g)(x) -g(x)   \right)\\
\nonumber &=\frac{1}{2\pi}\lim_{t\to0^{+}}\int_{\mathbb R}\widehat{g}(z)e^{iz x} \frac{1}{t}\left(e^{iz x(e^{-t}-1)}\phi_{t}(z)-1\right)dz \\ 
\nonumber&=\frac{1}{2\pi}\int_{\mathbb R}\widehat{g}(z)e^{iz x}\left(-x+\beta- \int_{\{|u|\leq 1 \}}u \nu_{\alpha}(du)+ \displaystyle\int_{\mathbb{R}}e^{iz u }u\nu(du) \right )(iz)dz \text{ (using Prop. \ref{PP1:appendixPro2} )}\\
\nonumber &=\frac{1}{2\pi}\int_{\mathbb R}\widehat{g}(z)e^{iz x}\left(-x+\beta_1+ \displaystyle\int_{\mathbb{R}}e^{iz u }u\nu(du) \right )(iz)dz \text{ (where $\beta_{1}=\beta- \int_{\{|u|\leq 1 \}}u \nu_{\alpha}(du)$)}\\ 
\nonumber&=(-x+\beta_{1})g^{\prime}(x)+\int_{\mathbb{R}}g^{\prime}(x+u)u\nu_{\alpha}(du)\\
\nonumber&=(-x+\beta)g^{\prime}(x)+\int_{\mathbb{R}}(g^{\prime}(x+u)-g^{\prime}(x)\mathbf{1}_{\{|u|\leq 1\}}(u))u\nu_{\alpha}(du),
\end{align} 
\noindent
where the last equality follows by splitting $\beta_{1}$ (see \eqref{PP1:u4}).\vspace{2mm}

\noindent
$\mathbf{(ii)~~\alpha\in(1,2)}:$ For all $g\in \mathcal{S}(\mathbb{R}),$
\begin{align}
\nonumber \mathcal{T}_{\alpha}(g)(x)&=\lim_{t\to0^{+}}\frac{1}{t}\left(P^{\alpha}_{t}(g)(x) -g(x)   \right)\\
\nonumber &=\frac{1}{2\pi}\lim_{t\to0^{+}}\int_{\mathbb R}\widehat{g}(z)e^{iz x} \frac{1}{t}\big(e^{iz x(e^{-t}-1)}\phi_{t}(z)-1\big)dz\\ 
\nonumber &= \frac{1}{2\pi}\int_{\mathbb R}\widehat{g}(z)e^{iz x}\left(-x+\beta+ \int_{\{|u|> 1 \}}u \nu_{\alpha}(du)+ \displaystyle\int_{\mathbb{R}}(e^{iz u }-1)u\nu(du) \right )(iz)dz \text{ (using Prop. \ref{PP1:appendixPro3})} \\
\nonumber &= \frac{1}{2\pi}\int_{\mathbb R}\widehat{g}(z)e^{iz x}\left(-x+\beta_{2}+ \displaystyle\int_{\mathbb{R}}(e^{iz u }-1)u\nu(du) \right )(iz)dz \text{ (where $\beta_{2}=\beta+ \int_{\{|u|> 1 \}}u \nu_{\alpha}(du)$)} \\
\nonumber&=(-x+\beta_{2})g^{\prime}(x)+\int_{\mathbb{R}}(g^{\prime}(x+u)-g^{\prime}(x))u\nu_{\alpha}(du)\\
\nonumber&=(-x+\beta)g^{\prime}(x)+\int_{\mathbb{R}}(g^{\prime}(x+u)-g^{\prime}(x)\mathbf{1}_{\{|u|\leq 1\}}(u))u\nu_{\alpha}(du),
\end{align}
\noindent
where the last equality follows by splitting $\beta_{2}$ (see \eqref{PP1:u5}).

\noindent
This completes the proof.\vspace{2mm}

\noindent
Next we handle the case $\alpha=1$, using tempered 1-stable random variable $Y_\gamma$. Recall that the characteristic exponent of $Y_\gamma$ is given in \eqref{PP1:cetemst}. Let $\phi_{1,\gamma}(z):=e^{\eta_{1,\gamma}(z)},$ $z\in\mathbb{{R}}$ be the characteristic function of $Y_\gamma$. Then, for all $z\in\mathbb{{R}}$, we define

%
\begin{align}\label{PP1:a25}
\phi_{1,t,\gamma}(z):=\frac{\phi_{1,\gamma}(z)}{\phi_{1,\gamma}(e^{-t}z)},~~  t \geq 0.
\end{align}
\noindent
Using \cite[Corollary 15.11]{k16}, it can be easily shown that $\phi_{1,t,\gamma}$ is a well-defined characteristic function.

\noindent
Now, using \eqref{PP1:a25}, define a family of operators $(P^{1,\gamma}_{t})_{t \geq 0},$ for all $x\in \mathbb{R},$ by
\begin{align}\label{PP1:a26}
P^{1, \gamma}_t(f)(x)&:=\frac{1}{2\pi}\int_{\mathbb R}\widehat{g}(z)e^{iz xe^{-t}}\phi_{1,t,\gamma}(z)dz=\int_{\mathbb{R}}g(u+xe^{-t})F_{Y_{(\gamma,t)}}(du) ,~~g\in\mathcal{F}.
\end{align} 
\noindent
Following similar steps as Proposition \ref{PP1:proSem}, one can show that $(P^{1,\gamma}_{t})_{t \geq 0},$ is a $\mathbb{C}_0$-semigroup on $\mathcal{F}$.\vspace{2mm} 

\noindent
Next, we obtain a generator for the semigroup \eqref{PP1:a26}.
\begin{lem}
	Let $(P^{1,\gamma}_{t})_{t\geq0}$ be a $\mathbb{C}_0$-semigroup as defined in \eqref{PP1:a26}. Then, its generator $\mathcal{T}_{1,\gamma}$ is given by
	
	\begin{equation*}
	\mathcal{T}_{1,\gamma}g(x)=-xg^{\prime}(x)+\int_{\mathbb{R}}\left(g^{\prime}(x+u)-g^{\prime}(x)\mathbf{1}_{\{|u|\leq1 \}}(u)\right)u\nu_{1,\gamma}(du),~~g\in\mathcal{S}(\mathbb{R}).
	\end{equation*}	
\end{lem}
\noindent
\textbf{Proof.} For all $g\in\mathcal{S}(\mathbb{R}),$

\noindent
we get,
\begin{align}
\nonumber \mathcal{T}_{1,\gamma}(g)(x)&=\lim_{t\to0^{+}}\frac{1}{t}\left(P^{1,\gamma}_{t}(g)(x) -g(x)   \right)
\end{align}
\begin{align}
\nonumber \mathcal{T}_{1,\gamma}(g)(x)&=\frac{1}{2\pi}\lim_{t\to0^{+}}\int_{\mathbb R}\widehat{g}(z)e^{iz x} \frac{1}{t}\left(e^{iz x(e^{-t}-1)}\frac{\phi_{1,\gamma}(z)}{\phi_{1,\gamma}(e^{-t}z)}-1\right)dz\\
\nonumber &=\frac{1}{2 \pi}\int_{\mathbb{R}} \widehat{g}(z)e^{izx}\left(-x+\beta+\displaystyle\int_{\mathbb{R}}(e^{i\xi u }-1_{\{|u|\leq 1 \}})u\nu_{1,\gamma}(du) \right)(iz)dz\\
\nonumber &=(-x+\beta)g^{\prime}(x)+\displaystyle\int_{\mathbb{R}}(g^{\prime}(x+u)-g^{\prime}(x)\mathbf{1}_{\{|u|\leq 1 \}}  )u\nu_{1,\gamma}(du),
\end{align}
\noindent
where the last but one equality follows by doing computations similar to Proposition \ref{PP1:appendixPro3}.

\noindent
This completes the proof. \vspace{2mm}

\noindent
	Now, observe that, $\lim _{\gamma \to 0^{+}}P_t^{1,\gamma}g(x) = P_t^{1}g(x)$, $g\in \mathcal{F}$, as defined in \eqref{PP1:a20}. Hence $\lim_{\gamma \to 0^{+}}\mathcal{T}_{1,\gamma} = \mathcal{T}_1$, where $\mathcal{T}_1$ is given by
	\begin{align*}
	\mathcal{T}_1g(x)=(-x+\beta)g^{\prime}(x)+\displaystyle\int_{\mathbb{R}}(g^{\prime}(x+u)-g^{\prime}(x)\mathbf{1}_{\{|u|\leq 1 \}}  )u\nu_{1}(du),~~g\in\mathcal{S}(\mathbb{R}).
	\end{align*}

\begin{rem}
	Note that, on careful adjustments of the integrals with respect to $\nu_\alpha$ for $\alpha\in (0,2)$, we see that the operator $\mathcal{T}_\alpha$ is also a Stein operator for an $\alpha$-stable random variable, see \cite{k19,k20,jin,k14}. In these articles, authors use $\mathcal{T}_\alpha$ to set up their Stein equations. However, we consider $\mathcal{A}^{\alpha}_{X}$ to set up our Stein equation.
\end{rem}
\noindent
Next, we provide the solution to our Stein equation.

\begin{thm}\label{PP1:ThmSol}
	\begin{enumerate}
		\item [(i)] For $\alpha\in (0,1]$, let $X\sim \mathcal{S}(\alpha,\beta,m_1,m_2)$ and $h\in \mathcal{H}_\delta,$ $0<\delta<\alpha.$ Also, let
		\begin{equation}\label{PP1:steinequ1}
		\mathcal{A}^{\alpha}_{X}(g)(x)=h(x)-\mathbb{E}h(X)
		\end{equation}
		\noindent
		be a Stein equation for $X$. Then, the function $g_{h}^{\alpha}:\mathbb{{R}}\to\mathbb{{R}}$ defined as
	 \begin{equation}\label{PP1:solution1}
		g^{\alpha}_{h}(x)=
		-\displaystyle\int_{0}^{\infty}e^{-t}\int_{\mathbb{R}}h^{\prime}(u+xe^{-t})F_{X_{(t)}}(du)dt,
	\end{equation}
		\noindent
		solves \eqref{PP1:steinequ1}.

		\item [(ii)] For $\alpha\in(1,2)$, $X\sim \mathcal{S}(\alpha,\beta,m_1,m_2)$ and $h\in \mathcal{H}_2.$ Also, let
		
		\begin{equation}\label{PP1:steinequ2}
		\mathcal{A}^{\alpha}_{X}(g)(x)=h(x)-\mathbb{E}h(X)
		\end{equation}
		\noindent
		be a Stein equation for $X$. Then, the function $g_{h}^{\alpha}:\mathbb{{R}}\to\mathbb{{R}}$ defined as
		\begin{equation}\label{PP1:solution2}
		g^{\alpha}_{h}(x)=
		-\displaystyle\int_{0}^{\infty}e^{-t}\int_{\mathbb{R}}h^{\prime}(u+xe^{-t})F_{X_{(t)}}(du)dt,
		\end{equation}
		\noindent
		solves \eqref{PP1:steinequ2}.
		\end{enumerate}
\end{thm}
\noindent
In the following remark, we review the techniques used by several authors to solve Stein equations under various constraints and justify our claim of unification.

\begin{rem}
	\begin{enumerate}

		\item[(i)]Chen et. al. \cite{k19,k20} derive the solution to Stein equation for an $\alpha$-stable random variable with $\alpha\in(0,1]$ and $\alpha\in(1,2)$ respectively using Barbour's generator approach \cite{k8}, and the transition density function of $\alpha$-stable processes. Xu $\cite{k14}$ also uses the Barbour's generator approach to solve Stein equation for a symmetric $\alpha$-stable random variable with $\alpha\in(1,2)$. Arras and Houdr\'e \cite{k1} provide the semigroup approach to solve Stein equation for an infinitely divisible random variable with the first finite moment. Recently, Arras and Houdr\'e \cite[Remark 4.3]{k15} show that semigroup approach for deriving the solution to Stein equation is also applicable for multivariate $\alpha$-stable random vectors with $\alpha\in(0,1)$,  and they also mention that the semigroup approach is also applicable for $\alpha=1$, but requires different estimates.
		
		\item[(ii)]Note that, for both parts of Theorem \ref{PP1:ThmSol}, we use only the semigroup approach to solve Stein equation for an $\alpha$-stable random variable with $\alpha\in(0,2)$, and this unifies the method of solving the Stein equation for $\alpha$-stable random variables.
		
	\end{enumerate}
\end{rem}

\subsection{Properties of solution to Stein equation}
\noindent 
Let us now study regularity estimates of the solution to our Stein equation. Recall that the L\'evy measure $\nu_{\alpha}$ for $\alpha$-stable distributions is given by $\nu_{\alpha}(du)=\left(m_1\frac{1}{|u|^{1+\alpha}}\mathbf{1}_{(0,\infty)}(u)+m_2\frac{1}{|u|^{1+\alpha}}\mathbf{1}_{(-\infty,0)}(u)\right)du,$ where $m_1,m_2\in [0,\infty)$, $m_1+m_2>0$ and $\alpha\in (0,2)$. In the following theorem, we establish estimates of $g^{\alpha}_{h}$, which play a crucial role in the $\alpha$-stable approximation problem.
\begin{thm}\label{PP1:Th4}
	$(i)$ Let $\alpha\in(0,1)$. For $h\in \mathcal{H}_{\delta}$, $0<\delta<\alpha$, let $g^{\alpha}_{h}$ be defined in \eqref{PP1:solution1}. Then, for any $x,y\in \mathbb{R}$  
	\begin{align}
	\label{PP1:a49}	\left\|g^{\alpha}_{h}\right\|&\leq\left\|h^{\prime}\right\|,\\
	\label{PP1:a50}	\left\|g^{\alpha}_{h}(x)-g^{\alpha}_{h}(y)  \right\|&\leq \frac{1}{\delta+1}|x-y|^{\delta}.
	\end{align}
	\noindent
	Define $A_0g^{\alpha}_{h}(x):=\displaystyle\int_{\mathbb{R}}ug^{\alpha}_{h}(x+u)\nu_{\alpha}(du).$ Then
	\begin{align}
	\label{PP1:a51}	\left\|A_0g^{\alpha}_{h} \right\|\leq C_{\alpha,\delta,m_1,m_2}:=\frac{\alpha(m_1+m_2)}{\delta(\alpha-\delta)}+\frac{\alpha(m_1-m_2)}{1-\alpha},
	\end{align}
	\noindent
	$(ii)$ Let $\alpha=1$. For $h\in \mathcal{H}_{\delta}$, $0<\delta<1$, let $g^{1}_{h}$ be defined in \eqref{PP1:solution1}. Then, for any $x,y\in \mathbb{R}$  
	
	\begin{align}
	\label{PP1:a52}	\left\|g^{1}_{h}\right\|&\leq\left\|h^{\prime}\right\|,\\
	\label{PP1:a53}	\left\|g^{1}_{h}(x)-g^{1}_{h}(y)  \right\|&\leq \frac{1}{\delta+1}|x-y|^{\delta}.
	\end{align}
	\noindent
	Define $A_1g^{1}_{h}(x):=\displaystyle\int_{\mathbb{R}}u(g^{1}_{h}(x+u)-ug^{1}_{h}\mathbf{1}_{\{|u|\leq 1 \}})\nu_{1}(du)$. Then
	\begin{align}
	\label{PP1:a54}	\left\|A_1g^{1}_{h} \right\|\leq C_{1,\delta,m_1,m_2}:=\frac{2(m_1+m_2)}{\delta(1-\delta^{2})}.
	\end{align}

	\noindent
	$(iii)$ Let $\alpha\in (1,2).$ For $h\in\mathcal{H}_{2}$, let $g^{\alpha}_{h}$ be defined in \eqref{PP1:solution2}. Then, $g^{\alpha}_{h}$ is differentiable on $\mathbb{R}$, 
	\begin{equation}\label{PP1:a55}
	\left\|g^{\alpha}_{h}\right\|\leq \left\|h^{\prime}\right\|~~\text{and}~~ \left\|(g^{\alpha}_{h})^{\prime}\right\|\leq\frac{1}{2} \left\|h^{\prime \prime}\right\|.
	\end{equation}
	Let $m_1=m_2=m$. Define $A_{2}g^{\alpha}_{h}(x):=\displaystyle\int_{\mathbb{{R}}}(g^{\alpha}_{h}(x+u)-g^{\alpha}_{h}(x))u\nu_{\alpha}(du)$. For any $x,y\in \mathbb{R}$
	\begin{align}\label{PP1:a56}
	\|	A_2g^{\alpha}_{h}(x)-A_2g^{\alpha}_{h}(y)\| \leq C_{\alpha,m}\|h^{\prime\prime}\||x-y|^{2-\alpha},
	\end{align}
	where $C_{\alpha,m}$ is a positive constant.	
\end{thm}

\noindent
In the following remark, we review the estimates of solution to Stein equation available in the literature. In $(iii)$ and $(iv)$, we compare them with our results.
\begin{rem}
	\begin{enumerate}
		
		\item [(i)] Chen et. al. \cite[Proposition 1 and Proposition 2]{k19} provide estimates of solution to Stein equation for $\alpha\in(0,1]$, using the properties of transition density of an $\alpha$-stable process. From the derivation of these estimates, we observe that the scaling property of transition density plays a crucial role in deriving the estimates. In particular, authors prove that for any $x,y\in \mathbb{R}$, $|f_{h}^{\prime}(x)-f^{\prime}_{h}(y)| \leq C|x-y|^{\alpha}$, where $f_h$ is the solution to Stein equation and $C$ is some positive constant, see \cite[Proposition 1 (Eq.3.2)]{k19}. This fact helps them to determine a good rate of convergence and this is discussed in more detail in Section \ref{PP1:application}.
		
		\item [(ii)] Chen et. al. \cite{k20} and Jin et. al. \cite{jin} derive estimates of solution to Stein equation for $\alpha\in (1,2)$ for an $\alpha$-stable random variable with $\alpha\in (1,2)$, using scaling property of transition density of an $\alpha$-stable process. Arras and Houdr\'e \cite[Section 5]{k1} also derive estimates of solution to Stein equation for an $\alpha$-stable random variable with $\alpha\in (1,2)$, using self-decomposability and properties of the semigroup used to solve their Stein equation. In \cite[Section 5]{k1}, we see that for any $h\in\mathcal{H}_2$, upper bounds of first and second derivative of solution to Stein equation heavily rely on first and second derivative of $h$, respectively (due to the choice of Stein equation $\mathcal{T}_\alpha g(x)=h(x)-\mathbb{E}h(X)$, where $X\sim \mathcal{S}(\alpha,\beta,m_1,m_2)$).
		
		\item [(iii)] Observe that, for any $h\in\mathcal{H}_{\delta}$ and $\alpha\in (0,1]$, the upper bound of the solution to our Stein equation $g^{\alpha}_{h}$ relies on $h^{\prime}$. Note that, for $\alpha\in(0,1]$, other estimates of $g^{\alpha}_{h}$ are derived by suitably adjusting the integrals, and using the properties of the semigroup defined in \eqref{PP1:a20}. These estimates are used to obtain the rate of convergence for $\alpha$-stable approximations with $\alpha\in(0,1)$, and we obtain a flexible rate of convergence (faster rate when $\alpha\in (0,0.5]$ and slower rate when $\alpha\in (0.5,1)$). This is discussed in more detail in Section \ref{PP1:application}. 
		
		\item [(iv)] Observe also that, for $h\in\mathcal{H}_2$ and $\alpha\in(1,2)$, the upper bound of the solution to Stein equation, $g^{\alpha}_{h}$ relies on first derivative of $h$, and upper bound for the first derivative of $g^{\alpha}_{h}$ relies on the second derivative $h$. Other estimates of $g^{\alpha}_{h}$ are derived by suitably adjusting the integrals, and using the properties of the semigroup defined in \eqref{PP1:a20}.  These estimates provide a better rate of convergence for $\alpha$-stable approximations, whenever $\alpha\in (1,2)$ and this is also discussed in more detail in Section \ref{PP1:application}.
		

	\end{enumerate}
\end{rem}

\noindent
Next, we provide bounds in the $d_{W_\delta}$ distance for $\alpha$-stable approximations of a partial sum of sequence of i.i.d random variables in the domain of normal attraction of an $\alpha$-stable distribution with $\alpha\in(0,1]$. Let us define the domain of normal attraction of an $\alpha$-stable distribution as follows.
\begin{defn}\label{PP1:def2} 
	\cite[ p.4]{k19}
	A real-valued random variable $Y$ is said to be in the domain of normal attraction of an $\alpha$-stable distribution with $\alpha\in(0,1]$ if its CDF, $F_{Y}$  satisfies
	\begin{equation}\label{PP1:a68}
	1-F_{Y}(y)=\frac{A+e(y)}{|y|^{\alpha}}(1+\theta)\; \text{and}\; F_{Y}(-y)=\frac{A+e(-y)}{|y|^{\alpha}}(1-\theta),
	\end{equation}
	where $y>1$, $\alpha\in(0,1]$, $\theta=\frac{m_1-m_2}{m_1+m_2}\in[-1,1]$, $A(>0)$ a constant and $e:\mathbb{R}\to\mathbb{R}$ is a bounded differentiable function vanishing at $\pm\infty$. Since $e$ is bounded, we denote $K=\sup_{y\in\mathbb{R}}e(y)$.
\end{defn}
\noindent
We denote $Y\in D_\alpha$, if $Y$ is in the domain of normal attraction of an $\alpha$-stable distribution, and for a positive constant $L$, the function $e$ defined in \eqref{PP1:a68} is $C^{2}$ with the domain $\{|y|>L \}$, and it satisfies $ye^{\prime}(y)\to 0$ and $y^{2}e^{\prime\prime}(y)\to 0$ as $|y|\to \infty$.


\begin{thm}\label{PP1:Th5}
	Let $Y_{1},Y_{2},\ldots,Y_{n}$ be a sequence of i.i.d random variables such that $Y_{i}\in D_\alpha$ and $X\sim \mathcal{S}(\alpha,\beta,m_1,m_2)$ with $\alpha\in (0,1]$. Define  $S_{n}=n^{-1/\alpha}(Y_{1}+Y_{2}+\cdots+Y_{n})$. Then,
	
	\noindent
	\textbf{(a)} for $\alpha\in (0,1)$
	\begin{align}
	\nonumber 	d_{W_{\delta}}(S_{n},X)&\leq C_{\alpha,\delta,m_1,m_2}^{A,K}n^{-1}+C_{1,\delta,L} n^{1-\frac{(1+\delta)}{\alpha}}+C_{2,\delta}n^{1-\frac{(1+\delta)}{\alpha}}\sup_{L<|y|<n^{\frac{1}{\alpha}}} \left( \alpha|e(y)|+|ye^{\prime}(y)| \right)\displaystyle\int_{L<|y|<n^{\frac{1}{\alpha}}}|y|^{\delta -\alpha}dy\\
\label{PP1:bound1}	&+C_{\alpha,\delta,m_1,m_2} n^{-\frac{(1-\alpha)}{\alpha}}+n^{-\frac{(1-\alpha)}{\alpha}}\displaystyle\int_{|y|<n^{\frac{1}{\alpha}}}|y|dF_{Y}(y)+R_{\alpha,n},
	\end{align}
	\noindent
	where $C_{\alpha,\delta,m_1,m_2}^{A,K}, C_{1,\delta,L},C_{2,\delta}$, $C_{\alpha,\delta,m_1,m_2}$ are positive constants, and
	$$R_{\alpha,n}=\beta_{1}+2\sup_{|y|>n^{\frac{1}{\alpha}}}\left(\alpha |e(y)|+|ye^{\prime}(y)| \right)\displaystyle\int_{|y|> 1 }\eta_{\alpha,\beta,\delta,m_1,m_2}(y)|y|^{-1-\alpha}dy.$$
	
	\noindent
	\textbf{(b)} For $\alpha=1$
	\begin{align}
	\nonumber	d_{W_{\delta}}(S_{n},X) &\leq  C_{1,\delta,m_1,m_2}^{A,K}n^{-1} +\frac{1}{\delta +1} n^{-\delta} \left(L^{2}+m_1+m_2  \right)+\frac{n^{-\delta}}{1+\delta}\displaystyle\int_{L<|u|<\frac{1}{a}} \frac{| e(u)-ue^{\prime}(u)|}{|u|^{1-\delta}}du\\
\label{PP1:bound2}	&+n^{-1} \displaystyle\int_{|u|>1}\frac{|e(nu)-nue^{\prime}(nu)|}{|u|}du+R_{1,n},
	\end{align}
	where $C_{1,\delta,m_1,m_2}^{A,K},C_{1,\delta,m_1,m_2}$ are positive constants, and
	$$R_{1,n}=\beta+2K+ 2 C_{1,\delta,m_1,m_2} \int_{0}^{n}dF_{|Y|}(y)+\left|\displaystyle\int_{0}^{n}\frac{e(y)-e(-y)}{y}dy \right|.$$
\end{thm}

\begin{rem}
	Note that, in view of Theorem \ref{PP1:Th5}, we consider only real-valued random variables $Y_i\in D_\alpha.$ Indeed, integer-valued random variables in general do not belong to $D_\alpha,$ see \cite{k19}. The problem for developing an approach that allow to handle integer-valued sums is still open. Recently, Chen et. al. \cite{k19} also provide bounds in $d_{W_\delta}$ distance for $\alpha$-stable approximations of a partial sum of sequence of i.i.d random variables, belong to $D_\alpha$. Our bounds given in \eqref{PP1:bound1} and \eqref{PP1:bound2} are similar to the bounds given in \cite[Theorem 4]{k19}. Note also that, our bounds include location parameter on $\alpha$-stable approximations with $\alpha\in(0,1]$. Chen et. al. \cite{k19}, Chen and Xu \cite{k18} give error bounds on $\alpha$-stable approximation with $\alpha\in(0,1]$ by choosing the location parameter to be 0. 
\end{rem}
\noindent
Now, we present bound in the $d_{W_2}$ distance  for $\alpha$-stable approximations of a partial sum of sequence of i.i.d random variables  with $\alpha\in(1,2)$. To derive this bound, we apply kernel decomposition method introduced by Xu \cite{k14} for symmetric $\alpha$-stable approximations. Later, Arras and Houdr\'e \cite{k1} generalized it for IDD with first finite moment. Before stating our result, define 
\begin{align*}
K_{\nu_{\alpha}}(t,N)&=\mathbf{1}_{[0,N]}(t)\displaystyle\int_{t}^{N}u\nu_{\alpha}(du)+\mathbf{1}_{[-N,0]}(t)\displaystyle\int_{-N}^{t}(-u)\nu_{\alpha}(du),\text{ and}\\
K_{i}(t,N)&=\mathbb{E}\left(Z_{i}\mathbf{1}_{\{0\leq t\leq Z_{i}\leq N\}} -Z_{i}\mathbf{1}_{\{-N\leq Z_{i}\leq t\leq0\}}  \right),
\end{align*}
\noindent
where $N>0$ is an arbitrary number. Our theorem is as follows.
\begin{thm}\label{PP1:th6}
Let $Y_{1},Y_{2},\ldots,Y_{n}$ be a sequence of i.i.d random variables with $\mathbb{E}Y_{i}=0$ and $\mathbb{E}|Y_{i}|<\infty$. Let $X\sim \mathcal{S}(\alpha,0,m,m)$ with $\alpha\in (1,2)$. Define $Z_i=n^{-\frac{1}{\alpha}}Y_{i}$ and $S_{n}=Z_1+Z_2+\ldots+Z_n$. Then,
	\begin{align}
	d_{W_{2}}(S_{n},X)&\leq\frac{1}{2}\sum_{i=1}^{n}\int_{-N}^{N}\left|\frac{K_{\nu_{\alpha}}(t,N)}{n}-K_{i}(t,N)\right|dt+R_{n,N},~~N>0 , \label{e23}
	\end{align}
	where $R_{n,N}=2\left( \int_{|u|>N}|u|\nu_{\alpha}(du)+\sum_{i=1}^{n}\mathbb{E}\left[ |Z_{i}|\mathbf{1}_{ \{|Z_{i}|>N \}}  \right]\right)  + \frac{C_{\alpha,m}}{n}\sum_{i=1}^{n}\mathbb{E}|Z_i|^{2-\alpha}$.
\end{thm}
\begin{rem}
	In the existing literature, several authors use this kernel discrepancy type bound for $\alpha$-stable approximations with $\alpha\in (1,2)$, see \cite{k1,jin,k14}. In these articles, we note that the derivation of this bound heavily depends on the upper bound for the second derivative of solution to Stein equation.
	 However, our kernel discrepancy bound \eqref{e23} mainly depends on upper bound of first derivative of solution to our Stein equation, as our Stein equation \eqref{PP1:a15} is an integral equation, and our bound given in \eqref{e23} is comparable to the bounds given in \cite{k1,jin,k14}. 
\end{rem}

\section{Applications}\label{PP1:application}
\noindent
In this section, we discuss the rates of convergence for $\alpha$-stable approximations using two examples and we compare them with existing literature.
 
 \begin{ex}[Pareto distribution with $\alpha\in(0,1)$\cite{k19,kuske}]

 \noindent
 Assume that $Y_1,Y_2,\ldots,Y_n$ be i.i.d random variables having a Pareto distribution with $\alpha\in(0,1),$ i.e.
 
 $$P(Y_1>y)=\frac{1}{2|y|^{\alpha}},y\geq 1,~~P(Y_1\leq y)=\frac{1}{2|y|^{\alpha}},y\leq -1.$$
 
 \noindent
 In this case $\theta=0, A=\frac{1}{2},e(y)=\frac{|y|^{\alpha}-1}{2}\textbf{1}_{(-1,1)}$ and $K=\frac{1}{2}.$ Observe that $L=1,e(y)=0$ for $|y|>1.$ By Theorem \ref{PP1:Th5}, Case 1, one can easily verify that
 $$d_{W_{\delta}}(S_n,X) \leq C n^{-(\frac{1}{\alpha}-1)},$$
 
 \noindent
 where $C$ is some positive constant. Moreover, $d_{W_\delta}(S_n,X)=O(n^{-(\frac{1}{\alpha}-1)}).$\vspace{2mm}
 
 \noindent
 Let us compare our result with the known results in literature.
 
 \noindent
 The reference \cite{kuske} shows a convergence rate $d_{\text{Kol}(S_n,X)}\leq C_\alpha n^{-1},$ for $\alpha\in (0,1],$ where an exact value of $C_\alpha$ was not given. Chen et. al. \cite{k19} proved that the rate $n^{-1}$ is valid for the $d_{W_\delta}$ distance, whenever $\alpha\in (0,1).$ For $\alpha\in (0,1)$, our rate is $n^{-(\frac{1}{\alpha}-1)}$, which is flexible. In comparison with the rates derived in \cite{k19,kuske}, we see that our rate is faster ($\alpha\in (0,0.5)$), same ($\alpha=0.5$) and slower ($\alpha\in (0.5,1)$).	
 \end{ex}

%
%
%
%
%

	
\begin{ex}[Pareto distribution with $\alpha\in(1,2)$\cite{john,k14}]
	
\noindent
Assume that $Y_1,Y_2,\ldots,Y_n$ be i.i.d random variables having a Pareto distribution with $\alpha\in(1,2),$ i.e.

$$P(Y_1>y)=\frac{1}{2|y|^{\alpha}},y\geq 1,~~P(Y_1\leq y)=\frac{1}{2|y|^{\alpha}},y\leq -1.$$
\noindent

\noindent
Assume also that $\beta=0$ and $m_1=m_2=m.$ Denote $Z_i=n^{-\frac{1}{\alpha}}Y_i$ and $S_n=Z_1+Z_2+\ldots+Z_n$. Now, using Theorem \ref{PP1:th6}, we show $d_{W_2}(S_n,X)=O(n^{-(\frac{2-\alpha}{\alpha})})$.
\noindent
Let us first compute the terms in Remainder $R_{N,n}$.

\noindent
We have,
\begin{align*}
2\sum_{i=1}^{n}\mathbb{E}\left(|Z_{i}|\mathbf{1}_{|Z_i|>N}  \right)&=2n^{-\frac{1}{\alpha}}\left(\displaystyle
\int_{n^{\frac{1}{\alpha}}N}^{\infty} xp(x)dx+\displaystyle
\int_{-\infty}^{n^{\frac{1}{\alpha}}N} xp(x)dx \right)\\
&=\frac{4}{\alpha-1}N^{1-\alpha}=D_{0}N^{1-\alpha}.
\end{align*}

\noindent
We also have,
\begin{align*}
\frac{C_{\alpha,m}}{n}\sum_{i=1}^{n}\mathbb{E}|Z_{i}|^{2-\alpha}&=\frac{C_{\alpha,m}}{n}\sum_{i=1}^{n}\displaystyle\int_{|x|>1}(n^{-\frac{1}{\alpha}}|x|)^{2-\alpha}p(x)dx\\
&=\frac{C_{\alpha,m}}{\alpha-1}n^{-\frac{2-\alpha}{\alpha}}=D_1n^{-\frac{2-\alpha}{\alpha}},
\end{align*}
hence $$R_{N,n}=D_{0}N^{1-\alpha}+D_1n^{-\frac{2-\alpha}{\alpha}}.$$
\noindent
For any $N>0,$ we have
\begin{align*}
\nonumber K_{\nu_{\alpha}}(t,N)&=\mathbf{1}_{[0,N]}(t)\displaystyle\int_{t}^{N}u\nu_{\alpha}(du)+\mathbf{1}_{[-N,0]}(t)\displaystyle\int_{-N}^{t}(-u)\nu_{\alpha}(du)\\
\nonumber &=\frac{m}{1-\alpha}\left(N^{1-\alpha}-t^{1-\alpha} \right)+\frac{m}{1-\alpha}\left(N^{1-\alpha}-(-t)^{1-\alpha}  \right)\\
&=\frac{m}{\alpha -1}\left(|t|^{1-\alpha}-N^{1-\alpha}  \right).
\end{align*}

\noindent
Using symmetry property of Pareto distribution, we have
\begin{align*}
K_{i}(t,N)=\frac{\alpha}{2n(\alpha -1)}\left( (|t|\land n^{-\frac{1}{\alpha}})^{1-\alpha}-N^{1-\alpha
} \right)
\end{align*}
\noindent
Then by Theorem \ref{PP1:th6}, we obtain 
\begin{align*}
d_{W_2}(S_n,X)\leq D_0N^{1-\alpha}+D_1 n^{-(\frac{2-\alpha}{\alpha})}+D_2 n^{-(\frac{2-\alpha}{\alpha})},
\end{align*}

\noindent
where $D_0,D_1$ and $D_2$ are positive constants. Since $N$ is arbitrary, let $N\to \infty.$ Hence, $d_{W_2}(S_n,X)=O(n^{-(\frac{2-\alpha}{\alpha})}).$	By \cite[Lemma A.4]{k1}, we have 
\begin{align}\label{PP1:cg}
d_{W_{1}}(S_n,X) &\leq D_{3} \sqrt{n^{-(\frac{2-\alpha}{\alpha})}}=D_{3}n^{-(\frac{1}{\alpha}-\frac{1}{2})},
\end{align}
\noindent
where $D_{3}$ is an another positive constant. Moreover, $d_{W_{1}}(S_{n},X)=O(n^{-(\frac{1}{\alpha}-\frac{1}{2})})$.

\noindent
Let us compare our result with the known results in literature.

\noindent
Johnson and Samworth \cite{john} show a convergence rate for $\alpha$-stable approximations in the Mallows distance $d_r$ with some $r>0$. They have shown that $S_n=n^{-\frac{1}{\alpha}}\sum_{i=1}^{n}Y_i$ converges to an $\alpha$-stable distribution with a rate $n^{(\frac{1}{\beta}-\frac{1}{\alpha})}$ in the distance $d_\beta$ for some $\beta\in (\alpha,2]$. Hence, at $\beta=2$, they show that the convergence  rate is at most $n^{\frac{1}{2}-\frac{1}{\alpha}}$.\vspace{1mm}
\noindent
Xu \cite{k14} proved that $S_n$ converges to an symmetric $\alpha$-stable distribution with a rate $n^{-\frac{(2-\alpha)}{\alpha}}$ in the Wasserstein-1 distance. From \cite[Example 1]{k14}, it is clear that the convergence rate $n^{\frac{1}{2}-\frac{1}{\alpha}}$ is not accessible. Note that, we obtain a rate $n^{-\frac{2-\alpha}{\alpha}}$ with $\alpha\in(1,2)$ in the $d_{W_{2}}$ distance, which is faster rate than the rate obtained in \cite{john}, whenever $\beta\in (\alpha,2)$. Observe also that, at $\beta=2,$ the rate obtained in \cite[Theorem 1.2]{john} becomes $n^{\frac{1}{2}-\frac{1}{\alpha}}$. From \eqref{PP1:cg}, it immediately follows that the rate $n^{\frac{1}{2}-\frac{1}{\alpha}}$ is accessible in the $d_{W_1}$ distance using our estimates.
\end{ex}

\section{Proofs}\label{PP1:proof}

\subsection{Proof of Theorem \ref{PP1:Th1}}
\noindent
Recall first that for $X \sim IDD(\beta,\sigma^{2},\nu)$, characteristic exponent $\eta$ is given by

\begin{equation}\label{PP1:a1}
\eta(z)=\log \phi_{X}(z)= iz\beta- \frac{\sigma^{2}z^{2}}{2} +\int_{\mathbb{R}}(e^{izu}-1-izu\mathbf{1}_{\{|u|\leq1\}}(u))\nu(du),~~z\in \mathbb{R}. 
\end{equation}

\noindent
Differentiating \eqref{PP1:a1} with respect to $z$, we have

%
%

\begin{equation}\label{PP1:a2}
\phi^{\prime}_{X}(z)=\left(i\beta-\sigma^{2}z+i\int_{\mathbb R}u(e^{izu}-\mathbf{1}_{\{|u|\leq1\}}(u))\nu(du)\right)\phi_X(z).
\end{equation}

\noindent
Recall from Section 2, $F_X$ is the distribution function (cumulative distribution function) of $X$ and if $\phi^{\prime}_{X}$ exists on $\mathbb{R}$, then, 
\begin{equation}\label{PP1:a3}
\phi_{X}(z)=\int_{\mathbb R}e^{izx}F_{X}(dx) ~~\text{and}~~\phi^{\prime}_{X}(z)=i\int_{\mathbb R}xe^{izx}F_{X}(dx),~~z\in\mathbb{R}.
\end{equation}
\noindent 
Using (\ref{PP1:a3}) in (\ref{PP1:a2}) and rearranging the integrals, we have
\begin{align}
\nonumber 0&=i\displaystyle\int_{\mathbb R}xe^{izx}F_{X}(dx)-\left(i\left(\beta +\int_{\mathbb R}u(e^{izu}-\mathbf{1}_{\{|u|\leq 1 \}}(u))\nu(du) \right)\phi_{X}(z)-\sigma^{2}z\phi_{X}(z)\right)\\
\nonumber &=i\left(\displaystyle\int_{\mathbb R}(x-\beta)e^{izx}F_{X}(dx)- \left(\int_{\mathbb R}ue^{izu}\nu(du)\right)\phi_{X}(z)+ \left(\int_{\mathbb R}u\mathbf{1}_{\{|u|\leq 1 \}}(u)\nu(du) \right)\phi_X(z)-i\sigma^{2}z\phi_X(z)\right)\\
&=\left(\displaystyle\int_{\mathbb R}(x-\beta)e^{izx}F_{X}(dx)- \left(\int_{\mathbb R}ue^{izu}\nu(du)\right)\phi_{X}(z)+ \left(\int_{\mathbb R}u\mathbf{1}_{\{|u|\leq 1 \}}(u)\nu(du) \right)\phi_{X}(z)-iz\sigma^{2}\phi_{X}(z)\right)\label{PP1:a4}
\end{align}

\noindent
The second integral of \eqref{PP1:a4} can be written as
\begin{align}
\nonumber \left(\int_{\mathbb R}ue^{izu}\nu(du)\right)\phi_X(z)&=\int_{\mathbb R}\int_{\mathbb R}ue^{izu}e^{izx}F_{X}(dx)\nu(du)\\
\nonumber &=\int_{\mathbb R}\int_{\mathbb R}ue^{iz(u+x)}\nu(du)F_{X}(dx)\\
\nonumber &=\int_{\mathbb R}\int_{\mathbb R}ue^{izy}\nu(du)F_{X}(d(y-u))\\
\nonumber &=\int_{\mathbb R}\int_{\mathbb R}ue^{izx}\nu(du)F_{X}(d(x-u))\\ 
&=\int_{\mathbb R}e^{izx}\int_{\mathbb R}uF_{X}(d(x-u))\nu(du).\label{PP1:a5}
\end{align}

\noindent
Substituting \eqref{PP1:a5} in \eqref{PP1:a4}, we have 
\begin{align}
\nonumber 0&=\left(\displaystyle\int_{\mathbb R}(x-\beta)e^{izx}F_{X}(dx)-\int_{\mathbb R}uF_{X}(d(x-u))\nu(du) + \left(\int_{\mathbb R}u\mathbf{1}_{\{|u|\leq 1 \}}(u)\nu(du) \right)\phi_X(z)-iz\sigma^{2}\phi_X(z)\right)\\
\label{PP1:a6}&=\displaystyle\int_{\mathbb R}e^{izx}\left((x-\beta)F_{X}(dx)-\int_{\mathbb R}uF_{X}(d(x-u))\nu(du) + \left(\int_{\mathbb R}u\mathbf{1}_{\{|u|\leq 1 \}}(u)\nu(du) \right)F_X(dx)-iz\sigma^{2}F_X(dx)\right)
\end{align}

\noindent
On applying Fourier transform to \eqref{PP1:a6}, multiplying with $g\in \mathcal{S}(\mathbb{R}),$ and integrating over $\mathbb{R},$ we get
\begin{align}\label{PP1:a7}
\displaystyle\int_{\mathbb R}g(x)\left((x-\beta)+ \int_{\mathbb R}u\mathbf{1}_{\{|u|\leq 1 \}}(u)\nu(du) \right)F_X(dx) -\int_{\mathbb R}ug(x)F_{X}(d(x-u))\nu(du) -\sigma^{2}\int_{\mathbb{R}}g^{\prime}(x)F_{X}(dx)=0.
\end{align}
\noindent
The third integral of \eqref{PP1:a7} can be seen as
\begin{align}
\nonumber \int_{\mathbb R}\int_{\mathbb R}ug(x)F_{X}(d(x-u))\nu(du)&=\int_{\mathbb R}\int_{\mathbb R}ug(y+u)F_{X}(dy)\nu(du)\\ 
\nonumber &=\int_{\mathbb R}\int_{\mathbb R}ug(x+u)F_{X}(dx)\nu(du)\\
&=\mathbb{E}\left(\int_{\mathbb R}ug(X+u)\nu(du)\right).\label{PP1:a8}
\end{align}
\noindent 
Substituting \eqref{PP1:a8} in \eqref{PP1:a7}, we have
\begin{equation*}
\mathbb{E}\left((X-\beta)g(X)-\int_{\mathbb R}u(g(X+u)-g(X)\mathbf{1}_{\{|u|\leq 1 \}}(u))\nu(du)-\sigma^{2}g^{\prime}(X)\right)=0.
\end{equation*}

\noindent
This proves the theorem.


\subsection{Proof of Theorem \ref{PP1:Th2}}
\noindent
Recall first that, for $X\sim\mathcal{S}(\alpha,\beta,m_1,m_2)$, characteristic exponent $\eta_{\alpha}$ is given by

\begin{equation*}
\eta_{\alpha}(z)=\log \phi_{\alpha}(z)= iz\beta +\int_{\mathbb{R}}(e^{izu}-1-izu\mathbf{1}_{\{|u|\leq1\}}(u))\nu_{\alpha}(du),~~z\in \mathbb{R}. 
\end{equation*}

%
%
%

\noindent
Following similar steps to proof of Theorem \ref{PP1:Th1}, we get the result for $\alpha \in (0,1)\cup (1,2).$\vspace{2mm}

\noindent
For $\alpha=1,$ as the characteristic exponent $\eta_{1}$ is not differentiable over $\mathbb{R}$. Let us consider tempered $1$-stable random variable $Y_{\gamma}$ with characteristic exponent (see, Section 2) given by
\begin{equation*}
\eta_{1,\gamma}(z)= iz\beta+\int_{\mathbb R}(e^{izu}-1-izu\mathbf{1} _{\{|u|\leq1\}}(u))\nu_{1,\gamma}(du),~~z \in \mathbb{R},
\end{equation*}
\noindent
where $\gamma\in(0,\infty),$ and $\nu_{1,\gamma}$ is the L\'evy measure defined as

$$	\nu_{1,\gamma}(du):=\left(m_{1}\frac{e^{-\gamma u}}{u^{2}}\mathbf{1}_{(0,\infty)}(u)+m_{2}\frac{e^{-\gamma |u|}}{|u|^{2}}\mathbf{1}_{(-\infty,0)}(u)\right) du.
$$

\noindent
Observe that $Y_\gamma$ is infinitely divisible and its characteristic exponent $\eta_{1,\gamma}$ is differentiable on $\mathbb{R}$. Also, it can be easily shown that as $\gamma\to 0^{+}$, $\eta_{1,\gamma}\to \eta_{1}$, the characteristic exponent of $1$-stable random variable $X$.

\noindent
Now, applying Theorem \ref{PP1:Th1}, we get the Stein identity for $Y_\gamma$ as follows.
\begin{equation}\label{PP1:a14}
\mathbb{E}\left[\left(Y_\gamma-\beta\right)g(Y_\gamma)-\int_{\mathbb R}\left(g(Y_\gamma+u)-g(Y_\gamma)\mathbf{1}_{\{|u|\leq1\}}(u)\right)u\nu_{1,\gamma}(du) \right]=0, ~~ g\in\mathcal{S}(\mathbb{R}) .
\end{equation}

\noindent
Now, taking limit as $\gamma\to 0^{+}$, \eqref{PP1:a14} reduces to

$$\mathbb{E}\left[\left(X-\beta\right)g(X)-\int_{\mathbb R}\left(g(X+u)-g(X)\mathbf{1}_{\{|u|\leq1\}}(u)\right)u\nu_{1}(du) \right]=0, ~~  g\in\mathcal{S}(\mathbb{R}).$$

\noindent
This proves the theorem.

\subsection{Proof of Theorem \ref{PP1:ThmSol}}

\noindent
For the proof of this theorem, we use the connection between the operators $\mathcal{A}_{X}^{\alpha}$ and $\mathcal{T}_\alpha.$

\noindent
\textbf{Proof of} $\mathbf{(i).}$ The proof of this part is split into two parts.

\noindent
\textbf{(a) $\alpha\in (0,1)$:} We have
\begin{align}
\nonumber\mathcal{A}_{X}^{\alpha}g_{h}^{\alpha}(x)&=(-x+\beta)g_{h}^{\alpha}(x)+\int_{\mathbb{R}}\left( g_{h}^{\alpha}(x+u)-g_{h}^{\alpha}(x)\mathbf{1}_{\{|u|\leq 1 \}}(u) \right)u\nu_{\alpha}(du)\\
\nonumber&=\mathcal{T}_{\alpha}(\tilde{g}_{h}^{\alpha})(x),~~(\text{where }\tilde{g}_{h}^{\alpha}(x)=-\displaystyle\int_{0}^{\infty}\left(P_{t}^{\alpha}(h)(x)-\mathbb{E}h(X)  \right)dt,h\in \mathcal{H}_{\delta},\delta\in(0,\alpha)) \\
\nonumber&=-\displaystyle\int_{0}^{\infty}\mathcal{T}_{\alpha}P_{t}^{\alpha}(h)(x)dt\\
\nonumber&=-\displaystyle\int_{0}^{\infty}\frac{d}{dt}P_{t}^{\alpha}(h)(x)dt\\
\nonumber&=P_{0}h(x)-P_{\infty}h(x)\\
\nonumber&=h(x)-\mathbb{E}h(X)~~(\text{by Proposition \ref{PP1:proSem}}).
\end{align}
\noindent
Hence, $g_{h}^{\alpha}$ is the solution to \eqref{PP1:steinequ1}. Now, it remains to show that, $g_{h}^{\alpha}$ is well-defined. Let us first consider $\tilde{g}_{h}^{\alpha}: \mathbb{R}\to \mathbb{R}$ be defined as
\begin{equation*}
\tilde{g}_{h}^{\alpha}(x)=-\displaystyle\int_{0}^{\infty}\left(P_{t}^{\alpha}(h)(x)-\mathbb{E}h(X)  \right)dt,~~h\in \mathcal{H}_{\delta}, \delta\in(0,\alpha),
\end{equation*}
\noindent
where $P_{t}^{\alpha}$ is the semigroup as defined in \eqref{PP1:a20}. We show that, for any $h\in \mathcal{H}_{\delta}$, $0<\delta<\alpha$, $ \tilde{g}^{\alpha}_{h}$ is well-defined.

\noindent
Using \eqref{PP1:a21}, we have
\begin{align}
\nonumber |P^{\alpha}_{t}(h)(x)-\mathbb{E}h(X)|&=\left|\int_{\mathbb R}h(r+e^{-t}x)F_{X_{(t)}}(dr)-\int_{\mathbb R}h(r)F_{X}(dr) \right|\\
\nonumber&=\left|\int_{\mathbb R}(h(r+e^{-t}x)-h(r))F_{X_{(t)}}(dr)+\int_{\mathbb{R}}h(r)F_{X_{(t)}}(dr)-\int_{\mathbb R}h(r)F_{X}(dr) \right|\\
\nonumber &\leq \min\left\{e^{-t}|x| ,(e^{-t}|x|)^{\delta} \right\}+\left|\int_{\mathbb{R}}\widehat{h}(z)\left(\phi_{t}(z)-\phi_{\alpha}(z)  \right)dz\right|\\
&\leq \min\left\{e^{-t}|x| ,(e^{-t}|x|)^{\delta} \right\}+\int_{\mathbb{R}}|\widehat{h}(z)||\phi_{t}(z)-\phi_{\alpha}(z) |dz  \label{PP1:well0}
\end{align}
\noindent
Now, let us calculate an upper bound between the difference of two characteristic functions $\phi_{t}$ and $\phi_{\alpha}$. For all $t>0$ and $z\in \mathbb{R},$
\begin{align*}
|\phi_{t}(z)-\phi_{\alpha}(z)|=\left|\frac{\phi_{\alpha}(z)}{\phi_{\alpha}(e^{-t}z)}-\phi_{\alpha}(z)\right|\leq \left|\phi_{\alpha}(e^{-t}z)-1   \right|=|e^{\omega_{t}(z)}-1|,
\end{align*}
\noindent
where $\omega_{t}(z)=iz\beta_{1}e^{-t}+e^{-t\alpha}\int_{\mathbb{R}}(e^{izu}-1)\nu_{\alpha}(du)$, and $\beta_1=\beta-\int_{\{|u|\leq 1 \}}u\nu_\alpha(du)$. Now, from \cite[Lemma 7.9]{k16}, the function $z\to e^{s\omega_{t}(z)}$ is a characteristic function for all $s\in (0,\infty)$, thus, for all $z\in \mathbb{R}$ and $t>0$

\begin{align}
\nonumber |\phi_{t}(z)-\phi_{\alpha}(z)|&\leq \left| \int_{0}^{1}\frac{d}{ds}(\exp (s\omega_{t}(z))) ds  \right|\\
\nonumber &\leq |\omega_{t}(z)|\\
\nonumber &\leq |z||\beta_{1}|e^{-t}+e^{-t\alpha}\left|\int_{\mathbb{R}}(e^{izu}-1)\nu_{\alpha}(du) \right|\\
\nonumber &= |z||\beta_{1}|e^{-t}+e^{-t\alpha}\left| \int_{\{|u|>1 \}}(e^{izu}-1)\nu_{\alpha}(du)+\int_{\{|u|\leq 1\}}(e^{izu}-1)\nu_{\alpha}(du) \right|\\
\nonumber&\leq  |z||\beta_{1}|e^{-t}+2e^{-t\alpha}\left|\int_{\{|u|> 1 \}}\nu_{\alpha}(du)\right|\\
\nonumber&+e^{-t\alpha}\left(\left|\int_{\{|u|\leq 1 \}}(\cos(zu)-1 )\nu_{\alpha}(du)\right|+\left|\int_{\{|u|\leq 1 \}}\sin(zu)\nu_{\alpha}(du)\right|\right)\\
&=|z||\beta_{1}|e^{-t}+2\frac{m_1+m_2}{\alpha}e^{-t\alpha}+|z|^{\alpha}\left(M_1+M_2 \right)e^{-t\alpha},\label{PP1:well1}
\end{align}
\noindent
where $M_1=\left|\int_{\{|u|\leq z \}}(\cos v -1)\nu_{\alpha}(dv) \right|$ and $M_2=\left|\int_{\{|u|\leq z \}}\sin v \nu_{\alpha}(dv) \right|$.

\noindent
 Using \eqref{PP1:well1} in \eqref{PP1:well0}, one can easily show that $\displaystyle\int_{0}^{\infty}|P_t^{\alpha}(h)(x)-\mathbb{E}h(X)|dt<\infty$. Hence, $\tilde{g}^{\alpha}_{h}$ is well-defined.

\noindent
By dominated convergence theorem, we see that $\tilde{g}_{h}^{\alpha}$ is differentiable and
\begin{align}
\nonumber(\tilde{g}_{h}^{\alpha})^{\prime}(x)&=-\lim_{\zeta\to\infty} \frac{d}{dx}\displaystyle\int_{0}^{\zeta}(P_t^{\alpha}(h)(x)-\mathbb{E}h(X))dt\\
\nonumber&=-\lim_{\zeta\to\infty}\displaystyle\int_{0}^{\zeta}\frac{d}{dx}\left(\int_{\mathbb{R}}h(xe^{-t}+u)F_{X_{(t)}}(du)  \right)dt\\
\nonumber&=-\displaystyle\int_{0}^{\infty}e^{-t}\int_{\mathbb{R}}h^{\prime}(u+xe^{-t})F_{X_{(t)}}(du)dt=g_{h}^{\alpha}(x),
\end{align}

\noindent
 the desired conclusion follows.

\noindent
\textbf{(b) $\alpha=1$:} To solve \eqref{PP1:steinequ1} for $\alpha=1$, consider a Stein equation (see, \eqref{PP1:a14}) for tempered 1-stable random variable $Y_\gamma$ is given by
\begin{align}\label{PP1:a41}
\left(-x+\beta\right)g(x)+\displaystyle\int_{\mathbb R}\left(g(x+u)-g(x)\mathbf{1}_{\{|u|\leq1\}}(u)\right)u\nu_{1,\gamma}(du)=h(x)-\mathbb{E}h(Y_\gamma),~~h\in\mathcal{H}_{\delta}.
\end{align}
\noindent
Following similar steps to proof as the Case 1 of \text{(i)}, and using \eqref{PP1:a26}, we see that the function $g_{h}^{1,\gamma}(x)=-\displaystyle\int_{0}^{\infty}e^{-t}\int_{\mathbb{R}}h^{\prime}(u+xe^{-t})F_{Y_{(\gamma,t)}}(du)dt$ solves \eqref{PP1:a41} i.e.
\begin{align}\label{PP1:a42}
\left(-x+\beta\right)g_{h}^{1,\gamma}(x)+\displaystyle\int_{\mathbb R}\left(g_{h}^{1,\gamma}(x+u)-g_{h}^{1,\gamma}(x)\mathbf{1}_{\{|u|\leq1\}}(u)\right)u\nu_{1,\gamma}(du)=h(x)-\mathbb{E}h(Y_\gamma).
\end{align} 
\noindent
Observe that,
\begin{align*}
\lim_{\gamma \to 0^{+}}g_{h}^{1,\gamma}(x)&=-\lim_{\gamma \to 0^{+}}\displaystyle\int_{0}^{\infty}e^{-t}\int_{\mathbb{R}}h^{\prime}(u+xe^{-t})F_{Y_{(\gamma,t)}}(du)dt\\
&=-\displaystyle\int_{0}^{\infty}e^{-t}\int_{\mathbb{R}}h^{\prime}(u+xe^{-t})F_{X_{(t)}}(du)dt\\
&=g_h^{1}(x).
\end{align*}
\noindent
Also,
\begin{align*}
\lim_{\gamma \to 0^{+}}\mathbb{E}h(Y_\gamma)&=\mathbb{E}h(X),~~ (\text{since }Y_\gamma\overset{d}{\to}X).
\end{align*}

\noindent
Hence taking limit as $\gamma\to0^{+}$ on \eqref{PP1:a42}, we get
\begin{align*}
\left(-x+\beta\right)g_{h}^{1}(x)+\displaystyle\int_{\mathbb R}\left(g_{h}^{1}(x+u)-g_{h}^{1}(x)\mathbf{1}_{\{|u|\leq1\}}(u)\right)u\nu_{1}(du)=h(x)-\mathbb{E}h(X).
\end{align*} 
\noindent
Hence, $g_{h}^{1}$ is the solution to \eqref{PP1:steinequ1} when $\alpha=1$. Note here that, on careful adjustments of the integrals and suitably adjustments of parameters as previous case, one can verify that the function $\tilde{g}_{h}^{1}(x)=-\displaystyle\int_{0}^{\infty}\left(P_{t}^{1}(h)(x)-\mathbb{E}h(X)  \right)dt,~~h\in \mathcal{H}_{\delta},\delta\in(0,1)$ is well-defined and $(\tilde{g}^{1}_{h})^{\prime}(x)=g^{1}_{h}(x)$ for all $x\in\mathbb{{R}}.$\vspace{1mm}

\noindent
$\mathbf{\textbf{Proof of }(ii).}$ Following similar steps to proof of Case 1 of \textbf{(i)}, it immediately shows that $g^{\alpha}_{h}$ (where $\alpha\in(1,2)$ and $h\in\mathcal{H}_2$) is the solution to \eqref{PP1:steinequ2}. So, it remains to show that $g_{h}^{\alpha}$ is well-defined. Let us consider a function $\tilde{g}_{h}^{\alpha}: \mathbb{R}\to \mathbb{R}$ defined as
\begin{equation*}
\tilde{g}_{h}^{\alpha}(x)=-\displaystyle\int_{0}^{\infty}\left(P_{t}^{\alpha}(h)(x)-\mathbb{E}h(X)  \right)dt,~~h\in \mathcal{H}_{2},
\end{equation*}
\noindent
where $P_{t}^{\alpha}$ is the semigroup as defined in \eqref{PP1:a20}. Now, we show that for any $h\in \mathcal{H}_{2}$ and $\alpha\in(1,2)$, $ \tilde{g}^{\alpha}_{h}$ is well-defined and $(\tilde{g}^{\alpha}_{h})^{\prime}(x)=g^{\alpha}_{h}(x)$ for all $x\in\mathbb{{R}}$.

%

\noindent
Using \eqref{PP1:a21}, we have
\begin{align}
\nonumber |P^{\alpha}_{t}(h)(x)-\mathbb{E}h(X)|&=\left|\int_{\mathbb R}h(r+e^{-t}x)F_{X_{(t)}}(dr)-\int_{\mathbb R}h(r)F_{X}(dr) \right|\\
\nonumber&=\left|\int_{\mathbb R}(h(r+e^{-t}x)-h(r))F_{X_{(t)}}(dr)+\int_{\mathbb{R}}h(r)F_{X_{(t)}}(dr)-\int_{\mathbb R}h(r)F_{X}(dr) \right|\\
\nonumber &\leq e^{-t}|x| |h^{\prime}|+\left|\int_{\mathbb{R}}\widehat{h}(z)\left(\phi_{t}(z)-\phi_{\alpha}(z)  \right)dz\right|\\
&\leq e^{-t}|x| |h^{\prime}|+\int_{\mathbb{R}}|\widehat{h}(z)||\phi_{t}(z)-\phi_{\alpha}(z) |dz  \label{PP1:ell00}
\end{align}
\noindent
Now, let us calculate an upper bound between the difference of two characteristic functions $\phi_{t}$ and $\phi_{\alpha}$. For all $t>0$ and $z\in \mathbb{R},$
\begin{align*}
|\phi_{t}(z)-\phi_{\alpha}(z)|=\left|\frac{\phi_{\alpha}(z)}{\phi_{\alpha}(e^{-t}z)}-\phi_{1}(z)\right|\leq \left|\phi_{\alpha}(e^{-t}z)-1   \right|=|e^{\omega_{t}(z)}-1|,
\end{align*}
\noindent
where $\omega_{t}(z)=e^{-t\alpha}(iz\tilde{\beta}+\int_{\mathbb{R}}(e^{izu}-1-izu\mathbf{1}_{\{|u|\leq1 \}})\nu_{\alpha}(du))$, $\tilde{\beta}=\beta e^{(\alpha-1)t}+\int_{\mathbb{R}}(u\mathbf{1}_{\{|u|\leq1 \}}-u\mathbf{1}_{\{|u|\leq e^{-t} \}})\nu_{\alpha}(du).$ Note that the function $z\to e^{s\omega_{t}(z)}$ is a characteristic function for all $s\in (0,\infty)$. Indeed, $e^{s\omega_{t}(z)}$ is a characteristic function of an $\alpha$-stable random variable with different parameters.

\noindent
Thus, for all $z\in \mathbb{R}$ and $t>0$,
\begin{align}
\nonumber |\phi_{t}(z)-\phi_{\alpha}(z)|&\leq \left| \int_{0}^{1}\frac{d}{ds}(\exp (s\omega_{t}(z))) ds  \right|\\
\nonumber &\leq |\omega_{t}(z)|\\
\nonumber &\leq e^{-t\alpha}\left|iz\tilde{\beta}+\int_{\mathbb{R}}(e^{izu}-1-izu\mathbf{1}_{\{|u|\leq 1 \}})\nu_{\alpha}(du)  \right|\\
& \leq C_\alpha e^{-t\alpha}(1+|z|^{2}),~~C_\alpha>0, \label{PP1:ell01}
\end{align}
\noindent
where the last inequality is followed by \cite[p.30, Ex. 1.2.16]{k27}. Using \eqref{PP1:ell01} in \eqref{PP1:ell00}, one can easily show that $\displaystyle\int_{0}^{\infty}|P_t^{\alpha}(h)(x)-\mathbb{E}h(X)|dt<\infty$. Hence, $\tilde{g}^{\alpha}_{h}(x)$ is well-defined. The rest of this part follows from similar computations as Case 1 of \textbf{(i)}. 

\subsection{Proof of Theorem \ref{PP1:Th4}}
\noindent
Recall the definition of $(P_{t}^{\alpha})_{t\geq0},$ 
\begin{equation*}
P_{t}^{\alpha}(g)(x)=\int_{\mathbb R}g(r+e^{-t}x)F_{X_{(t)}}(dr),~~g\in\mathcal{F},
\end{equation*} 
where $\alpha\in (0,2)$ and $F_{X_{(t)}}$ is the distribution function of $X_{(t)}$ (see, \eqref{PP1:a19}).

\noindent
\textbf{Proof of (i).} Suppose $\alpha\in (0,1)$ and $h\in \mathcal{H}_\delta,$ $\delta\in(0,\alpha)$.

\noindent
Let \begin{equation*}
g^{\alpha}_{h}(x)=
-\displaystyle\int_{0}^{\infty}e^{-t}\int_{\mathbb R}h^{\prime}(xe^{-t}+u)F_{X_{(t)}}(du)dt .
\end{equation*}

\noindent
It is clear that
\begin{align*}
\|g^{\alpha}_{h}\| &\leq \left|\int_{0}^{\infty}e^{-t}dt|| \int_{\mathbb{R}}\|h^{\prime}\|F_{X_{(t)}}(du)\right|\\
&=\|h^{\prime}\|,
\end{align*}
the desired conclusion follows.

\noindent
Now observe that, for any $x,y\in\mathbb{R}$ and $h\in \mathcal{H}_{\delta}$,

\begin{align*}
\left|g^{\alpha}_{h}(x)-g^{\alpha}_{h}(y) \right| & \leq \displaystyle\int_{0}^{\infty}e^{-t}\int_{\mathbb{R}}\left|h^{\prime}(xe^{-t}+z)-h^{\prime}(ye^{-t}+z)  \right|F_{X_{(t)}}(dz)dt\\
&\leq \displaystyle\int_{0}^{\infty}e^{-t}\int_{\mathbb{R}}\left|x-y\right|^{\delta}e^{-t\delta}F_{X_{(t)}}(dz)dt\\
&=\left|x-y  \right|^{\delta}\displaystyle\int_{0}^{\infty}e^{-t(1+\delta)}dt
\end{align*}
\begin{align*}
\left|g^{\alpha}_{h}(x)-g^{\alpha}_{h}(y) \right|&\leq\frac{1}{1+\delta}\left|x-y \right|^{\delta},
\end{align*}
\noindent
the desired conclusion follows.

\noindent
For $\alpha\in (0,1),$ we have

\begin{align*}
\left|\displaystyle\int_{\mathbb{R}}ug^{\alpha}_{h}(x+u)\nu_{\alpha}(du)  \right|&=\alpha\left|\int_{\mathbb{{R}}}\int_{0}^{\infty}(P_{t}^{\alpha}h(x+u)-P_{t}^{\alpha}h(x))dt\nu_{\alpha}(du)\right|\text{ (using Proposition \ref{PP1:appendixPro4})}\\
& \leq\alpha\left| \int_{|u|>1 }\int_{0}^{\infty}(P_{t}^{\alpha}h(x+u)-P_{t}^{\alpha}h(x))dt\nu_{\alpha}(du) \right|\\
&+\alpha\left|\int_{|u|\leq 1 }\int_{0}^{\infty}(P_{t}^{\alpha}h(x+u)-P_{t}^{\alpha}h(x))dt\nu_{\alpha}(du)\right |\\
&:=\text{I+II}.
\end{align*}
Now observe that
\begin{align}
\nonumber \text{I}&= \alpha\left|\displaystyle\int_{|u|>1}\int_{0}^{\infty}(P_{t}^{\alpha}h(x+u)-P_{t}^{\alpha}h(x))dt\nu_{\alpha}(du)  \right|\\
\nonumber&=\alpha\left|\displaystyle\int_{|u|>1}\int_{0}^{\infty}\int_{\mathbb{R}}\left(h((x+u)e^{-t}+y)-h(xe^{-t}+y)  \right)F_{X_{(t)}} (dy) dt \nu_{\alpha}(du)  \right|\\
\nonumber&\leq \alpha\displaystyle\int_{|u|>1}|u|^{\delta}\int_{0}^{\infty}e^{-t\delta}dt\nu_{\alpha}(du)\\
\label{PP1:a64}&=\frac{\alpha(m_1+m_2)}{\delta(\alpha- \delta)}.
\end{align}
\noindent
 Consider,
\begin{align}
\nonumber\text{II}=& \alpha\left|\displaystyle\int_{|u|\leq1}\int_{0}^{\infty}(P_{t}^{\alpha}h(x+u)-P_{t}^{\alpha}h(x))dt\nu_{\alpha}(du)  \right|\\
\nonumber&=\alpha\left|\displaystyle\int_{|u|\leq1}\int_{0}^{\infty}\int_{\mathbb{R}}\left(h((x+u)e^{-t}+y)-h(xe^{-t}+y)  \right)F_{X_{(t)}} (dy) dt \nu_{\alpha}(du)  \right|\\
\nonumber&\leq \alpha\displaystyle\int_{|u|\leq1}|u|\int_{0}^{\infty}e^{-t}dt\nu_{\alpha}(du)\\
\label{PP1:a65}
&\leq\frac{\alpha(m_1-m_2)}{1-\alpha}.
\end{align}
Hence, by \eqref{PP1:a64} and \eqref{PP1:a65}, we get

\begin{align*}
|A_0g^{\alpha}_{h}(x)| \leq \frac{\alpha(m_1+m_2)}{\delta(\alpha- \delta)}+\frac{\alpha(m_1-m_2)}{1-\alpha}:=C_{\alpha,\delta,m_1,m_2},
\end{align*}
\noindent
the desired conclusion follows.\vspace{1mm}

\noindent
\textbf{Proof of (ii).} The proofs of first two properties are similar to previous case. To prove the third property, we split  $g^{1}_{h}$ in terms of the semigroup $P^{1}_t$ defined in \eqref{PP1:a20}. We write
\begin{align*}
A_1g^{1}_{h}&:=\displaystyle\int_{\mathbb{R}}u\left(g^{1}_{h}(x+u)-g^{1}_{h}(x)\mathbf{1}_{\{|u|\leq 1 \}} \right)\nu_{1}(du)\\
&=\displaystyle\int_{\{|u|\leq 1 \}}u\left(g^{1}_{h}(x+u)-g^{1}_{h}(x) \right)\nu_{1}(du)+\displaystyle\int_{\{|u|>1\}}ug^{1}_{h}(x+u)\nu_{1}(du)\\
&=\displaystyle\int_{\{|u|\leq 1 \}}u\int_{0}^{\infty}\left(e^{-t}\left(\int_{\mathbb{R}}h^{\prime}((x+u)e^{-t}+y)-h^{\prime}(xe^{-t}+y)  \right)F_{X_{(t)}}(dy)  \right)dt\nu_{1}(du)\\
&+\displaystyle\int_{\{|u|>1\}}\int_{0}^{\infty}(P^{1}_{t}h(x+u)-P^{1}_{t}h(x))dt\nu_{1}(du) \text{ (using Fubini's theorem)}\\
&:=\text{I+II}.
\end{align*}

\noindent
By similar computation as Case 1 of \textbf{(i)}, it is easy to show that 

\begin{align}
\nonumber \text{I}&=\displaystyle\int_{\{|u|\leq 1 \}}u\int_{0}^{\infty}\left(e^{-t}\left(\int_{\mathbb{R}}h^{\prime}((x+u)e^{-t}+y)-h^{\prime}(xe^{-t}+y)  \right)F_{X_{(t)}}(dy)  \right)dt\nu_{1}(du)\\
& \leq \frac{1}{1+\delta}\displaystyle\int_{\{|u|\leq 1 \}}|u|^{1+\delta}\nu_{1}(du)=\frac{m_1+m_2}{\delta(1+\delta)}.\label{PP1:es1}
\end{align}

\noindent
Consider,
\begin{align}
\nonumber\text{II}&=\displaystyle\int_{\{|u|>1\}}\int_{0}^{\infty}(P^{1}_{t}h(x+u)-P^{1}_{t}h(x))dt\nu_{1}(du)\\
&\leq \frac{1}{\delta}\int_{\{|u|>1\}}|u|^{\delta}\nu_{1}(du)=\frac{m_1+m_2}{\delta(1-\delta)}\label{PP1:es2}
\end{align}

\noindent
Hence, by \eqref{PP1:es1} and \eqref{PP1:es2}, we get

 $$\|A_1g^{1}_{h}\| \leq \frac{2(m_1+m_2)}{\delta(1-\delta^{2})}:=C_{1,\delta,m_1,m_2},$$

\noindent
 the desired conclusion follows.\vspace{2mm}


\noindent
\textbf{Proof of (iii).} Suppose $\alpha\in (1,2)$ and $h\in \mathcal{H}_2$.

\noindent
Let \begin{equation*}
g^{\alpha}_{h}(x)=
-\displaystyle\int_{0}^{\infty}e^{-t}\int_{\mathbb R}h^{\prime}(xe^{-t}+u)F_{X_{(t)}}(du)dt .
\end{equation*}
\noindent
Then,
\begin{align*}
(g^{\alpha}_{h})^{\prime}(x)&=-\displaystyle\int_{0}^{\infty}e^{-2t}\int_{\mathbb R}h^{\prime}(xe^{-t}+u)F_{X_{(t)}}(du)dt .
\end{align*}

\noindent
It is also easy to show that
$$\|(g^{\alpha}_{h})\| \leq \|h^{\prime}\|, \text{ and }\|(g^{\alpha}_{h})^{\prime}\| \leq \frac{1}{2} \|h^{\prime\prime}\|.$$

\noindent
Let $m_1=m_2=m$. Let $A_{2}g^{\alpha}_{h}(x)=\displaystyle\int_{\mathbb{{R}}}(g^{\alpha}_{h}(x+u)-g^{\alpha}_{h}(x))u\nu_{\alpha}(du)$. Then, for any $x,y\in \mathbb{R}$
\begin{align*}
\left|A_{2}g^{\alpha}_{h}(x)-A_{2}g^{\alpha}_{h}(y) \right|&\leq  \displaystyle\int_{\mathbb{{R}}}\left|(g^{\alpha}_{h}(x+u)-g^{\alpha}_{h}(y+u))-(g^{\alpha}_{h}(x)-g^{\alpha}_{h}(y))\right||u|\nu_{\alpha}(du)\\
&=m\left(\int_{|u|>|x-y|}+\int_{-|x-y|}^{|x+y|}  \right)\left|(g^{\alpha}_{h}(x+u)-g^{\alpha}_{h}(y+u))-(g^{\alpha}_{h}(x)-g^{\alpha}_{h}(y))\right|\frac{du}{|u|^{\alpha}}\\
&=: \text{I+II}
\end{align*}
\noindent
Now observe that
\begin{align}
\nonumber\text{I}&=m\int_{|u|>|x-y|}  \left|(g^{\alpha}_{h}(x+u)-g^{\alpha}_{h}(y+u))-(g^{\alpha}_{h}(x)-g^{\alpha}_{h}(y))\right|\frac{du}{|u|^{\alpha}}\\
\nonumber&\leq m\int_{|u|>|x-y|}\left(\left| g^{\alpha}_{h}(x+u)-g^{\alpha}_{h}(y+u) \right|+|(g^{\alpha}_{h}(x)-g^{\alpha}_{h}(y))|\right)\frac{du}{|u|^{\alpha}} \\
\nonumber&\leq 4m\|(g^{\alpha}_{h})^{\prime}\| |x-y|\int_{|x-y|}^{\infty}u^{-\alpha}du\\
&\leq 2m\|h^{\prime\prime}\| \frac{|x-y|^{2-\alpha}}{\alpha-1}.\label{PP1:es3}
\end{align}
Consider,
\begin{align}
\nonumber\text{II}&=m\int_{-|x-y|}^{|x+y|}  \left|(g^{\alpha}_{h}(x+u)-g^{\alpha}_{h}(y+u))-(g^{\alpha}_{h}(x)-g^{\alpha}_{h}(y))\right|\frac{du}{|u|^{\alpha}}\\
\nonumber&\leq m\int_{-|x-y|}^{|x-y|}\left(\left| g^{\alpha}_{h}(x+u)-g^{\alpha}_{h}(x) \right|+|(g^{\alpha}_{h}(y+u)-g^{\alpha}_{h}(y))|\right)\frac{du}{|u|^{\alpha}}\\
\nonumber&\leq 4m\|(g^{\alpha}_{h})^{\prime}\|\displaystyle\int_{0}^{|x-y|}u^{1-\alpha}du\\
&\leq 2m\|h^{\prime\prime
}\| \frac{|x-y|^{2-\alpha}}{2-\alpha}.\label{PP1:es4}
\end{align}
\noindent
Hence, by \eqref{PP1:es3} and \eqref{PP1:es4}, we get
\begin{align*}
\left\|A_{2}g^{\alpha}_{h}(x)-A_{2}g^{\alpha}_{h}(y) \right\|\leq C_{\alpha,m} \|h^{\prime\prime}\||x-y|^{2-\alpha},~~\text{where }C_{\alpha,m}=2m\left(\frac{1}{2-\alpha}+\frac{1}{\alpha-1} \right).
\end{align*}

\subsection{Proof of Theorem \ref{PP1:Th5}}
\noindent
Recall that $(Y_{n})_{n\geq1}$ is a sequence of i.i.d. random variables such that $Y_{1}\in D_{\alpha}$. Let us denote
\begin{align*}
S_{n}&=Y_{1}+Y_{2}+\ldots+Y_{n},\text{ and}\\
S_{n,i}&=S_{n}-n^{-\frac{1}{\alpha}}Y_{i}.
\end{align*}
\noindent
 Note that, $S_{n,i}$ and $Y_{i}$ are independent. To prove this theorem, we first derive some lemmas. With the help of these lemmas, we obtain our bounds in the $d_{W_{\delta}}$ distance for $\alpha$-stable approximations with $\alpha\in(0,1]$. 
\subsubsection{Proof of (a).}
\noindent
 To prove this part of Theorem \ref{PP1:Th5}, we use the following lemmas. Recall the L\'evy measure $\nu_{\alpha}$ for $\alpha$-stable distributions is given by $\nu_{\alpha}(du)=\left(m_1\frac{1}{|u|^{1+\alpha}}\mathbf{1}_{(0,\infty)}(u)+m_2\frac{1}{|u|^{1+\alpha}}\mathbf{1}_{(-\infty,0)}(u)\right)du,$ where $m_1,m_2\in [0,\infty)$, $m_1+m_2>0$ and $\alpha\in (0,2)$.

\begin{lem}\label{PP1:lembd0}
	Let $\alpha\in(0,1).$ Let $g^{\alpha}_{h}$ is defined in \eqref{PP1:solution1}. Then, for any $a>0$, 
	\begin{equation*}
	\displaystyle\int_{\mathbb{R}}ug^{\alpha}_{h}(x+u)\nu_{\alpha}(du)=a^{1-\alpha}\int_{\mathbb{R}}ug^{\alpha
	}_{h}(x+au)\nu_{\alpha}(du)
	\end{equation*}
\end{lem}

\noindent
$\mathbf{Proof.}$ We write 
\begin{align*}
\displaystyle\int_{\mathbb{R}}ug^{\alpha}_{h}(x+u)\nu_{\alpha}(du)&=\int_{\mathbb{R}}ug^{\alpha
}_{h}(x+u)\frac{(m_{1}1_{(0,\infty)}(u)+m_{2}1_{(-\infty,0)}(u))}{|u|^{\alpha+1}}du\\
&=a^{1-\alpha}\int_{\mathbb{R}}ug^{\alpha
}_{h}(x+au)\frac{(m_{1}1_{(0,\infty)}(u)+m_{2}1_{(-\infty,0)}(u))}{|u|^{\alpha+1}}du\\
&=a^{1-\alpha}\int_{\mathbb{R}}ug^{\alpha
}_{h}(x+au)\nu_{\alpha}(du),
\end{align*}

\noindent
the desired conclusion follows.

\begin{lem}\label{PP1:lembd1}
	Let $\alpha\in(0,1)$. Let $Y\in D_\alpha$ and $g_{h}^{\alpha}$ is defined in \eqref{PP1:solution1}. Then, for $0<a<1$ and $z\in\mathbb{R},$
	\begin{equation*}
	\mathbb{E}\left(\left|\int_{\mathbb{R}} u\left(g_{h}^{\alpha}(z+aY+u)- g^{\alpha
	}_{h}(z+u)  \right)\nu_{\alpha}(du)\right|\right)
	\leq C_{\alpha,\delta,m_1,m_2}^{A,K}a^{\alpha}
	\end{equation*}
\end{lem}

\noindent
\textbf{Proof.} We write 

$$	\mathbb{E}\left(\left|\int_{\mathbb{R}} u\left(g_{h}^{\alpha}(z+aY+u)- g^{\alpha
}_{h}(z+u)  \right)\nu_{\alpha}(du)\right|\right)
:=\text{I+II,}
$$
where
\begin{align*}
\text{I:}&=\mathbb{E}\left(\left|\int_{\mathbb{R}} u\left(g_{h}^{\alpha}(z+aY+u)- g^{\alpha
}_{h}(z+u)  \right)\nu_{\alpha}(du)\right|\mathbf{1}_{|Y|>a^{-1}}\right),\\
\text{II:}&=\mathbb{E}\left(\left|\int_{\mathbb{R}} u\left(g_{h}^{\alpha}(z+aY+u)- g^{\alpha
}_{h}(z+u)  \right)\nu_{\alpha}(du)\right|\mathbf{1}_{|Y|\leq a^{-1}}\right).
\end{align*}
\noindent
For $\alpha\in (0,1),$ one can write by \eqref{PP1:a51} and \eqref{PP1:a68},
\begin{align}
\nonumber\text{I}&\leq 2C_{\alpha,\delta,m_1,m_2}P(|Y|\geq a^{-1})\\
\nonumber&\leq 4 C_{\alpha,\delta,m_1,m_2}\left(A+\text{sup}_{|y|\geq a^{-1}}|e(y)| \right)a^{\alpha}\\
&\leq 4C_{\alpha,\delta,m_1,m_2}(A+K)a^{\alpha}.\label{PP1:pflem0}
\end{align}
\noindent
It is also easy to show that
\begin{align}\label{PP1:pflem1}
\text{II}\leq C_{\alpha}a^{\alpha}.
\end{align}

\noindent
Hence, by \eqref{PP1:pflem0} and \eqref{PP1:pflem1}, we have
\begin{align*}
\nonumber\mathbb{E}\left(\left|\int_{\mathbb{R}} u\left(g_{h}^{\alpha}(z+aY+u)- g^{\alpha
}_{h}(z+u)  \right)\nu_{\alpha}(du)\right|\right)&\leq \left(4C_{\alpha,\delta,m_1,m_2}+C_\alpha \right)a^{\alpha}\\
&\leq C_{\alpha,\delta,m_1,m_2}^{A,K}a^{\alpha},
\end{align*}
\noindent
the desired conclusion follows.\vspace{2mm}

\noindent
Recall the definition of $D_\alpha$ in Definition \ref{PP1:def2}. We see that the function $e$ satisfies certain conditions with the domain $\{|y|>L\}$. These conditions play an important role for proving the following lemma. 
\begin{lem}\label{PP1:lembd2}
	Let $\alpha\in (0,1)$. Let $Y\in D_\alpha$ and $X$ be a random variable with finite $\delta$-th moment, which is independent of $Y$. For any $0<a<\frac{1}{L}$ and $g^{\alpha}_{h}$ defined in \eqref{PP1:solution1}, define
	\begin{align*}
	\mathcal{J}_1:=\left|\mathbb{E}\left(Yg^{\alpha}_{h}(X+aY)\right)-\mathbb{E}\left(Y\mathbf{1}_{(-1,1)}(aY)  \right)\mathbb{E}\left(g^{\alpha}_{h} (X)\right)    \right|,
	\end{align*}
	Then, 
	\begin{align*}
\nonumber	\mathcal{J}_1 &\leq C_{1,\delta,L}a^{\delta}+C_{2,\delta}a^{\delta}\sup_{L<|y|<\frac{1}{a}} \left( \alpha|e(y)|+|ye^{\prime}(y)| \right)\displaystyle\int_{L<|y|<\frac{1}{a}}|y|^{\delta -\alpha}dy\\
	&+2a^{\alpha-1}\sup_{|y|>a^{-1}}\left(\alpha |e(y)|+|ye^{\prime}(y)| \right)\displaystyle\int_{|y|> 1 }\eta_{\alpha,\beta,\delta,m_1,m_2}(y)|y|^{-1-\alpha}dy,
	\end{align*}
	where $C_{1,\delta,L}=\frac{2}{1+\delta}L^{1+\delta}$ and $C_{2,\delta}=\frac{2}{1+\delta}.$
\end{lem}
\noindent
\textbf{Proof.} We have by \eqref{PP1:a68},
\begin{align*}
\mathcal{J}_1=\left|\mathbb{E}\left(\displaystyle\int_{\mathbb{R}}(yg^{\alpha}_{h}(X+ay)-y\mathbf{1}_{(-1,1)}(ay)g^{\alpha}_{h}(X))dF_{Y}(y)   \right)\right|.
\end{align*}

\noindent
Since $e$ is in $C^{2},$ for any $|y|>L,$
$$dF_{Y}(y)=\frac{A\alpha+\alpha e(y)-ye^{\prime}(y)}{|y|^{1+\alpha}}\kappa _{\theta}(y)dy,$$
where $\kappa_{\theta}(y)=(1+\theta)\mathbf{1}_{(0,\infty)}(y)+(1-\theta)\mathbf{1}_{(-\infty,0)}(y).$

\noindent
Thus, we have
\begin{align*}
\mathcal{J}_1 &\leq \mathbb{E}\displaystyle\int_{|y|<L}|y|\left|g^{\alpha}_{h}(X+ay)-g^{\alpha}_{h}(X)\right|dF_{Y}(y)   \\
&+\mathbb{E}\displaystyle\int_{L<|y|<\frac{1}{a}}|y|\left|g^{\alpha}_{h}(X+ay)-g^{\alpha}_{h}(X)\right|dF_{Y}(y)   \\
&+\mathbb{E}\displaystyle\int_{|y|>\frac{1}{a}}|yg^{\alpha}_{h}(X+ay)|dF_{Y}(y)\\
&   :=\text{I+II+III}.
\end{align*}

\noindent
It is easy to verify by \eqref{PP1:a50},
\begin{align*}
\text{I}&\leq \frac{1}{1+\delta}a^{\delta}\displaystyle\int_{|y|<L}|y|^{1+\delta}dF_{Y}(y)\\
&\leq \frac{2}{1+\delta} L^{1+\delta} a^{\delta},
\end{align*}
\noindent
and 
\begin{align*}
\text{II}&\leq \frac{2a^{\delta}}{1+\delta}\displaystyle\int_{L<|y|<\frac{1}{a}} \frac{|\alpha e(y)-ye^{\prime}(y)|}{|y|^{\alpha -\delta}}dy\\
&\leq\frac{2a^{\delta}}{1+\delta}\displaystyle\int_{L<|y|<\frac{1}{a}}|y|^{\delta -\alpha}\left(\alpha|e(y)|+|ye^{\prime}(y)| \right)dy\\
&\leq \frac{2a^{\delta}}{1+\delta} \sup_{L<|y|<\frac{1}{a}} \left( \alpha|e(y)|+|ye^{\prime}(y)| \right)\displaystyle\int_{L<|y|<\frac{1}{a}}|y|^{\delta -\alpha}dy.
\end{align*}

\noindent
For the third term, we have,
\begin{align*}
\text{III}& \leq 2 \displaystyle\int_{|y|>\frac{1}{a}}\frac{|\alpha e(y)-ye^{\prime}(y)||yg^{\alpha}_{h}(X+ay)|}{|y|^{\alpha+1}}dy\\
&\leq 2\sup_{|y|>a^{-1}}\left(\alpha |e(y)|+|ye^{\prime}(y)| \right)\displaystyle\int_{|y|>\frac{1}{a}}\frac{|yg^{\alpha}_{h}(X+ay)|}{|y|^{\alpha+1}}dy\\
&\leq 2a^{\alpha-1}\sup_{|y|>a^{-1}}\left(\alpha |e(y)|+|ye^{\prime}(y)| \right)\displaystyle\int_{|y|> 1 }\eta_{\alpha,\beta,\delta,m_1,m_2}(y)|y|^{-1-\alpha}dy,
\end{align*}

\noindent
where the last inequality follows by Lemma \ref{PP1:lembd0} and Proposition \ref{PP1:proappendix5}. Combining the estimates obtained in I, II and III, the desired conclusion follows.\vspace{1mm}

\noindent
\textbf{Proof of (a).}  With the help of above lemmas, we now find bound in the $d_{W_\delta}$ distance for $\alpha$-stable approximation with $\alpha\in(0,1)$.

\noindent
By \eqref{PP1:a15}, we have 
$$\left|\mathbb{E}[h(S_{n})-h(X)]\right|= \left|\mathbb{E}\left[-S_{n}g^{\alpha
}_{h}(S_{n})+\beta_1 g^{\alpha
}_{h}(S_{n})+\int_{\mathbb{R}}g^{\alpha
}_{h}(S_{n}+u)u\nu_{\alpha}(du)\right]\right|\leq \text{I+II+III} $$  

\noindent
where,

\begin{align*}
\text{I}&:=\frac{1}{n}\sum_{i=1}^{n}\mathbb{E}\left|\int_{\mathbb{R}}u\left(g^{\alpha}_{h}(S_{n,i}+n^{-\frac{1}{\alpha}}Y_i+u)-g^{\alpha}_{h}(S_{n,i}+u)  \right)\nu_{\alpha}(du)  \right|\\
\text{II}&:=n^{-\frac{1}{\alpha}}\sum_{i=1}^{n}\mathbb{E}\left|-Y_{i}g^{\alpha}_{h}(S_{n})+Y_{i}\mathbf{1}_{(-1,1)}(|n^{-\frac{1}{\alpha}}Y_{i}|) \mathbb{E}g^{\alpha}_{h}(S_{n,i}) \right|\\
\text{III}&:=\frac{1}{n}\sum_{i=1}^{n}\mathbb{E}\left|\displaystyle\int_{\mathbb{R}}ug^{\alpha}_{h}(S_{n,i}+u)\nu_{\alpha}(du)+\beta_1g^{\alpha}_{h}(S_{n}) -n^{1-\frac{1}{\alpha}}Y_{i}\mathbf{1}_{(-1,1)}(|n^{-\frac{1}{\alpha}}Y_{i}|)\mathbb{E}g^{\alpha}_{h} (S_{n,i}) \right|.
\end{align*}

\noindent
For $\alpha \in (0,1),$ we have by Lemma \ref{PP1:lembd1} with $a=n^{-\frac{1}{\alpha}}$, 

\begin{align*}
\text{I} \leq \frac{C^{A,K}_{\alpha,\delta,m_1,m_2}}{n}.
\end{align*}

\noindent
By Lemma \ref{PP1:lembd2} with $a=n^{-\frac{1}{\alpha}}$, we have

\begin{align*}
\text{II} &\leq C_{1,\delta,L} n^{1-\frac{(1+\delta)}{\alpha}}+C_{2,\delta}n^{1-\frac{(1+\delta)}{\alpha}}\sup_{L<|y|<n^{\frac{1}{\alpha}}} \left( \alpha|e(y)|+|ye^{\prime}(y)| \right)\displaystyle\int_{L<|y|<n^{\frac{1}{\alpha}}}|y|^{\delta -\alpha}dy\\
&+2\sup_{|y|>n^{\frac{1}{\alpha}}}\left(\alpha |e(y)|+|ye^{\prime}(y)| \right)\displaystyle\int_{|y|> 1 }\eta_{\alpha,\beta,\delta,m_1,m_2}(y)|y|^{-1-\alpha}dy.
\end{align*}
\noindent
Using Lemma \ref{PP1:lembd0} with $a=n^{-\frac{1}{\alpha}}$, we have

\begin{align*}
\text{III} \leq C_{\alpha,\delta,m_1,m_2} n^{-\frac{(1-\alpha)}{\alpha}} +\beta_1+n^{-\frac{(1-\alpha)}{\alpha}}\displaystyle\int_{|y|<n^{\frac{1}{\alpha}}}|y|dF_{Y}(y).
\end{align*}

\noindent
Combining the estimates obtained in I, II and III, the desired conclusion follows.

\subsubsection{Proof of (b).}
 The following two lemmas play an important role for finding bound in the $d_{W_\delta}$ distance for 1-stable approximation. Recall that $A_1g^{1}_{h}(x):=\displaystyle\int_{\mathbb{R}}u(g^{1}_{h}(x+u)-ug^{1}_{h}\mathbf{1}_{\{|u|\leq 1 \}})\nu_{1}(du)$, where $g_{h}^{1}$ is defined in \eqref{PP1:solution1} and $\nu_{1}$ is the L\'evy measure (see \eqref{u1}).

\begin{lem}\label{PP1:lembd3}
	Let $\alpha=1$. Let $Y\in D_\alpha$ and $g_{h}^{1}$ is defined in \eqref{PP1:solution1}. Then, for any $0<a<1$ and $z\in \mathbb{R}$,
	\begin{align*}
	\mathbb{E} \left( \left| A_{1}g^{1}_{h}(z)-A_{1}g^{1}_{h}(z+aY) \right| \right)    \leq C_{1,\delta,m_1,m_2}^{A,K} a + 2 C_{1,\delta,m_1,m_2} \int_{0}^{\frac{1}{a}}dF_{|Y|}(y),
	\end{align*}
	where $C_{1,\delta,m_1,m_2}^{A,K}$ and $C_{1,\delta,m_1,m_2}$ are constants.
\end{lem}

\noindent
\textbf{Proof.} We write
\begin{align*}
\mathbb{E} \left( \left| A_{1}g^{1}_{h}(z)-A_{1}g^{1}_{h}(z+aY) \right| \right):=\text{I+II},
\end{align*}
\noindent
where
\begin{align*}
\text{I}&:=\mathbb{E} \left( \left| A_{1}g^{1}_{h}(z)-A_{1}g^{1}_{h}(z+aY) \right|\mathbf{1}_{|y|>\frac{1}{a}} \right),\\
\text{II}&:=\mathbb{E} \left( \left| A_{1}g^{1}_{h}(z)-A_{1}g^{1}_{h}(z+aY) \right|\mathbf{1}_{|y|\leq\frac{1}{a}} \right)
\end{align*}

\noindent
When $\alpha=1,$ one can write by \eqref{PP1:a54} and \eqref{PP1:a68},

\begin{align}
\nonumber\text{I}&\leq 2C_{1,\delta,m_1,m_2}P(|Y|\geq a^{-1})\\
\nonumber&\leq 4 C_{1,\delta,m_1,m_2}\left(A+\text{sup}_{|y|\geq a^{-1}}|e(y)| \right)a\\
\nonumber&\leq 4C_{1,\delta,m_1,m_2}(A+K)a,
\end{align}
and 
\begin{align}
\nonumber\text{II}& \leq 2 C_{1,\delta,m_1,m_2}P(|Y|<\frac{1}{a})\\
\nonumber&=C_{1,\delta,m_1,m_2} \int_{0}^{\frac{1}{a}}dF_{|Y|}(y),
\end{align}
\noindent
Combining the estimates obtained in I and II, the desired conclusion follows.
\begin{lem}\label{PP1:lembd4}
	Let $\alpha=1.$ Let $Y\in D_{\alpha}$ and $X$ be a random variable with $\delta$-th finite moment such that $X$ and $Y$ are independent. For any $0<a<\frac{1}{L},$ define
	
	\begin{align*}
	\mathcal{J}_{2}:=\left|\mathbb{E}\left(Yg^{1}_{h}(X+aY)  \right)-\mathbb{E}\left( Y\mathbf{1}_{(-1,1)}(aY)\right)\mathbb{E}\left( g^{1}_{h}(X)\right)-\mathbb{E}\left(A_{0}g^{1}_{h}(X) \right)  \right|.
	\end{align*}
	Then, 
	\begin{align*}
	\mathcal{J}_2 \leq \frac{1}{\delta +1} a^{\delta} \left(L^{2}+m_1+m_2  \right)+\frac{a^{\delta}}{1+\delta}\displaystyle\int_{L<|u|<\frac{1}{a}} \frac{| e(u)-ue^{\prime}(u)|}{|u|^{1-\delta}}du+a \displaystyle\int_{|u|>1}\frac{|e(u/a)-u/ae^{\prime}(u/a)|}{|u|}du.
	\end{align*}
\end{lem}
\noindent
\textbf{Proof.} We have 
\begin{align*}
\mathbb{E}A_1g^{1}_{h}(X)&=\mathbb{E}\displaystyle\int_{\mathbb{R}}u\left(g^{1}_{h}(X+u)-g^{1}_{h}(X)\mathbf{1}_{\{|u|\leq 1 \}}(u) \right)\nu_{1}(du)\\
&=\mathbb{E}\displaystyle\int_{\mathbb{R}}u\left(g^{1}_{h}(X+aY)-g^{1}_{h}(X)\mathbf{1}_{\{|aY|\leq 1 \}}(aY) \right)\nu_{1}(du)
\end{align*}
\noindent
and

\begin{align*}
&\mathbb{E}\left(Yg^{1}_{h}(X+aY)  \right)-\mathbb{E}\left( Y\mathbf{1}_{(-1,1)}(aY)\right)\mathbb{E}\left( g^{1}_{h}(X)\right)\\
=&\mathbb{E}\left(\displaystyle\int_{\mathbb{R}} \left( ug^{1}_{h}(X+au)-u\mathbf{1}_{(-1,1)}(au)g^{1}_{h}(X) \right)dF_{Y}(u)  \right).
\end{align*}
\noindent
Since $e$ is $C^{2}$, for any $|y|>L$

$$dF_Y(y)=\frac{A\alpha + e(y)-ye^{\prime}(y)}{|y|^{2}}\kappa_{\theta}(y)dy,$$

\noindent
where $\kappa_{\theta}(y)=(1+\theta)\mathbf{1}_{(0,\infty)}(y)+(1-\theta)\mathbf{1}_{(-\infty,0)}(y)$.

\noindent
Thus we have,
\begin{align*}
\mathcal{J}_2 &\leq \mathbb{E} \displaystyle\int_{|u|<L} \left| ug^{1}_{h}(X+au)-ug^{1}_{h}(X)  \right|\left(dF_{y}(u)+\nu_{1}(du) \right)\\
&+\mathbb{E}\displaystyle\int_{L<|u|<\frac{1}{a}} \left|ug^{1}_{h}(X+au)-ug^{1}_{h}(X)\right|\frac{|\alpha e(u)-ue^{\prime}(u)|}{|u|^{2}}du\\
&+\mathbb{E}\displaystyle\int_{|u|>\frac{1}{a}}\left|ug^{1}_{h}(X+au)  \right|\frac{|\alpha e(u)-ue^{\prime}(u)|}{|u|^{2}}du\\
&:=\text{I+II+III}
\end{align*}
\noindent
Moreover, by \eqref{PP1:a53}, it is easy to verify
\begin{align*}
\text{I} & \leq \frac{1}{\delta +1} a^{\delta} \displaystyle\int_{|u|<L}u^{1+\delta}(dF_{Y}(u)+\nu_{1}(du))\\
& \leq \frac{1}{\delta +1} a^{\delta} \left(L^{2}+m_1+m_2  \right).
\end{align*}

\noindent
Using \eqref{PP1:a53}, we also have
\begin{align*}
\text{II} \leq \frac{a^{\delta}}{1+\delta}\displaystyle\int_{L<|u|<\frac{1}{a}} \frac{| e(u)-ue^{\prime}(u)|}{|u|^{1-\delta}}du.
\end{align*}

\noindent
For the third term, using \eqref{PP1:a52}, it can be immediately shown that

\begin{align*}
\text{III} \leq a \displaystyle\int_{|u|>1}\frac{|e(u/a)-u/ae^{\prime}(u/a)|}{|u|}du.
\end{align*} 

\noindent
Combining the estimates obtained in I, II and III, the desired conclusion follows.

\paragraph{Proof of (b).} With the help of above lemmas, we now find bound in the $d_{W_\delta}$ distance for 1-stable approximation. By \eqref{PP1:a15}, we have
$$\left|\mathbb{E}[h(S_{n})-h(X)]\right|= \left|\mathbb{E}\left[(-S_{n}+\beta)g^{1
}_{h}(S_{n})+\int_{\mathbb{R}}(g^{1
}_{h}(S_{n}+u)-g^{1}_{h}(S_n)\mathbf{1}_{\{|u|\leq 1 \}}(u))u\nu_{1}(du)\right]\right|\leq \text{I+II+III} ,   $$

\noindent
where 
\begin{align*}
\text{I}&:=\frac{1}{n}\sum_{i=1}^{n}\left|\mathbb{E} A_{1}g^{1}_{h}(S_{n,i})-\mathbb{E}A_{1}g^{1}_{h}(S_n) \right|\\
\text{II}&:=\frac{1}{n}\sum_{i=1}^{n}\left|\mathbb{E}\left( Y_ig^{1}_{h}(S_{n,i}+\frac{1}{n}Y_i)\right)- \mathbb{E}\left(Y_{i}\mathbf{1}_{(-1,1)}(|\frac{1}{n}Y_i|)g^{1}_{h}(S_{n,i}) \right)- \mathbb{E}\left(A_{1}g^{1}_{h}(S_{n,i})  \right)  \right|\\
\text{III}&:=\frac{1}{n}\sum_{i=1}^{n}\left|\mathbb{E}\left(Y_{i}\mathbf{1}_{(-1,1)}(|\frac{1}{n}Y_{i} |)g^{1}_{h}(S_{n,i})  \right) -\beta\mathbb{E}g^{1}_{h}(S_n) \right|
\end{align*}
\noindent
For $\alpha=1,$ we have by Lemma \ref{PP1:lembd3} with $a=\frac{1}{n}$,

\begin{align*}
\text{I}\leq C_{1,\delta,m_1,m_2}^{A,K} \frac{1}{n} + 2 C_{1,\delta,m_1,m_2} \int_{0}^{n}dF_{|Y|}(y)
\end{align*}

\noindent
By Lemma \ref{PP1:lembd4} with $a=\frac{1}{n}$, we have

\begin{align*}
\text{II} \leq \frac{1}{\delta +1} n^{-\delta} \left(L^{2}+m_1+m_2  \right)+\frac{n^{-\delta}}{1+\delta}\displaystyle\int_{L<|u|<\frac{1}{a}} \frac{| e(u)-ue^{\prime}(u)|}{|u|^{1-\delta}}du+\frac{1}{n} \displaystyle\int_{|u|>1}\frac{|e(nu)-nue^{\prime}(nu)|}{|u|}du .
\end{align*}

\noindent
Using \eqref{PP1:a68} and \eqref{PP1:a52}, we have

\begin{align*}
\text{III} \leq \left|\displaystyle\int_{0}^{n}\frac{e(y)-e(-y)}{y}dy \right|+2K+\beta.
\end{align*}

\noindent
Combining the estimates obtained in I, II and III, the desired conclusion follows.

\subsection{Proof of Theorem \ref{PP1:th6}}

Recall that $( Y_{n})_{n\geq1}$ is a sequence of i.i.d. random variables with $\mathbb{E}Y_{i}=0$  and $\mathbb{E}|Y_{i}|<\infty$ for $1\leq i\leq n$. Let $Z_i=n^{-\frac{1}{\alpha}}Y_i$ and define,

\begin{align*}
S_{n}&=Z_{1}+Z_{2}+\ldots+Z_{n} \text{ and}\\
S_{n}(i)&=S_{n}-Z_{i}.
\end{align*}
\noindent
Note that $S_{n}(i)$ and $S_{n}$ are independent. To derive an error bound in the $d_{W_{2}}$ distance for $\alpha$-stable approximations with $\alpha\in(1,2)$, we need to go through three important lemmas.

\begin{lem}\label{l1}
	Let $\nu_{\alpha}$ be a L$\acute{e}$vy measure for $\alpha$-stable distributions with $\alpha\in(1,2)$. Let $g^{\alpha}_{h}$ is defined in \eqref{PP1:solution2}. Then for any $N>0$, 
	\begin{equation*}
	\int_{\mathbb{R}}\big(g^{\alpha
	}_{h}(S_{n}+u)-g^{\alpha
	}_{h}(S_{n})\big)u\nu_{\alpha}(du)=\int_{-N}^{N}K_{\nu_{\alpha}}(t,N)(g^{\alpha}_{h})^{\prime}(S_{n}+t)dt+ R_{N}(S_{n}),
	\end{equation*}
	where
	\begin{align*}
	\begin{split}
	K_{\nu_{\alpha}}(t,N)&=\mathbf{1}_{[0,N]}(t)\int_{t}^{N}u\nu_{\alpha}(du)+\mathbf{1}_{[-N,0]}(t)\int_{-N}^{t}(-u)\nu_{\alpha}(du), ~~\text{and}\\
	R_{N}(S_{n})&=\int_{|u|>N}\big(g^{\alpha}(S_{n}+u)-g^{\alpha}(S_{n})\big)u\nu_{\alpha}(du).
	\end{split}
	\end{align*}	
\end{lem}
\noindent
The proof of this lemma follows by similar computations \cite[Lemma 5.3]{k1}.

\begin{lem}\label{l2}
	Let $g^{\alpha}_{h}$ is defined in \eqref{PP1:solution2}. Then for any $N>0$, we have,
	\begin{equation*}
	\mathbb{E}\big[ S_{n}g^{\alpha}(S_{n})  \big]=\sum_{i=1}^{n}\int_{-N}^{N}\mathbb{E}\big[K_{i}(t,N)(g^{\alpha
	}_{h})^{\prime}(S_{n}(i)+t)   \big]dt+R_{1},
	\end{equation*}
	
	where
	\begin{align*}
	K_{i}(t,N)&=\mathbb{E}\big[Z_{i}1_{\{0\leq t\leq Z_{i}\leq N\}} -Z_{i}1_{\{-N\leq Z_{i}\leq t\leq0\}}  \big] , \text{and} \\
	R_{1}&=\sum_{i=1}^{n}\mathbb{E}\big[\xi_{i}\{g^{\alpha
	}_{h}(S_{n})-g^{\alpha}(S_{n}(i)))  \}   \big]1_{\{|\xi_{i}|\geq N\}}  .
	\end{align*}	
\end{lem} 
\noindent
The proof of this lemma follows by similar computations \cite[Lemma 4.5]{k14}. \vspace{2mm}

\noindent
Next, we derive a result using the above two lemmas which is as follows.
\begin{lem}\label{l3}
	Let $g^{\alpha}_{h}$ is defined in \eqref{PP1:solution2}. Then,
	\begin{equation*}
	\begin{split}
	\mathbb{E}\left[\displaystyle\int_{\mathbb{R}}(g^{\alpha
	}_{h}(S_{n}+u)-g^{\alpha}_{h}(S_{n}))u\nu_{\alpha}(du)- S_{n}g^{\alpha
	}_{h}(S_{n})\right]&=\sum_{i=1}^{n}\int_{-N}^{N}\mathbb{E}\big(\frac{K_{\nu_{\alpha}}(t,N)}{n}-K_{i}(t,N) \big)(g^{\alpha
	}_{h})^{\prime}(S_{n}(i)+t)dt\\
	&	+\frac{1}{n}\sum_{i=1}^{n}\mathbb{E}(R_{N}(S_{n}(i)))+R_{1}+R_{2},
	\end{split}
	\end{equation*}
	where $R_{N}(x)$ and $R_{1}$ are defined in Lemmas \ref{l1} and \ref{l2} respectively,
	\[
	R_{2}=\frac{1}{n}\sum_{i=1}^{n}\mathbb{E}\left[\int_{\mathbb{R}}(g^{\alpha
	}_{h}(S_{n}+u) -g^{\alpha
	}_{h}(S_{n}))u\nu_{\alpha}(du)
	-\int_{\mathbb{R}}(g^{\alpha
	}_{h}(S_{n}(i)+u) -g^{\alpha
	}_{h}(S_{n}(i)))u\nu_{\alpha}(du)\right] .
	\]
\end{lem}

\noindent
$\mathbf{Proof.}$ We have,
\begin{align*}
\begin{split}
\mathbb{E}\left[\displaystyle\int_{\mathbb{R}}(g^{\alpha
}_{h}(S_{n}+u)-g^{\alpha}_{h}(S_{n}))u\nu_{\alpha}(du)- S_{n}g^{\alpha
}_{h}(S_{n})\right]
&=\frac{1}{n}\sum_{i=1}^{n}\mathbb{E}\big[\int_{\mathbb{R}}(g^{\alpha
}_{h}(S_{n}(i)+u)-g^{\alpha
}_{h}(S_{n}(i)))u\nu_{\alpha}(du)\\
&- S_{n}g^{\alpha}_{h}(S_{n})\big] +R_{1}+R_{2}+\frac{1}{n}\sum_{i=1}^{n}\mathbb{E}\big[R_{N}(S_{n}(i))\big]\\
&=\sum_{i=1}^{n}\int_{-N}^{N}\mathbb{E}\big(\frac{K_{\nu_{\alpha}}(t,N)}{n}-K_{i}(t,N) \big)(g^{\alpha
}_{h})^{\prime}(S_{n}(i)+t)dt\\
&+\frac{1}{n}\sum_{i=1}^{n}\mathbb{E}(R_{N}(S_{n}(i)))+R_{1}+R_{2},
\end{split}
\end{align*}
\noindent
the desired conclusion follows.

\paragraph{Proof of Theorem \ref{PP1:th6}} By (\ref{PP1:a15}), we have
\begin{align*}
\begin{split}
\mathbb{E}[h(S_{n})-h(X)]&=\mathbb{E}\left(- S_{n}g^{\alpha
}_{h}(S_{n})+ \int_{\mathbb{{R}}}(g^{\alpha
}_{h}(S_{n}+u)-g^{\alpha}_{h}(S_{n}))u\nu_{\alpha}(du)\right)+\mathbb{E}\left(\beta+\int_{|u|>1}u\nu_{\alpha}(du)   \right)g^{\alpha}_{h}(S_{n}).\\
\end{split}
\end{align*}

\noindent
To bound $\mathbb{E}[h(S_{n})-h(X)]$, it is sufficient to find an upper bound of right hand side of the above result. By Lemma \ref{l3} and \eqref{PP1:a55}, we have

\[
\begin{split}
\left|\sum_{i=1}^{n}\int_{-N}^{N}\mathbb{E}\big(\frac{K_{\nu_{\alpha}}(t,N)}{n}-K_{i}(t,N) \big)(g^{\alpha}_{h})^{\prime}(S_{n}(i)+t)dt\right|&\leq\frac{1}{2}||h^{\prime\prime}||  \sum_{i=1}^{n}\int_{-N}^{N}\left|\frac{K_{\nu_{\alpha}}(t,N)}{n}-K_{i}(t,N)\right|  . \\
\end{split}
\]

\noindent
Note that,
\begin{align*}
\left|\frac{1}{n}\sum_{i=1}^{n}\mathbb{E}(R_{N}(S_{n}(i)))\right|
&\leq\frac{1}{n}\sum_{i=1}^{n}\mathbb{E}\int_{|u|>N}\left|f^{\prime}(S_{n}(i)+u)-f^{\prime}(S_{n}(i))\right|u\nu_{\alpha}(du)\\
&\leq2||h^{\prime}||\int_{|u|>N}|u|\nu_{\alpha}(du), ~~\text{and}
\end{align*}
\begin{align*}
|R_{1}|&=\left| \sum_{i=1}^{n}\mathbb{E}\big[Z_{i}\{f^{\prime}(S_{n})-f^{\prime}(S_{n}(i)))  \}   \big]1_{\{|Z_{i}|\geq N\}}\right|\\
&\leq 2||h^{\prime}||\sum_{i=1}^{n}\mathbb{E}\big[ |Z_{i}|1_{ \{|Z_{i}|>N \}}  \big].
\end{align*}

\noindent
Using \eqref{PP1:a56}, we have
\begin{align*}
|R_{2}|&\leq\frac{1}{n}\sum_{i=1}^{n}\left| \mathbb{E}\left[\int_{\mathbb{R}}(g^{\alpha
}_{h}(S_{n}+u) -g^{\alpha
}_{h}(S_{n}))u\nu_{\alpha}(du)
-\int_{\mathbb{R}}(g^{\alpha
}_{h}(S_{n}(i)+u) -g^{\alpha
}_{h}(S_{n}(i)))u\nu_{\alpha}(du)\right] \right| \\
&\leq\frac{1}{n}\sum_{i=1}^{n}\mathbb{E}\left|\displaystyle\int_{\mathbb{R}}\left[(g^{\alpha
}_{h}(S_{n}+u)-g^{\alpha
}_{h}(S_{n}(i)+u))-(g^{\alpha
}_{h}(S_{n})-g^{\alpha
}_{h}(S_{n}(i))) \right]u\nu_{\alpha}(du)\right|\\
&\leq C_{\alpha,m} \frac{1}{n}\sum_{i=1}^{n}\mathbb{E}|Z_i|^{2-\alpha}.
\end{align*}

%

\noindent
Also, for $m_1=m_{2}=m$ and $\beta=0$, we have
\begin{align*}
\left|\mathbb{E}\left(\beta+\int_{|u|>1}u\nu_{\alpha}(du)   \right)g^{\alpha}_{h}(S_{n})\right|=0.
\end{align*}
\noindent
Combining all the estimates above, we get the inequality of the theorem, as desired.

\appendix
\section{Appendix}\label{appendix}

\noindent
In this section, we prove some technical results used in the previous sections.

\begin{pro}
	Let $X\sim \mathcal{S}(\alpha,\beta,m_1,m_2)$. Then, its characteristic exponent $\eta_{\alpha}$ given in \eqref{u3} can be written in the following form.
	\begin{align*}
	\eta_{\alpha}(z)=
	\begin{cases}
	iz \gamma_{\alpha}-d_{\alpha}|z|^{\alpha}\left(1-i\theta\frac{z}{|z|}\tan\frac{\pi}{2}\alpha\right),~~\alpha\in(0,2)\setminus \{1\},\\
	iz\gamma_{1}-d_{1}|z|(1+i\theta\frac{z}{|z|}\frac{2}{\pi}\log|z|),~~\alpha=1,
	\end{cases}
	\end{align*}
	where $\alpha\in(0,2),$ $\gamma_{\alpha}\in\mathbb{R}$, $d_\alpha\geq 0$ and $\theta\in [-1,1] .$
\end{pro}
\noindent
\textbf{Proof.} Recall that for  $X\sim \mathcal{S}(\alpha,\beta,m_1,m_2)$, the characteristic exponent is given by
\begin{align*}
\eta_{\alpha}(z)= iz\beta +\int_{\mathbb{R}}(e^{izu}-1-izu\mathbf{1}_{\{|u|\leq1\}}(u))\nu_{\alpha}(du),~~z\in \mathbb{R},
\end{align*}
\noindent
where $\nu_\alpha$ is the L\'evy measure given by
$$\nu_{\alpha}(du)=\left(m_1\frac{1}{u^{1+\alpha}}\mathbf{1}_{(0,\infty)}(u)+m_2\frac{1}{|u|^{1+\alpha}}\mathbf{1}_{(-\infty,0)} (u)  \right)du.$$

\noindent
Now, we have to consider three different cases to proceed to the derivations of these expressions.

\noindent
\textbf{(i)} $\alpha \in (0,1)$

\noindent
As noted in Section 2, for $\alpha\in (0,1),$ the integral $\int_{\{|u|\leq 1\}}u\nu_{\alpha}(du)<\infty.$ Indeed $\int_{\{|u|\leq 1\}}u\nu_{\alpha}(du)=\frac{m_1-m_2}{1-\alpha}$  Denote $\beta_{1}=\beta-\frac{m_1-m_2}{1-\alpha}.$ So, one can write $\eta_{\alpha}$ as
\begin{align}\label{appen1}
\eta_{\alpha}(z)= iz\beta_{1} +\int_{\mathbb{R}}(e^{izu}-1)\nu_{\alpha}(du),~~z\in \mathbb{R}
\end{align}
\noindent
Suppose $z>0,$ then from \eqref{appen1}
\begin{align}
\nonumber \eta_{\alpha}(z)&=iz\beta_{1} +\int_{0}^{\infty}(e^{izu}-1)\nu_{\alpha}(du)+\int_{-\infty}^{0}(e^{izu}-1)\nu_{\alpha}(du)\\
\nonumber &=iz\beta_{1}+m_1\int_{0}^{\infty}(e^{izu}-1)\frac{du}{u^{1+\alpha}}+m_2\int_{-\infty}^{0}(e^{izu}-1)\frac{du}{|u|^{1+\alpha}}\\
&=iz\beta_{1}+z^{\alpha}\left(m_1\int_{0}^{\infty}(e^{iv}-1)\frac{du}{v^{1+\alpha}}+m_2\int_{0}^{\infty}(e^{-iv}-1)\frac{dv}{v^{1+\alpha}}  \right) \label{appen2}
\end{align}
\noindent
Applying Cauchy's Theorem of contour integration on \eqref{appen2}, we have
\begin{align}
\nonumber \eta_{\alpha}(z)&=iz\beta_{1}+z^{\alpha}\left(m_1e^{-i\frac{\pi}{2}\alpha}L(\alpha)+m_2 e^{i\frac{\pi}{2}\alpha}L(\alpha) \right), 
\end{align}
\noindent
where $L(\alpha)=\int_{0}^{\infty}(e^{-y}-1)\frac{dy}{y^{1+\alpha}}<0$, see \cite[p.164]{k17}.

\noindent
Thus,
\begin{align*}
\nonumber\eta_{\alpha}(z)&=iz\beta_{1}+z^{\alpha}L(\alpha)\left((m_1+m_2)\cos(\frac{\pi}{2}\alpha)+i(m_2-m_1)\sin(\frac{\pi}{2}\alpha)  \right)\\
&=iz\beta_{1}+z^{\alpha}L(\alpha)(m_1+m_2)\cos(\frac{\pi}{2}\alpha)\left(1+ i\frac{m_2-m_1}{m_1+m_2}\tan(\frac{\pi}{2}\alpha)  \right)
\end{align*}

\noindent
For $z<0$, $$\eta_{\alpha}(z)=\overline{\eta_{\alpha}(-z)}=iz\beta_{1}+(-z)^{\alpha}L(\alpha)(m_1+m_2)\cos(\frac{\pi}{2}\alpha)\left(1+ i\frac{m_2-m_1}{m_1+m_2}\frac{z}{|z|}\tan(\frac{\pi}{2}\alpha)  \right).$$
\noindent
Therefore, for any $z\in \mathbb{R}$
\begin{align*}
 \eta_{\alpha}(z)&=iz\beta_{1}+|z|^{\alpha}L(\alpha)(m_1+m_2)\cos(\frac{\pi}{2}\alpha)\left(1+ i\frac{m_2-m_1}{m_1+m_2}\frac{z}{|z|}\tan(\frac{\pi}{2}\alpha)  \right)\\
&=iz\gamma_{\alpha}-d_\alpha|z|^{\alpha}\left(1-i\theta\frac{z}{|z|}\tan(\frac{\pi}{2}\alpha)   \right),
\end{align*}
\noindent
where $\gamma_{\alpha}=\beta_{1}=\beta-\frac{m_1-m_2}{1-\alpha},$ $d_\alpha=(m_1+m_2)\cos(\frac{\pi}{2}\alpha)\int_{0}^{\infty}(1-e^{-y})\frac{dy}{y^{1+\alpha}}$ and $\theta=\frac{m_1-m_2}{m_1+m_2}.$

\noindent
\textbf{(ii)} $\alpha\in (1,2)$

\noindent
As noted in Section 2, for $\alpha\in (1,2),$ the integral $\int_{\{|u|> 1\}}u\nu_{\alpha}(du)<\infty.$ Indeed $\int_{\{|u|> 1\}}u\nu_{\alpha}(du)=\frac{m_1-m_2}{1-\alpha}$  Denote $\beta_{2}=\beta-\frac{m_1-m_2}{1-\alpha}.$ So, one can write $\eta_{\alpha}$ as
\begin{align}\label{appen3}
\eta_{\alpha}(z)= iz\beta_{2} +\int_{\mathbb{R}}(e^{izu}-1-izu)\nu_{\alpha}(du),~~z\in \mathbb{R}
\end{align}
\noindent
Suppose $z>0,$ then from \eqref{appen3}
\begin{align}
\nonumber \eta_{\alpha}(z)&=iz\beta_{2} +\int_{0}^{\infty}(e^{izu}-1-izu)\nu_{\alpha}(du)+\int_{-\infty}^{0}(e^{izu}-1-izu)\nu_{\alpha}(du)\\
\nonumber &=iz\beta_{2}+m_1\int_{0}^{\infty}(e^{izu}-1-izu)\frac{du}{u^{1+\alpha}}+m_2\int_{-\infty}^{0}(e^{izu}-1-izu)\frac{du}{|u|^{1+\alpha}}\\
&=iz\beta_{2}+z^{\alpha}\left(m_1\int_{0}^{\infty}(e^{iv}-1-iv)\frac{du}{v^{1+\alpha}}+m_2\int_{0}^{\infty}(e^{-iv}-1+iv)\frac{dv}{v^{1+\alpha}}  \right) \label{appen4}
\end{align}
\noindent
Applying Cauchy's Theorem of contour integration on \eqref{appen4}, we have
\begin{align}
\nonumber \eta_{\alpha}(z)&=iz\beta_{2}+z^{\alpha}\left(m_1e^{-i\frac{\pi}{2}\alpha}M(\alpha)+m_2 e^{i\frac{\pi}{2}\alpha}M(\alpha) \right), 
\end{align}
\noindent
where $M(\alpha)=\int_{0}^{\infty}(e^{-y}-1+y)\frac{dy}{y^{1+\alpha}}>0$, see \cite[p.164]{k17}.

\noindent
Thus, for any $z>0$
\begin{align*}
\nonumber\eta_{\alpha}(z)&=iz\beta_{2}+z^{\alpha}M(\alpha)\left((m_1+m_2)\cos(\frac{\pi}{2}\alpha)+i(m_2-m_1)\sin(\frac{\pi}{2}\alpha)  \right)\\
&=iz\beta_{2}+z^{\alpha}M(\alpha)(m_1+m_2)\cos(\frac{\pi}{2}\alpha)\left(1+ i\frac{m_2-m_1}{m_1+m_2}\tan(\frac{\pi}{2}\alpha)  \right)
\end{align*}

\noindent
For $z<0$, $$\eta_{\alpha}(z)=\overline{\eta_{\alpha}(-z)}=iz\beta_{2}+(-z)^{\alpha}M(\alpha)(m_1+m_2)\cos(\frac{\pi}{2}\alpha)\left(1+ i\frac{m_2-m_1}{m_1+m_2}\frac{z}{|z|}\tan(\frac{\pi}{2}\alpha)  \right).$$
\noindent
Therefore, for any $z\in \mathbb{R}$
\begin{align*}
 \eta_{\alpha}(z)&=iz\beta_{2}+|z|^{\alpha}M(\alpha)(m_1+m_2)\cos(\frac{\pi}{2}\alpha)\left(1+ i\frac{m_2-m_1}{m_1+m_2}\frac{z}{|z|}\tan(\frac{\pi}{2}\alpha)  \right)\\
&=iz\gamma_{\alpha}-d_\alpha|z|^{\alpha}\left(1-i\theta\frac{z}{|z|}\tan(\frac{\pi}{2}\alpha)   \right),
\end{align*}
\noindent
where $\gamma_{\alpha}=\beta_{1}=\beta-\frac{m_1-m_2}{1-\alpha},$ $d_\alpha=(m_1+m_2)\cos(\frac{\pi}{2}\alpha)\int_{0}^{\infty}(1-e^{-y}-y)\frac{dy}{y^{1+\alpha}}$ and $\theta=\frac{m_1-m_2}{m_1+m_2}.$

\noindent
\textbf{(iii)} $\alpha=1$

\noindent
For $z\in\mathbb{R},$ it is easy to show that

\begin{align*}
\int_{0}^{\infty}\frac{\cos zu -1}{u^{2}}du=-\frac{\pi}{2}z
\end{align*}
\noindent
Now, suppose $z>0,$ then
\begin{align}\label{appen5}
\eta_{1}(z)=iz\beta+\int_{0}^{\infty}(e^{izu}-1-izu\mathbf{1}_{\{|u|\leq 1\}})\nu_{1}(du)+\int_{-\infty}^{0}(e^{izu}-1-izu\mathbf{1}_{\{|u|\leq 1\}})\nu_{1}(du)
\end{align}
\noindent
Let us consider second integral of \eqref{appen5}. Then, we have 
\begin{align}
\nonumber \int_{0}^{\infty}(e^{izu}-1-izu\mathbf{1}_{\{|u|\leq 1\}})\nu_{1}(du)&=m_1\left(\int_{0}^{\infty}\frac{\cos zu -1}{u^{2}}du+i\int_{0}^{\infty}(\sin zu - zu\mathbf{1}_{\{|u|\leq 1 \}} )\frac{du}{u^{2}} \right)\\
&=m_1\left(-\frac{\pi}{2}z+i\lim_{\epsilon \to 0^{+}}\int_{\epsilon}^{\infty}\left(\frac{\sin zu}{u^{2}}-z\frac{u\mathbf{1}_{\{|u|\leq 1 \}}}{u^{2}}   \right)du  \right)\label{appen6}
\end{align}
\noindent
Using the transformation $zu=v$ and changing suitably the limit of integration on \eqref{appen6}, we have

\begin{align}
\nonumber \int_{0}^{\infty}(e^{izu}-1-izu\mathbf{1}_{\{|u|\leq 1\}})\nu_{1}(du)&=m_1\left(-\frac{\pi}{2}z+\lim_{\epsilon \to 0^{+}}\left( -z\int_{\epsilon}^{\epsilon z}\frac{\sin v}{v^{2}}dv+z\int_{\epsilon}^{\infty}\left(\frac{\sin v}{v^{2}}-\frac{\mathbf{1}_{\{|v|\leq1 \}}}{v}   \right)dv  \right)   \right)\\
&=m_1\left( -\frac{\pi}{2}z-iz\log z+iz\int_{0}^{\infty}\left(\frac{\sin v}{v^{2}}-\frac{\mathbf{1}_{\{|v|\leq1 \}}}{v}   \right)dv  \right)  \label{appen7}
\end{align}
The last equality of \eqref{appen7} follows, since $\lim_{\epsilon \to 0^{+}} \int_{\epsilon}^{\epsilon z}\frac{\sin v}{v^{2}}dv=\lim_{\epsilon \to 0^{+}}\int_{\epsilon}^{\epsilon z}\frac{1}{v}dv=\log z.$ If we set $\Gamma=\int_{0}^{\infty}\left(\frac{\sin v}{v^{2}}-\frac{\mathbf{1}_{\{|v|\leq1 \}}}{v}   \right)dv, $ then \eqref{appen7} simplifies to

\begin{align}
\nonumber \int_{0}^{\infty}(e^{izu}-1-izu\mathbf{1}_{\{|u|\leq 1\}})\nu_{1}(du)=m_1\left( -\frac{\pi}{2}z-iz\log z+iz\Gamma \right) 
\end{align}
\noindent
Similarly, the last integral of \eqref{appen5} leads to
$$\int_{-\infty}^{0}(e^{izu}-1-izu\mathbf{1}_{\{|u|\leq 1\}})\nu_{1}(du)=m_2\left(-\frac{\pi}{2}z+iz\log (-z)-iz\Gamma  \right)$$

\noindent
Thus, for any $z>0$

\begin{align*}
\nonumber \eta_{1}(z)&=iz\beta-(m_1+m_2)\frac{\pi}{2}z+i(m_2-m_1)z\log z+iz(m_1-m_2)\Gamma\\
&=iz(\beta +(m_1-m_2)\Gamma)-(m_1+m_2)\frac{\pi}{2}z\left( 1-i \frac{(m_2-m_1)}{m_1+m_2}\frac{2}{\pi} \log z\right)
\end{align*}

\noindent
For any $z<0$,

$$\eta_{1}(z)=\overline{\eta_{1}(-z)}=iz(\beta +(m_1-m_2)\Gamma)-(m_1+m_2)\frac{\pi}{2}(-z)\left( 1-i \frac{(m_2-m_1)}{m_1+m_2}\frac{z}{|z|}\frac{2}{\pi} \log (-z)\right).$$
\noindent
Therefore, for any $z\in \mathbb{R}$
\begin{align*}
\nonumber \eta_{1}(z)&=iz(\beta +(m_1-m_2)\Gamma)-(m_1+m_2)\frac{\pi}{2}|z|\left( 1-i \frac{(m_2-m_1)}{m_1+m_2}\frac{z}{|z|}\frac{2}{\pi} \log |z|\right)\\
&=iz\gamma_{1}-d_1 |z|\left(1+i\theta \frac{z}{|z|}\frac{2}{\pi} \log |z|  \right),
\end{align*}
\noindent
where $\gamma_{1}=\beta +(m_1-m_2)\Gamma,$ $d_1=(m_1+m_2)\frac{\pi}{2}$ and $\theta=\frac{(m_1-m_2)}{m_1+m_2}.$

\noindent
This completes the proof.


\begin{pro}\label{PP1:appendixPro2}
	Let $x,z \in \mathbb{R}$ and $\alpha\in (0,1)$. Then, for all $t\geq 0,$
	\begin{equation}\label{e48}
	\lim_{t\to0^{+}}\frac{1}{t}\left(e^{iz x(e^{-t}-1)}\phi_{t}(z)-1\right)=\left(-x+\beta_{1}+ \displaystyle\int_{\mathbb{R}}ue^{iz u}\nu_{\alpha}(du) \right )(iz),
	\end{equation}
	where $\beta_{1}=\beta- \int_{\{|u|\leq 1 \}}u \nu_{\alpha}(du).$
\end{pro}
\noindent
\textbf{Proof.} Recall from Section 2, if $X$ be a $\alpha$-stable random variable with $\alpha\in (0,1)$ one can write 
\begin{align*}
\phi_{t}(z)=\frac{\phi_{\alpha}(z)}{\phi_{\alpha}(e^{-t}z)}=\exp\left(iz\beta_{1}(1-e^{-t})+\displaystyle\int_{\mathbb{R}}(e^{iz u}-e^{iue^{-t}z})\nu_\alpha(du)\right),~~t\geq 0,
\end{align*}
where $\beta_{1}=\beta- \int_{\{|u|\leq 1 \}}u \nu_{\alpha}(du)$ (see \eqref{PP1:u4}).

\noindent
Now, let us consider LHS of \eqref{e48},
\begin{align}\label{smy1}
\nonumber &	\lim_{t\to0^{+}}\frac{1}{t}\left(e^{iz x(e^{-t}-1)}\phi_{t}(z)-1\right)\\
\nonumber	=&\lim_{t\to0^{+}}\frac{1}{t}\left(\exp\left(izx(e^{-t}-1)+iz\beta_{1}(1-e^{-t})+ \displaystyle\int_{\mathbb{R}}(e^{iz u}-e^{iue^{-t}z })\nu_\alpha(du) \right) -1     \right)\\
=&\lim_{t\to0^{+}}\frac{1}{t}\left(\exp \left(A+iB\right)-1    \right),
\end{align}
\noindent
where 
\begin{align*}
A&=\int_{\mathbb{R}}(\cos(zu) -\cos(zue^{-t}))\nu_\alpha(du) \text{ and}\\
B&=\left(zx(e^{-t}-1)+z\beta_{1}(1-e^{-t})+\int_{\mathbb{R}}(\sin(zu)- \sin(zue^{-t}) ) \nu_\alpha(du)  \right).
\end{align*}

\noindent
Applying Euler's formula for complex exponential to \eqref{smy1}, and rearranging the limits, we have
\begin{align}\label{smy2}
\lim_{t\to0^{+}}\frac{1}{t}\left(e^{iz x(e^{-t}-1)}\phi_{t}(z)-1\right)&=\lim_{t\to0^{+}}\frac{e^{A}\cos(B)-1}{t}+i\lim_{t\to0^{+}}\frac{e^{A}\sin(B)}{t}.
\end{align}

\noindent
It is easy to show that at $t=0$, $e^{A}\cos(B)-1=0$ and $e^{A}\sin(B)=0.$ Thus, on applying L'Hospital rule on \eqref{smy2}, taking limit as $t$ tend to $0^{+}$, and using dominated convergence theorem, we have

\begin{align}
\nonumber \lim_{t\to0^{+}}\frac{1}{t}\left(e^{iz x(e^{-t}-1)}\phi_{t}(z)-1\right)&=\left(\int_{\mathbb{R}}iu\sin(zu)\nu_\alpha(du) -x+\beta_{1}+\int_{\mathbb{R}}u\cos(zu)\nu_{\alpha}(du)  \right)(iz)\\
\nonumber &=\left(-x+\beta_{1}+\int_{\mathbb{R}}u(\cos(zu)+i\sin(zu) )\nu_{\alpha}(du)   \right)(iz)\\
\nonumber &=\left(-x+\beta_{1}+\int_{\mathbb{R}}ue^{izu}\nu_{\alpha}(du)  \right)(iz)
\end{align}

\noindent
This completes the proof.

\begin{pro}\label{PP1:appendixPro3}
	Let $x,z \in \mathbb{R}$ and $\alpha\in (1,2).$ Then, for all $t\geq 0,$
	\begin{equation*}
	\lim_{t\to0^{+}}\frac{1}{t}\left(e^{iz x(e^{-t}-1)}\phi_{t}(z)-1\right)=\left(-x+\beta_{2}+ \displaystyle\int_{\mathbb{R}}u(e^{iz u}-1)\nu_\alpha(du) \right )(iz),
	\end{equation*}
	where $\beta_{2}=\beta+ \int_{\{|u|> 1 \}}u \nu_{\alpha}(du)$ (see \eqref{PP1:u5}).
\end{pro}
\noindent
\textbf{Proof.} Recall from Section 2, if $X$ be a $\alpha$-stable random variable with $\alpha \in (1,2),$ one can write 
\begin{align}\label{sm9}
\phi_{t}(z)=\frac{\phi_{\alpha}(z)}{\phi_{\alpha}(e^{-t}z)}=\exp\left(iz\beta_{2}(1-e^{-t})+\displaystyle\int_{\mathbb{R}}(e^{iz u}-e^{iue^{-t}z }-iuz(1-e^{-t}))\nu_\alpha(du)\right),~~t\geq 0,
\end{align}
where $\beta_{2}=\beta+ \int_{\{|u|> 1 \}}u \nu_{\alpha}(du).$

\noindent
Now, let us consider LHS of \eqref{sm9},
\begin{align}\label{sm10}
\nonumber &	\lim_{t\to0^{+}}\frac{1}{t}\left(e^{iz x(e^{-t}-1)}\phi_{t}(z)-1\right)\\
\nonumber	=&\lim_{t\to0^{+}}\frac{1}{t}\left(\exp\left(izx(e^{-t}-1)+iz\beta_{2}(1-e^{-t})+ \displaystyle\int_{\mathbb{R}}(e^{iz u}-e^{iue^{-t}z}-iuz(1-e^{-t}))\nu_\alpha(du) \right) -1     \right)\\
=&\lim_{t\to0^{+}}\frac{1}{t}\left(\exp \left(C+iD\right)-1    \right),
\end{align}
\noindent
where
\begin{align*}
C&=\int_{\mathbb{R}}(\cos(zu) -\cos(zue^{-t}))\nu_\alpha(du) \text{ and}\\
D&=\left(zx(e^{-t}-1)+z\beta_{2}(1-e^{-t})+\int_{\mathbb{R}}(\sin(zu)- \sin(zue^{-t})-zu(1-e^{-t}) ) \nu_\alpha(du)  \right).
\end{align*}

\noindent
Applying Euler's formula for complex exponential to \eqref{sm10}, and rearranging the limits, we have
\begin{align}\label{sm11}
\lim_{t\to0^{+}}\frac{1}{t}\left(e^{iz x(e^{-t}-1)}\phi_{t}(z)-1\right)&=\lim_{t\to0^{+}}\frac{e^{C}\cos(D)-1}{t}+i\lim_{t\to0^{+}}\frac{e^{C}\sin(D)}{t}.
\end{align}

\noindent
It is easy to show that at $t=0$, $e^{C}\cos(D)-1=0$ and $e^{C}\sin(D)=0.$ Thus, on applying L'Hospital rule on \eqref{sm11}, taking limit as $t$ tend to $0^{+}$, and using dominated convergence theorem, we have

\begin{align}
\nonumber \lim_{t\to0^{+}}\frac{1}{t}\left(e^{iz x(e^{-t}-1)}\phi_{t}(z)-1\right)&=\left(\int_{\mathbb{R}}iu\sin(zu)\nu_\alpha(du) -x+\beta_{2}+\int_{\mathbb{R}}u(\cos(zu)-1 )\nu_{\alpha}(du)  \right)(iz)\\
\nonumber &=\left(-x+\beta_{2}+\int_{\mathbb{R}}u(\cos(zu)+i\sin(zu)-1 )\nu_{\alpha}(du)   \right)(iz)\\
\nonumber &=\left(-x+\beta_{2}+\int_{\mathbb{R}}u(e^{izu}-1)\nu_{\alpha}(du)  \right)(iz)
\end{align}

\noindent
This completes the proof.

\begin{pro}\label{PP1:appendixPro4}
	Let $\alpha\in (0,2)$. Then,
	\begin{align*}
	\frac{1}{\alpha}\displaystyle\int_{{\mathbb R}}ug^{\prime}(x+u)\nu_{\alpha}(du)=\displaystyle\int_{\mathbb{{R}}}\left(g(x+u)-g(x)\right)\nu_{\alpha}(du),~~g\in \mathcal{S}(\mathbb{R}),
	\end{align*}
	where $\nu_{\alpha}$ is the L\'evy measure defined in \eqref{u1}.
\end{pro}
\noindent
\textbf{Proof.} We use Fubini's theorem and change in the order of integration in the following proof.

\noindent
 For $\alpha\in (0,2),$ we have
 
\begin{align*}
\displaystyle\int_{{\mathbb R}}ug^{\prime}(x+u)\nu_{\alpha}(du)&=m_1\displaystyle\int_{0}^{\infty}\frac{ug^{\prime}(x+u)}{u^{1+\alpha}}du+m_2\displaystyle\int_{-\infty}^{0}\frac{ug^{\prime}(x+u)}{(-u)^{1+\alpha}}du\\
&=m_1\displaystyle\int_{0}^{\infty}\frac{g^{\prime}(x+u)}{u^{\alpha}}du-m_2\displaystyle\int_{-\infty}^{0}\frac{g^{\prime}(x+u)}{(-u)^{\alpha}}du\\
&=\alpha m_1 \displaystyle\int_{0}^{\infty}g^{\prime}(x+u)\int_{u}^{\infty}\frac{1}{z^{1+\alpha}}dzdu\\
&-\alpha m_2 \displaystyle\int_{-\infty}^{0}g^{\prime}(x+u)\int_{-\infty}^{u}\frac{1}{(-z)^{1+\alpha}}dzdu\\
&=\alpha m_1\displaystyle\int_{0}^{\infty}\frac{1}{z^{1+\alpha}}\int_{0}^{z}g^{\prime}(x+u)dudz\\
&-\alpha m_2\displaystyle\int_{-\infty}^{0}\frac{1}{(-z)^{1+\alpha}}\int_{z}^{0}g^{\prime}(x+u)dudz\\
&=\alpha\displaystyle\int_{0}^{\infty}\left(g(x+z)-g(x) \right)\frac{m_1}{z^{1+\alpha}}+\alpha\displaystyle\int_{-\infty}^{0}\left(g(x+z)-g(x) \right)\frac{m_2}{(-z)^{1+\alpha}}\\
&=\alpha\displaystyle\int_{{\mathbb R}}\left(g(x+u)-g(x) \right)\nu_{\alpha}(du)
\end{align*}
\begin{pro}\label{PP1:proappendix5}
	Let $\alpha \in (0,1)$ and $h\in\mathcal{H}_\delta$ with $\delta\in (0,\alpha)$. Then,
	\begin{align*}
	\left|xg^{\alpha}_{h}(x)\right| \leq\eta_{\alpha,\beta,\delta,m_1,m_2}(x) :=|\beta_{1}|\|h^{\prime}\|+C_{\alpha,\delta,m_1,m_2}+|x|\wedge |x|^{\delta}+\mathbb{E}|X|^{\delta}.
	\end{align*}
	
\end{pro}
\noindent
\textbf{Proof.} For $\alpha\in (0,1)$, we have by \eqref{PP1:a15},

$$\left|xg^{\alpha}_{h}(x)\right|=\left|-\beta_{1}g^{\alpha}_{h}(x)-\int_{{\mathbb R}}ug^{\alpha}_{h}(x+u)\nu_{\alpha}(du) +(h(x)-h(0))-(\mathbb{E}h(X)-\mathbb{E}h(0)) \right|$$

\noindent
Thus, by \eqref{PP1:a49} and \eqref{PP1:a51}, we have 

\begin{align*}
\left|xg^{\alpha}_{h}(x)\right| \leq |\beta_{1}|\|h^{\prime}\|+C_{\alpha,\delta,m_1,m_2}+|x|\wedge |x|^{\delta}+\mathbb{E}|X|^{\delta}:=\eta_{\alpha,\beta,\delta,m_1,m_2}(x),
\end{align*}
\noindent
the desired conclusion follows.

\section*{Acknowledgments}
\noindent
The authors are thankful to the anonymous referee for some valuable comments that help to improve the article. The first author acknowledges the support of research grant(SB20210848MAMHRD008558) from Ministry of Education through IIT Madras.
The second author gratefully acknowledges the financial support of HTRA fellowship at IIT Madras.

\end{document}